\documentclass[12pt]{article}
\usepackage{amsmath,mathrsfs,amsfonts}
\usepackage{mathtools, bm}
\usepackage{amssymb, bm}
\usepackage{fullpage}
\usepackage{epic}
\usepackage{eepic}
\usepackage[english]{babel}
\usepackage[all]{xy}
\usepackage{enumerate}
\usepackage[T1]{fontenc}
\usepackage[latin1]{inputenc}

\usepackage{soul}
\usepackage{mathabx}
\usepackage{dsfont}
\newtheorem{theo}{Theorem}[subsection]
\newtheorem{conj}[theo]{Conjecture}
\newtheorem{lem}[theo]{Lemma}
\newtheorem{prop}[theo]{Proposition}
\newtheorem{cor}[theo]{Corollary}
\newtheorem{rem}[theo]{Remark}
\def\A{\mathbb{A}}
\def\AF{\mathbb{A}_{\overline{F}}}
\def\lV{\lVert}
\def\rV{\rVert}
\def\lv{\lvert}
\def\rv{\rvert}
\def\Q{\mathbb{Q}}
\def\oF{\overline{F}}
\def\R{\mathbb{R}}
\def\C{\mathbb{C}}
\def\E{\mathbb{E}}
\def\F{\mathbb{F}}
\def\V{\mathbb{V}}
\def\W{\mathbb{W}}
\def\H{\mathbb{H}}
\def\G{\mathbb{G}}
\title{Comparison of local spherical characters and the Ichino-Ikeda conjecture for unitary groups}
\author{Rapha\"{e}l Beuzart-Plessis \protect\footnote{Universit\'e d'Aix-Marseille, I2M-CNRS(UMR 7373), Campus de Luminy, 13288 Marseille C\'edex 9, France rbeuzart@gmail.com}}
\begin{document}

\maketitle

\begin{abstract}
In this paper, we prove a conjecture of Wei Zhang on comparison of certain local spherical characters from which we draw some consequences for the Ichino-Ikeda conjecture for unitary groups. 
\end{abstract}

\tableofcontents

\section{Introduction}

Let $E/F$ be a quadratic extension of number fields. Let $V$ be a $(n+1)$-dimensional hermitian space over $E$ and let $W\subset V$ be a nondegenerate hyperplane. Set $G=U(W)\times U(V)$ and $H=U(W)$. We view $H$ as a subgroup of $G$ via the natural diagonal embedding. Let $\pi$ be a cuspidal automorphic representation of $G(\A)$. Define the $H$-period of $\pi$ to be the linear form $\mathcal{P}_H:\pi\to\C$ given by

$$\displaystyle \mathcal{P}_H(\phi)=\int_{H(F)\backslash H(\A)} \phi(h) dh,\;\;\; \phi\in\pi$$

\noindent where $dh$ stands for the Tamagawa Haar measure on $H(\A)$ (the integral is absolutely convergent by cuspidality of $\pi$). Let $BC(\pi)$ be the base change of $\pi$ to $GL_n(\A_E)\times GL_{n+1}(\A_E)$ (known to exist thanks to the recent work of Mok [Mok] and Kaletha, Minguez, Shin and White [KMSW]). We may decompose $\pi=\pi_n\boxtimes \pi_{n+1}$ with $\pi_n$, $\pi_{n+1}$ cuspidal automorphic representations of $U(W)$ and $U(V)$ respectively. We have a similar decomposition $BC(\pi)=BC(\pi_n)\boxtimes BC(\pi_{n+1})$ with $BC(\pi_n)$, $BC(\pi_{n+1})$ two automorphic representations of $GL_{n,E}$ and $GL_{n+1,E}$ respectively. Let $L(s,BC(\pi))$ denote the $L$-function of pair $L(s,BC(\pi_n)\times BC(\pi_{n+1}))$ defined by Jacquet, Piatetskii-Shapiro and Shalika. If $\pi$ is tempered everywhere (meaning that for all place $v$ the local representation $\pi_v$ is tempered), a famous conjecture of Gan, Gross and Prasad links the nonvanishing of the period $\mathcal{P}_H$ to the nonvanishing of the central value $L(1/2,BC(\pi))$ (see [GGP,conjecture 24.1] for a precise statement). In the influential paper [II], Ichino and Ikeda have proposed a refinement of this conjecture for orthogonal groups in the form of an exact formula relating these two invariants. This conjecture has been suitably extended to unitary groups by N. Harris in his Ph.D. thesis ([Ha]). These formulas are modeled on the celebrated work of Waldspurger ([Wald3]) on toric periods for $GL_2$.

\vspace{2mm}

\noindent In two recent papers ([Zh1], [Zh2]), W.Zhang has proved both the Gan-Gross-Prasad and the Ichino-Ikeda conjectures for unitary groups under some local assumptions on $\pi$. More precisely, Zhang proves the Gan-Gross-Prasad conjecture under some mild local assumptions (mainly that $\pi$ is supercuspidal at one place of $F$ which splits in $E$, see [Zh1,Theorem 1.1]) but he only gets the Ichino-Ikeda conjecture under far more stringent assumptions (see [Zh2,Theorem 1.2]). This discrepancy is due to some local difficulties that we shall discuss shortly. In [Zh2], Zhang makes a series of conjectures (one for every place of $F$) which if true would allow to considerably weaken the assumptions of [Zh2, Theorem 1.2]. The goal of this paper is to prove this conjecture at all nonarchimedean place of $F$. Thus, it will allow us to derive new cases of the Ichino-Ikeda conjecture.

\vspace{2mm}

\noindent Let us now formulate the Ichino-Ikeda conjecture in a form suitable to our purpose. We assume from now on that $\pi$ is everywhere tempered. Set

$$\displaystyle \mathcal{L}(s,\pi):=\Delta_{n+1}\frac{L(s,BC(\pi))}{L(s+\frac{1}{2},\pi,Ad)}$$

\noindent where $\Delta_{n+1}$ is the following product of special values of Hecke $L$-functions

$$\displaystyle \Delta_{n+1}:=\prod_{i=1}^{n+1} L(i,\eta_{E/F}^i)$$

\noindent $\eta_{E/F}$ being the idele class character associated to the extension $E/F$ and where the adjoint $L$-function of $\pi$ is defined by

$$\displaystyle L(s,\pi,Ad):=L(s,BC(\pi_n),As^{(-1)^n})L(s,BC(\pi_{n+1}),As^{(-1)^{n+1}})$$

\noindent (see [GGP, \S 7] for the definition of the Asai L-functions). For all place $v$ of $F$, we will denote by $\mathcal{L}(s,\pi_v)$ the corresponding quotient of local $L$-functions. To the period $\mathcal{P}_H$ we associate a global spherical character $J_\pi$. It is a distribution on the Schwartz space $\mathcal{S}(G(\A))$ of $G(\A)$ given by

$$\displaystyle J_\pi(f)=\sum_{\phi\in \mathcal{B}_\pi} \mathcal{P}_H(\pi(f)\phi)\overline{\mathcal{P}_H(\phi)}$$

\noindent for all $f\in \mathcal{S}(G(\A))$ and where $\mathcal{B}_\pi$ is a (suitable) orthonormal basis of $\pi$ for the Petersson inner product

$$\displaystyle (\phi,\phi')_{Pet}=\int_{G(F)\backslash G(\A)}\phi(g)\overline{\phi(g)}dg$$

\noindent (where $dg$ is the Tamagawa Haar measure on $G(\A)$). We also define local spherical characters as follows. Fix factorizations $dg=\prod_v dg_v$ and $dh=\prod_v dh_v$ of the Tamagawa Haar measures on $G(\A)$ and $H(\A)$ respectively. For all place $v$ of $F$, we define a local spherical character $J_{\pi_v}:\mathcal{S}(G(F_v))\to \C$ (where $\mathcal{S}(G(F_v))$ denotes the Schwartz space of $G(F_v)$) by

$$\displaystyle J_{\pi_v}(f_v)=\int_{H(F_v)} Trace(\pi_v(h)\pi_v(f_v))dh_v,\;\;\; f_v\in\mathcal{S}(G(F_v))$$

\noindent (the integral is absolutely convergent by temperedness of $\pi_v$). For almost all place $v$ of $F$, if $f_v$ is the characteristic function of $G(\mathcal{O}_v)$ we have

$$\displaystyle J_{\pi_v}(f_v)=\mathcal{L}(\frac{1}{2},\pi_v) vol(H(\mathcal{O}_v))vol(G(\mathcal{O}_v))$$

\noindent We define a normalized spherical character $J_{\pi_v}^\natural$ by

$$\displaystyle J_{\pi_v}^\natural(f_v)=\frac{J_{\pi_v}(f_v)}{\mathcal{L}(\frac{1}{2},\pi_v)}$$

\noindent Finally, we will write $S_\pi$ for the component group associated to the $L$-parameter of $\pi$. It is a $2$-abelian group and if $BC(\pi)$ is cuspidal we have $S_\pi\simeq \left( \mathbf{Z}/2\mathbf{Z}\right)^2 $. We can now state the Ichino-Ikeda conjecture as follows:

\begin{conj}[Ichino-Ikeda]\label{II conj}
Assume that $\pi$ is everywhere tempered. Then, for all factorizable test function $f=\prod_v f_v\in \mathcal{S}(G(\A))$ we have

$$\displaystyle J_\pi(f)=\lv S_\pi\rv^{-1}\mathcal{L}(\frac{1}{2},\pi)\prod_v J_{\pi_v}^\natural(f_v)$$
\end{conj}

\noindent Note that the Ichino-Ikeda conjecture is not usually stated this way but rather in a form involving directly the (square of the absolute value of the) period $\mathcal{P}_H$ and some local periods (see [II,conjecture 1.5] and [Ha,conjecture 1.2]), see however [Zh2, lemma 1.7] for the equivalence between the two formulations.

\vspace{2mm}

\noindent The main tool used by Zhang to attack conjecture \ref{II conj} is a comparison of certain (simple) relative trace formulae that have been proposed by Jacquet and Rallis ([JR]). To carry this comparison, we need a fundamental lemma and the existence of smooth matching. The fundamental lemma for the case at hand has been proved by Yun ([Yu]) in positive characteristic and extended by J. Gordon to characteristic $0$ in the appendix to [Yu]. The existence of smooth matching at nonarchimedean places is one of the main achievements of Zhang in [Zh1]. It has been recently extended in a weak form by Xue ([Xue]) to archimedean places. The comparison between the two trace formulae has been done by Zhang in [Zh1]. The output is an identity relating the spherical $J_\pi$ (under some mild local assumptions on $\pi$) to certain periods on the base-change of $\pi$. More precisely, there is a certain spherical character $I_{BC(\pi)}$ attached to these periods and we get an equality between $J_\pi(f)$ and $I_{BC(\pi)}(f')$ up to an explicit factor for {\it nice} matching functions $f$ and $f'$ (see [Zh2,Theorem 4.3] and Theorem \ref{comparison RTF} below). Thanks to the work of Jacquet, Piatetskii-Shapiro and Shalika on Rankin-Selberg convolutions we know an explicit factorization for $I_{BC(\pi)}$ in terms of local (normalized) spherical characters $I^\natural_{BC(\pi_v)}$ (see [Zh2,Proposition 3.6]). As a consequence, we also get an explicit factorization of $J_\pi$. However, this factorization is still in terms of the local spherical characters $I^\natural_{BC(\pi_v)}$ which are living on (products of) general linear groups. In order to get the Ichino-Ikeda conjecture we need to compare them with the our original local spherical characters $J^\natural_{\pi_v}$. It is precisely the content of the following conjecture of Zhang (see [Zh2, conjecture 4.4] and conjecture \ref{Zhang conj} for precise statements):

\begin{conj}[Zhang]\label{conj}
Let $v$ be a place of $F$. Then for all {\it matching} functions $f_v\in \mathcal{S}(G(F_v))$ and $f'_v\in \mathcal{S}(G'(F_v))$ we have

$$\displaystyle I^\natural_{BC(\pi)_v}(f'_v)=C(\pi_v) J^\natural_{\pi_v}(f_v)$$

\noindent where $C(\pi_v)$ is some explicit constant.
\end{conj}

\noindent Together with the above-mentioned comparison of relative trace formulae, this conjecture implies the Ichino-Ikeda conjecture under mild local assumptions (see [Zh2,Proposition 4.5]). Zhang was able to verify his conjecture in certain particular cases. More precisely, in [Zh2] the above conjecture is proved for split places or when the representation $\pi_v$ is unramified (and the residual characteristic is sufficiently large) or supercuspidal (see Theorem 4.6 of {\it loc.cit}). This explains the very strong conditions that are imposed on $\pi$ in [Zh2, Theorem 1.2]. The main purpose of this paper is to prove conjecture \ref{conj} at every nonarchimedean place. Our main result thus reads as follows (see Theorem \ref{theo principal}):

\begin{theo}\label{theo1}
For every nonarchimedean place $v$ of $F$, conjecture \ref{conj} holds at $v$. 
\end{theo}

\noindent As a consequence of this theorem we obtain the following result towards conjecture \ref{II conj} (see Theorem \ref{theo II}):

\begin{theo}\label{theo2}
Let $\pi$ be cuspidal automorphic representation of $G(\A)$ which is everywhere tempered. Assume that all the archimedean places of $F$ split in $E$ and that there exists a nonarchimedean place $v_0$ of $F$ such that $BC(\pi_{v_0})$ is supercuspidal. Then conjecture \ref{II conj} holds for $\pi$.
\end{theo}

\vspace{2mm}

\noindent The main new ingredient in the proof of Theorem \ref{theo1} is a group analog of the local relative trace formula for Lie algebras developed by Zhang in [Zh1,\S 4.1]. Actually, this local trace formula can be derived directly from results contained in [Zh1] and [Beu] so that the proof of it is rather brief (see \S \ref{local TF}). We then deduce Theorem \ref{theo1} from a combination of this local trace formula with certain results of Zhang on truncated local expansion of spherical characters (see [Zh2,\S 8] and \S \ref{truncated local expansion}).

\vspace{2mm}

\noindent We now briefly describe the content of each section. In section 1, we set up the notations, fix the measures and recall a number of results (in particular concerning global and local base-change for unitary groups and the local Gan-Gross-Prasad conjecture) that will be needed in the sequel. In section 2 we mainly recall the work of Zhang on comparison of global relative trace formulae, we state precisely conjecture \ref{conj} as well as the main results (Theorem \ref{theo principal} and Theorem \ref{theo II}). Section 3 is devoted to the proofs of Theorem \ref{theo1} and Theorem \ref{theo2}. In section 4, we explain how we can remove the temperedness assumption in Theorem \ref{theo II}. Finally, we have included an appendix to prove that the simple Jacquet-Rallis trace formulae are still absolutely convergent for test functions which are not necessarily compactly supported (but nevertheless rapidly decreasing). For this, we define certain norms on the automorphic quotient $[G]:=G(F)\backslash G(\mathbb{A})$ and establish their basic properties. This material is certainly classical but the author hasn't be able to find a convenient reference, hence we provide complete proofs. It has however interesting consequences e.g. for $H$ a closed subgroup of $G$ we can give a criterion under which every cuspidal form on $[G]$ is integrable on $[H]$ (see Proposition \ref{prop norms} (ix), the criterion simply being that the variety $H\backslash G$ is quasi-affine).

\vspace{2mm}

\noindent\textbf{Aknowledgement}: I thank Volker Heiermann, Wee-Teck Gan and Hang Xue for useful comments on an earlier draft of this paper. This work has been done while the author was a Senior Research Fellow at the National University of Singapore and the author would like to thank this institution fo its warm hospitality.

\section{Preliminaries}\label{section 1}

\subsection{General notations and conventions}\label{general notations}

In this paper $E/F$ will always be a quadratic extension of number fields or of local fields of characteristic zero. We will always denote by $Tr_{E/F}$ the corresponding trace and by $x\mapsto \overline{x}$ the nontrivial $F$-automorphism of $E$. Moreover, we will fix a nonzero element $\tau\in E$ such that $Tr_{E/F}(\tau)=0$. The notation $R_{E/F}$ will stand for the Weil restriction of scalars from $E$ to $F$. For every finite dimensional hermitian space $V$ over $E$ we will denote by $U(V)$ the corresponding unitary group and we will write $\mathfrak{u}(V)$ for its Lie algebra. The standard maximal unipotent subgroup of $GL_n$ will be denoted by $N_n$. For all connected reductive group $G$ over $F$ we will write $Z_G$ for the center of $G$. For all $n\geqslant 1$ we define a variety $S_n$ over $F$ by

$$\displaystyle S_n:=\{s\in R_{E/F}GL_n; s\overline{s}=1 \}$$

\noindent and its "Lie algebra" $\mathfrak{s}_n$ by

$$\displaystyle \mathfrak{s}_n:=\{X\in R_{E/F}M_n; X+\overline{X}=0 \}$$

\noindent We have a surjective map $\nu:R_{E/F}GL_n/GL_n\to S_n$ given by $\nu(g)=g\overline{g}^{-1}$ which, by Hilbert $90$, is surjective at the level of $k$-points for any field $k$. We will denote by $\mathfrak{c}$ the Cayley map $\mathfrak{c}:X\mapsto (X+1)(X-1)^{-1}$ which realizes a birational isomorphism between $\mathfrak{s}_n$ and $S_n$ and also between $\mathfrak{u}(V)$ and $U(V)$ for all finite dimensional hermitian space $V$ over $E$.

\vspace{2mm}

\noindent Assume that the fields $E$ and $F$ are local. We will then denote by $\lv .\rv_F$ the normalized absolute value on $F$ (and similarly for $E$) and by $\eta_{E/F}$ the quadratic character of $F^\times$ associated to the extension $E/F$. We will also fix an extension $\eta'$ of $\eta_{E/F}$ to $E^\times$ and a nontrivial additive character $\psi:F\to \C^\times$. We will set $\psi_E(z)=\psi(\frac{1}{2}Tr_{E/F}(z))$ for all $z\in E$. Let $G$ be a reductive connected group over $F$. By a representation of $G(F)$ we will always mean a smooth representation if $F$ is $p$-adic and an admissible smooth Fr\'echet representation of moderate growth if $F$ is archimedean (see [BK], [Ca] and [Wall, section 11]). We will denote by $Irr(G)$, $Irr_{unit}(G)$ and $Temp(G)$ the set of isomorphism classes of irreducible, irreducible unitary, irreducible tempered representations of $G(F)$ respectively. We will endow these sets with the Fell topology (see [Tad]). For any parabolic subgroup $P=MU$  of $G$ ($U$ denoting the unipotent radical of $P$ and $M$ a Levi factor) and for any irreducible representation $\sigma$ of $M(F)$ we will denote by $i_P^G(\sigma)$ the normalized parabolic induction of $\sigma$. The notation $\Psi_{unit}(G)$ will stand for the group of unitary unramified characters of $G(F)$. The space of Schwartz functions $\mathcal{S}(G(F))$ consists of locally constant compactly supported functions if $F$ is $p$-adic oand functions rapidly decreasing with all their derivatives if $F$ is archimedean (see [Beu,\S 1.4]). If $F$ is $p$-adic and $\Omega$ is a finite union of Bernstein components of $G(F)$ (see [BD]), we will denote by $\mathcal{S}(G(F))_\Omega$ the corresponding summand of $\mathcal{S}(G(F))$ (for the action by left translation). Finally, if $\pi$ is an irreducible generic representation of $GL_n(E)$ we will denote by $\mathcal{W}(\pi,\psi_E)$ the Whittaker model of $\pi$ with respect to $\psi_E$. It is a space of smooth functions $W:G(F)\to \C$ satisfying the relation

$$\displaystyle W(ug)=\psi_E(\sum_{i=1}^{n-1}u_{i,i+1})W(g)$$

\noindent for all $u\in N_n(E)$ and such that $\pi$ is isomorphic to $\mathcal{W}(\pi,\psi_E)$ equipped with the $G(F)$-action by right translation.

\noindent In the number field case, we will denote by $\A$ and $\A_E$ the adele rings of $F$ and $E$ respectively and by $\eta_{E/F}$ the idele class character associated to the extension $E/F$. We will fix an extension $\eta'$ of $\eta_{E/F}$ to $\A_E^\times$. For every place $v$ of $F$ we will denote by $F_v$ the corresponding completion, $\mathcal{O}_v\subset F_v$ the ring of integers (if $v$ is nonarchimdean) and we will set $E_v=E\otimes_F F_v$, $\mathcal{O}_{E,v}=\mathcal{O}_E\otimes_{\mathcal{O}_F}\mathcal{O}_v$ where $\mathcal{O}_F$, $\mathcal{O}_E$ denote the ring of integers in $F$ and $E$ respectively. If $S$ is a finite set of places of $F$, we define $F_S=\prod_{v\in S}F_v$. If $\Sigma$ is a (usually infinite) set of places of $F$, we will write $\A_\Sigma$ for the restricted product of the $F_v$ for $v\in \Sigma$. We will also fix a nontrivial additive character $\psi:\A/F\to \C^\times$ and we will set $\psi_E(z)=\psi(\frac{1}{2}Tr_{E/F}(z))$ for all $z\in \A_E$. For all place $v$ of $F$, will denote by $\psi_v$, $\psi_{E,v}$ and $\eta'_v$ the local components at $v$ of $\psi$, $\psi_E$ and $\eta'$ respectively. Let $G$ be a connected reductive group over $F$. We will set $[G]=G(F)\backslash G(\A)$ and for all place $v$ of $F$ we will denote by $G_v$ the base-change of $G$ to $F_v$. The Schwartz space $\mathcal{S}(G(\A))$ of $G(\A)$ is by definition the restricted tensor product of the local Schwartz spaces $\mathcal{S}(G(F_v))$. We will denote by $\mathcal{U}(\mathfrak{g}_\infty)$ the enveloping algebra of the complexification of the Lie algebra $\mathfrak{g}_\infty$ of $\prod_{v\mid \infty} G(F_v)$ and by $C_G\in \mathcal{U}(\mathfrak{g}_\infty)$ the Casimir element. If a maximal compact subgroup $K=\prod_v K_v$ of $G(\A)$ has been fixed, we will also denote by $C_K\in \mathcal{U}(\mathfrak{g}_\infty)$ the Casimir element of $K_\infty:=\prod_{v\mid \infty} K_v$. Finally if $\eta:\A^\times/F^\times\to\C^\times$ is an idele class character and $g\in GL_n(\A)$ we will usually abbreviate $\eta(\det g)$ by $\eta(g)$.

\subsection{Analytic families of distributions}\label{analytic families}

\noindent Assume that $F$ is a local field. Let $G$ be a connected reductive group over $F$ and let $\pi\mapsto L_\pi$ be a family of (continuous if $F$ is archimedean) linear forms on $\mathcal{S}(G(F))$ indexed by the set $Temp(G(F))$ of all irreducible tempered representations of $G(F)$. Assume that the following condition is satisfied:

\vspace{4mm}

\hspace{12mm} For all parabolic subgroup $P=MU$ of $G$ and for all square-integrable

\vspace{0.4pt}

\hspace{12mm} representation $\sigma$ of $M(F)$ there is at most one irreducible subrepresentation $\pi$ of

\vspace{0.4pt}

\hspace{12mm} $i_P^G(\sigma)$ such that $L_\pi\neq 0$.

\vspace{4mm}

\noindent This condition is for example automatically satisfied if $G=GL_n$ (as in this case the representation $i_P^G(\sigma)$ is always irreducible). If this condition is satisfied, we may extend the family of distributions $\pi\mapsto L_\pi$ to any induced representation $i_P^G(\sigma)$ as above by setting $L_{i_P^G(\sigma)}=L_\pi$ if $\pi$ is the unique irreducible subrepresentation of $i_P^G(\sigma)$ such that $L_\pi\neq 0$ and $L_{i_P^G(\sigma)}=0$ if no such subrepresentation exists. We then say that this family is {\it analytic} if for all $f\in \mathcal{S}(G(F))$, all parabolic subgroup $P=MU$ and all square-integrable representation $\sigma$ of $M(F)$ the function

$$\displaystyle \chi\in \Psi_{unit}(M)\mapsto L_{i_P^G(\sigma\otimes \chi)}(f)$$

\noindent is analytic (recall that $\Psi_{unit}(M)$ being a compact real torus has a natural structure of analytic variety).

\subsection{Base Change for unitary groups}\label{base change}

Let $E/F$ be a quadratic extension of local fields of characteristic zero (either archimedean or $p$-adic). Let $V$ be a $n$-dimensional hermitian space over $E$. Recall that the set of Langlands parameters for $U(V)$ is in one-to-one correspondence with the set of $(-1)^{n+1}$-conjugate dual continuous semisimple representations $\varphi$ of the Langlands group $\mathcal{L}_E$ of $E$ (see [GGP, \S 3] for a definition of $\epsilon$-conjugate dual representations). In what follows, by a Langlands parameter for $U(V)$ we shall mean a representation $\varphi$ of this sort. By the recent results of Mok [Mok] and Kaletha-Minguez-Shin-White [KMSW] on the local Langlands correspondence for unitary groups together with the work of Langlands [Lan] for real groups, we know that there exists a canonical decomposition

$$\displaystyle Irr(U(V))=\bigsqcup_{\varphi} \Pi^{U(V)}(\varphi)$$

\noindent indexed by the set of all Langlands parameters for $U(V)$. The sets $\Pi^{U(V)}(\varphi)$ are finite (some of them may be empty) and called $L$-packets. By the Langlands classification, the above decomposition boils down to an analog decomposition of the tempered dual

$$\displaystyle Temp(U(V))=\bigsqcup_{\varphi}\Pi^{U(V)}(\varphi)$$

\noindent where the union is over the set of {\it tempered} Langlands parameters for $U(V)$ i.e. the parameters $\varphi$ whose image is bounded. This last decomposition admits a characterization in terms of endoscopic relations (see [Mok, Theorem 3.2.1] and [KMSW, Theorem 1.6.1]) and of the (known) Langlands correspondence for $GL_d(E)$ ([He], [HT], [S]). By this Langlands correspondence, every parameter $\varphi$ of $U(V)$ determines an irreducible representation $\pi(\varphi)$ of $GL_n(E)$. If $\pi$ is in the $L$-packet corresponding to $\varphi$ we will write $BC(\pi):=\pi(\varphi)$. If $\pi$ is tempered then so is $BC(\pi)$ and conversely. However it might happen that $\pi$ is supercuspidal or square-integrable but $BC(\pi)$ is not. Aubert, Moussaoui and Solleveld [AMS] have recently proposed a very general conjecture on how to detect supercuspidal representations in $L$-packets. Moreover, Moussaoui [Mou] has been able to verify this conjecture for orthogonal and symplectic groups. Most probably his work will soon cover unitary groups too. We will need the following particular case of the Aubert-Moussaoui-Solleveld conjecture for which however we can give a direct proof.

\begin{lem}\label{lemma supercuspidal}
Assume that $F$ is $p$-adic. Let $\pi\in Irr(U(V))$ and assume that $BC(\pi)$ is supercuspidal. Then so is $\pi$.
\end{lem} 

\noindent\ul{Proof}: We will use the following characterization of supercuspidal representations

\vspace{4mm}

\hspace{5mm} (1) $\pi$ is supercuspidal if and only if the Harish-Chandra character $\Theta_\pi$ of $\pi$ is compactly

\vspace{0.5pt}

\hspace{12mm} supported modulo conjugation.

\vspace{4mm}

\noindent The necessity is an old result of Deligne ([De]). The sufficiency follows for example from Clozel's formula for the character ([Cl1, Proposition 1]).

\vspace{2mm}

\noindent Let $\varphi$ be the Langlands parameter of $\pi$. Then by our assumption the $L$-packet $\Pi^{U(V)}(\varphi)$ is a singleton. Introduce the {\it twisted} group $\widetilde{GL_n(E)}=GL_n(E)\theta_n$ where $\theta_n(g)={}^t \overline{g}^{-1}$. It is the set of $F$-points of the nonneutral connected component of the non-connected group $G^+=R_{E/F}GL_n\rtimes\{1,\theta_n\}$. Since $\varphi$ is a conjugate-dual representation of $\mathcal{L}_E$, it follows that $BC(\pi)$ may be extended to a representation $BC(\pi)^+$ of $G^+(F)$. Denote by $\widetilde{BC(\pi)}$ the restriction of $BC(\pi)^+$ to $\widetilde{GL_n(E)}$ and denote by $\Theta_{\widetilde{BC(\pi)}}$ the Harish-Chandra character of $\widetilde{BC(\pi)}$ (the Harish-Chandra theory of characters has been extended to twisted groups by Clozel [Cl2]). Since $BC(\pi)$ is supercuspidal, the character $\Theta_{\widetilde{BC(\pi)}}$ is compactly supported modulo conjugation (this follows for example from the equality up to a factor between $\Theta_{\widetilde{BC(\pi)}}$ and weighted orbital integrals of coefficients of $\widetilde{BC(\pi)}$ see [Wald2, th\'eor\`eme 7.1]). By the endoscopic characterization of the local Langlands correspondence for unitary groups, there is a relation between $\Theta_\pi$ and $\Theta_{\widetilde{BC(\pi)}}$. More precisely there is a correspondence between (stable) regular conjugacy classes in $U(V)(F)$ and $\widetilde{GL_n(E)}$ (see [Beu2,\S 3.2], in this particular case the correspondence takes the form of an injective map $U(V)_{reg}(F)/stab\hookrightarrow \widetilde{GL_n(E)}_{reg}/stab$) and for all regular elements $y\in U(V)(F)$, $\widetilde{x}\in \widetilde{GL_n(E)}$ that correspond to each other we have (see [Mok, Theorem 3.2.1] and [KMSW, Theorem 1.6.1])

$$\displaystyle \Theta_{\widetilde{BC(\pi)}}(\widetilde{x})=\Delta(y,\widetilde{x}) \Theta_\pi(y)$$

\noindent where $\Delta(y,\widetilde{x})$ is (up to a sign) a certain transfer factor. From this relation we easily infer that $\Theta_\pi$ is compactly supported modulo conjugation and hence that $\pi$ is supercuspidal by (1). $\blacksquare$

\vspace{2mm}

\noindent We now move on to a global setting. Thus $E/F$ is a quadratic extension of number fields and $V$ is a $n$-dimensional hermitian space over $E$. If $v$ is a place of $F$ which splits in $E$ then we have isomorphisms $U(V)(F_v)\simeq GL_{n}(F_v)$ and $(R_{E/F}GL_n)(F_v)\simeq GL_n(F_v)\times GL_n(F_v)$ and we define a base change map $BC:Irr(U(V)_v)\to Irr((R_{E/F}GL_n)_v)$ by $\pi\mapsto \pi\boxtimes \pi^\vee$. By Theorem 2.5.2 of [Mok] and Theorem 1.7.1/Corollary 3.3.2 of [KMSW] we may associate to any cuspidal automorphic representation $\pi$ of $U(V)$ an isobaric conjugate-dual automorphic representation $BC(\pi)$ of $GL_n(E)$, the base-change of $\pi$, satisfying the following properties:

\begin{enumerate}[(1)]
\item The Asai $L$-function

$$\displaystyle L(s,BC(\pi),As^{(-1)^{n+1}})$$

\noindent has a pole at $s=1$ and moreover if $BC(\pi)$ is cuspidal this pole is simple (see [GGP, \S 7] for the definition of the Asai L-functions);

\item Let $v$ be a place of $F$. Then, if $BC(\pi)$ is generic or $v$ splits in $E$ we have $BC(\pi_v)=BC(\pi)_v$;

\item If $BC(\pi)$ is generic then the multiplicity of $\pi$ in $L^2([U(V)])$ is one (see Theorem 2.5.2/Remark 2.5.3 of [Mok] and Theorem 5.0.5, Theorem 1.7.1 and the discussion thereafter of [KMSW]).
\end{enumerate}

\vspace{2mm}

\noindent Let $v$ be a place of $F$ and $\pi\in Irr(U(V)(F_v))$. Assume first that $v$ is inert in $E$. By the Langlands classification there exist

\begin{itemize}
\item a parabolic subgroup $P=MN$ of $U(V)_{v}$ with

$$\displaystyle M\simeq R_{E_{v}/F_{v}}GL_{n_1}\times\ldots\times R_{E_{v}/F_{v}}GL_{n_r}\times U(V')$$

\noindent where $V'\subset V_{v}$ is a nondegenerate subspace;

\item tempered representations $\pi_i\in Temp(GL_{n_i}(E_{v}))$, $1\leqslant i\leqslant r$, and $\pi'\in Temp(U(V'))$;

\item real numbers $\lambda_1>\ldots>\lambda_r>0$,
\end{itemize}

\noindent such that $\pi$ is the unique irreducible quotient of

$$\displaystyle i_P^{U(V)_{v}}\left(\lv \det\rv_{E_{v}}^{\lambda_1}\pi_1\boxtimes\ldots\boxtimes\lv \det\rv_{E_{v}}^{\lambda_r}\pi_r\boxtimes\pi'\right)$$

\noindent The $r$-uple $(\lambda_1,\ldots,\lambda_r)$ only depends on $\pi$ and we will set $\mathfrak{c}(\pi)=\lambda_1$ if $r\geqslant 1$, $\mathfrak{c}(\pi)=0$ if $r=0$ (i.e. if $\pi$ is tempered). Assume now that $v$ splits in $E$. Then, we have an isomorphism $U(V)_v\simeq GL_{n,F_v}$ and there exists a $r$-uple $(n_1,\ldots,n_r)$ of positive integers such that $n_1+\ldots+n_r=n$, tempered representations $\pi_i\in Temp(GL_{n_1}(F_v))$ $i=1,\ldots,r$ and real numbers $\lambda_1>\lambda_2>\ldots>\lambda_r$ such that $\pi$ is the unique irreducible quotient of

$$\displaystyle i_P^{GL_n}\left(\lv \det\rv_{F_{v}}^{\lambda_1}\pi_1\boxtimes\ldots\boxtimes\lv \det\rv_{F_{v}}^{\lambda_r}\pi_r\right)$$

\noindent where $P$ denotes the standard parabolic subgroup of $GL_n$ with Levi $GL_{n_1}\times\ldots\times GL_{n_r}$. In this case, we set $\mathfrak{c}(\pi)=\max(\lv\lambda_1\rv,\lv \lambda_r\rv)$. This depends only on $\pi$ and in particular not on the choice of the isomorphism $U(V)_v\simeq GL_{n,F_v}$ (which is only defined up to an automorphism of $GL_{n,F_v}$ since it involves the choice of a place of $E$ above $v$).

\vspace{2mm}

\noindent In any case, for $c>0$ we define $Irr_{\leqslant c}(U(V)_{v})$ to be the set of irreducible representations $\pi\in Irr(U(V)_v)$ such that $\mathfrak{c}(\pi)\leqslant c$. Combining the above global results of Mok and Kaletha-Minguez-Shin-White with the bounds toward the Ramanujan conjecture for $GL_n$ of Luo-Rudnick-Sarnak [LRS] suitably extended to ramified places independently by M\"{u}ller-Speh and Bergeron-Clozel ([MS], [BC]), we get the following:

\begin{lem}\label{lemma 2}
Set $c=\frac{1}{2}-\frac{1}{n^2+1}$. Let $\pi$ be a cuspidal automorphic representation of $U(V)(\A)$ such that $BC(\pi)$ is generic. Then, for all place $v$ of $F$ we have

$$\displaystyle \pi_{v}\in Irr_{\leqslant c}(U(V)_{v})$$
\end{lem}

\subsection{The local Gan-Gross-Prasad conjecture for unitary groups}\label{local GGP}

Let $E/F$ be a quadratic extension of local fields of characteristic zero (either archimedean or $p$-adic). Let $W$ be a $n$-dimensional hermitian space over $E$ and define the hermitian space $V$ by $V=W\oplus^\perp e$ where $(e,e)=1$. Set $H=U(W)$ and $G=U(W)\times U(V)$. We view $H$ as a subgroup of $G$ via the diagonal embedding. We will say than an irreducible representation $\pi$ of $G(F)$ is $H$-{\it distinguished} if the space $Hom_H(\pi,\C)$ of $H(F)$-invariant (continuous in the archimedean case) linear forms on $\pi$ is nonzero. By multiplicity one results (see [AGRS], [JSZ]) we always have $\dim Hom_H(\pi,\C)\leqslant 1$. We will denote by $Irr_H(G)$ and $Temp_H(G)$ the subsets of $H$-distinguished representations in $Irr(G)$ and $Temp(G)$ respectively. Let $\varphi$ be a generic Langlands parameter for $G$. We have the following conjecture of Gan, Gross and Prasad ([GGP,conjecture 17.1])

\begin{conj}\label{conjecture GGP}
The $L$-packet $\Pi^G(\varphi)$ contains at most one $H$-distinguished representation.
\end{conj}

\noindent By [Beu, Theorem 12.4.1] and [GI, Proposition 9.3], the following cases of this conjecture are known.

\begin{theo}[Beuzart-Plessis, Gan-Ichino]\label{theorem GGP}.
\begin{enumerate}[(i)]
\item Let $\varphi$ be a tempered Langlands parameter for $G$. Then conjecture \ref{conjecture GGP} holds for $\varphi$.
\item Assume that $F$ is $p$-adic. Then conjecture \ref{conjecture GGP} holds for any generic Lamglands parameter $\varphi$ of $G$.
\end{enumerate}
\end{theo}

\subsection{Measures}\label{measures}

We will use the same normalization of measures as in [Zh2, \S 2]. Let us recall these choices. We actually define two sets of Haar measures: the {\it normalized} and the {\it unnormalized}. We will use the normalized Haar measures apart in section \ref{section 3} where we will use the unnormalized one. From now on and until section \ref{section 3}, where we will switch to a local setting, we fix a quadratic extension $E/F$ of number fields. We will denote by $\eta_{E/F}$ the idele class character corresponding to this extension. We will also fix a nonzero character $\psi:\A/F\to \C^\times$ and a nonzero element $\tau\in E$ such that $Tr_{E/F}(\tau)=0$. We will denote by $\psi_E$ the character of $\A_E$ given by $\psi_E(z)=\psi(\frac{1}{2}Tr_{E/F}(z))$.

\vspace{2mm}

\noindent Let $v$ be a place of $F$. We endow $F_v$ with the self-dual Haar measure for $\psi_v$. Similarly, we endow $E_v$ with the self-dual Haar measure for $\psi_{E,v}$. On $F_v^\times$, we define a normalized measure

$$\displaystyle d^\times x_v=\zeta_{F_v}(1)\frac{dx_v}{\lv x_v\rv_{F_v}}$$

\noindent and an unnormalized one

$$\displaystyle d^* x_v=\frac{dx_v}{\lv x_v\rv_{F_v}}$$

\noindent More generally, for all $n\geqslant 1$, we equip $GL_n(F_v)$ with the following normalized Haar measure

$$\displaystyle dg_v=\zeta_{F_v}(1)\frac{\prod_{ij}dg_{v,ij}}{\lv \det g_v\rv_{F_v}^n}$$

\noindent as well as with the following unnormalized one

$$\displaystyle d^*g_v=\frac{\prod_{ij}dg_{v,ij}}{\lv \det g_v\rv_{F_v}^n}$$

\noindent and similarly for $GL_n(E_v)$. Recall that $N_n$ denotes the standard maximal unipotent subgroup of $GL_n$. We will give $N_n(F_v)$ and $N_n(E_v)$ the Haar measures

$$\displaystyle du_v=\prod_{1\leqslant i<j\leqslant n} du_{v,ij}$$

\noindent We equip $\A^\times$, $N_n(\A)$, $N_n(\A_E)$, $GL_n(\A)$ and $GL_n(\A_E)$ with the global Tamagawa Haar measures given by

$$\displaystyle d^\times x=\prod_v d^\times x_v,\;\; du=\prod_v du_v,\;\; dg=\prod_v dg_v$$

\noindent Recall that $S_n=\{s\in R_{E/F}GL_n; s\overline{s}=1\}$ and its {\it Lie algebra} $\mathfrak{s}_n=\{X\in R_{E/F}M_n; X+\overline{X}=0 \}$. Let $V$ be a $n$-dimensional hermitian space over $E$ and denote by $\mathfrak{u}(V)$ the Lie algebra of $U(V)$. Choosing a basis of $V$ we get an embedding $\mathfrak{u}(V)\hookrightarrow R_{E/F}M_n$. Let us denote by $\langle .,\rangle$ the $GL_n(E_v)$-invariant bilinear pairing on $M_n(E_v)$ given by

$$\displaystyle \langle X,Y\rangle:=Trace(XY)$$

\noindent Note that the restrictions of $\langle .,.\rangle$ to $\mathfrak{s}_n(F_v)$ and $\mathfrak{u}(V)(F_v)$ are $F_v$-valued and nondegenerate. We define a Haar measure $dX$ on $\mathfrak{u}(V)(F_v)$ such that the Fourier transform

$$\displaystyle \widehat{\varphi}(Y)=\int_{\mathfrak{u}(V)(F_v)} \varphi(X) \psi_v(\langle X,Y\rangle) dX$$

\noindent and its dual

$$\displaystyle \widecheck{\varphi}(X)=\int_{\mathfrak{u}(V)(F_v)} \varphi(Y) \psi_v(-\langle Y,X\rangle) dX$$

\noindent are inverse of each other. We define similarly a Haar measure and Fourier transforms $\varphi\mapsto \widehat{\varphi}$, $\varphi\mapsto \widecheck{\varphi}$ on $\mathfrak{s}_n(F_v)$.

\vspace{2mm}

\noindent The Cayley map $\mathfrak{c}:X\mapsto \mathfrak{c}(X)=(1+X)(1-X)^{-1}$ induces birational isomorphisms from $\mathfrak{s}_n$ to $S_n$ and from $\mathfrak{u}(V)$ to $U(V)$. We define the {\it unnormalized} Haar measure $d^*g_v$ on $U(V)(F_v)$ to be the unique Haar measure such that the Jacobian of $\mathfrak{c}$ at the origin is $1$. The normalized Haar measure on $U(V)(F_v)$ is defined by

$$\displaystyle dg_v=L(1,\eta_{E_v/F_v})d^*g_v$$

\noindent Similarly, we endow $S_n(F_v)$ with an unnormalized measure $d^* s_v$ which is the unique $GL_n(E_v)$-invariant measure for which the Jacobian of the Cayley map $\mathfrak{c}$ at the origin is $1$. The corresponding normalized measure is given by

$$\displaystyle ds_v=L(1,\eta_{E_v/F_v})d^*s_v$$

\noindent Note that $d^* s_v$ (resp. $d s_v$) can also be identified with the quotient of the unnormalized (resp.normalized) Haar measures on $GL_n(E_v)$ and $GL_n(F_v)$ via the isomorphism \\ \noindent $\nu:GL_n(E_v)/GL_n(F_v)\simeq S_n(F_v)$, $\nu(g)=g\overline{g}^{-1}$. Finally, we equip $U(V)(\A)$ with the global Haar measure given by

$$\displaystyle dg=\prod_v dg_v$$

\noindent It is not the Tamagawa measure since there is a factor $L(1,\eta_{E/F})^{-1}$ missing. Note that the local normalized Haar measure $dg_v$ can be identified with the quotient of the normalized Haar measures on $E_v^\times$ and $F_v^\times$ via the isomorphism $E_v^\times/F_v^\times\simeq U(1)(F_v)$, $x\mapsto x/\overline{x}$. Hence, as the Tamagawa number of $U(1)$ is $2$, we have

$$\displaystyle vol\left(E^\times \A^\times\backslash \A_E^\times \right)=vol([U(1)])=2L(1,\eta_{E/F}) \leqno (1)$$

\section{Spherical characters, the Ichino-Ikeda conjecture and Zhang's conjecture}\label{section 2}

In this section $E/F$ will be a quadratic extension of number fields and we will use normalized Haar measures (see \S \ref{measures}). Let $W$ be an hermitian space of dimension $n$ over $E$. We will set $V=W\oplus^\perp Ee$ where $(e,e)=1$, $G=U(W)\times U(V)$ and $H=U(W)$. We view $H$ as a subgroup of $G$ via the diagonal embedding. We will fix a maximal compact subgroup $K=\prod_v K_v$ of $G(\A)$. We will say that an irreducible representation $\pi=\otimes'_v \pi_v$ of $G(\A)$ is {\it abstractly} $H$-{\it distinguished} if for all place $v$ of $F$ the representation $\pi_v$ is $H_v$-distinguished i.e. if $Hom_{H_v}(\pi_v,\C)\neq 0$. Set $G'=R_{E/F}\left(GL_n\times GL_{n+1}\right)$. We define two subgroups $H_1'=R_{E/F} GL_n$ and $H_2'=GL_n\times GL_{n+1}$ of $G'$ ($H_1'$ is embedded diagonally). We also define a character $\eta$ of $H_2'(\A)$ by

$$\displaystyle \eta(g_1,g_2)=\eta_{E/F}(g_1)^{n+1} \eta_{E/F}( g_2)^n$$

\noindent for all $(g_1,g_2)\in H_2'(\A)=GL_n(\A)\times GL_{n+1}(\A)$. We will also fix a maximal compact subgroup $K'=\prod_v K'_v$ of $G'(\A)$ such that $K'_v=GL_n(\mathcal{O}_{E,v})\times GL_{n+1}(\mathcal{O}_{E,v})$ for all nonarchimedean place $v$ of $F$. Finally, if $\pi$ and $\Pi$ are cuspidal automorphic representations of $G(\A)$ and $G'(\A)$ respectively then we endow them with the following Petersson inner products

$$\displaystyle  (\phi_1,\phi_2)_{Pet}:=\int_{[G]}\phi_1(g) \overline{\phi_2(g)} dg,\;\;\; \phi_1,\phi_2\in \pi$$

$$\displaystyle (\phi'_1,\phi'_2)_{Pet}:=\int_{[Z_{G'}\backslash G']}\phi'_1(g') \overline{\phi'_2(g')} dg',\;\;\; \phi'_1,\phi'_2\in \Pi$$

\subsection{Global spherical characters}\label{global spherical characters}

\noindent For any cuspidal automorphic representation $\pi$ of $G(\A)$ we define the $H$-period $\mathcal{P}_H:\pi\to \C$ by

$$\displaystyle \mathcal{P}_H(\phi)=\int_{[H]} \phi(h) dh,\;\;\; \phi\in \pi$$

\noindent The integral is absolutely convergent (see Proposition \ref{prop norms}(ix)). We will say that the cuspidal automorphic representation $\pi$ is {\it globally} $H$-{\it distinguished} if the period $\mathcal{P}_H$ is not identically zero on $\pi$. We may associate to this period a (global) spherical character $J_\pi: \mathcal{S}(G(\A))\to \C$ defined as follows. Let $f\in \mathcal{S}(G(\A))$ and choose a compact-open subgroup $K_f\subset G(\A_f)$ by which $f$ is biinvariant. Let $\mathcal{B}_\pi^{K_f}$ be an orthonormal basis for the Petersson inner product of $\pi^{K_f}$ whose elements are $C_G$ and $C_K$ eigenvectors. Then we set

$$\displaystyle J_\pi(f)=\sum_{\phi\in \mathcal{B}_\pi^{K_f}} \mathcal{P}_H(\pi(f)\phi) \overline{\mathcal{P}_H(\phi)}$$

\noindent The sum is absolutely convergent and does not depend on the choice of the basis $\mathcal{B}_\pi^{K_f}$ (see Proposition \ref{prop period}).

\vspace{2mm}

\noindent Let $\Pi$ be a cuspidal automorphic representation of $G'(\A)$ whose central character is trivial on $Z_{H_2'}(\A)=\A^\times\times\A^\times$. We define two periods $\lambda:\Pi\to \C$ and $\beta:\Pi\to \C$ by

$$\displaystyle \lambda(\phi)=\int_{[H_1']}\phi(h_1) dh_1$$

$$\displaystyle \beta(\phi)=\int_{[Z_{H_2'}\backslash H_2']} \phi(h_2)\eta(h_2)dh_2$$

\noindent for all $\phi\in \Pi$. The two above integrals are absolutely convergent (see Proposition \ref{prop norms}(ix)). We also define a (global) spherical character $I_\Pi:\mathcal{S}(G'(\A))\to \C$ as follows. Let $f'\in \mathcal{S}(G'(\A))$ and choose a compact-open subgroup $K_{f'}\subset G(\A_f)$ by which $f'$ is biinvariant. Let $\mathcal{B}_\Pi^{K_{f'}}$ be an orthonormal basis for the Petersson inner product of $\Pi^{K_{f'}}$ whose elements are $C_{G'}$ and $C_{K'}$ eigenvectors. Then we set

$$\displaystyle I_\Pi(f')=\sum_{\phi\in \mathcal{B}_\Pi^{K_{f'}}} \lambda(\Pi(f')\phi) \overline{\beta(\phi)}$$

\noindent \noindent The sum is absolutely convergent and does not depend on the choice of the basis $\mathcal{B}_\Pi^{K_{f'}}$ (see Proposition \ref{prop period}).

\subsection{Local spherical characters}\label{local spherical characters}

\noindent Let $v$ be a place of $F$. Let $\pi_v$ be a tempered representation of $G(F_v)$. We define a distribution $J_{\pi_v}:\mathcal{S}(G(F_v))\to \C$ (the local spherical character associated to $\pi_v$) by

$$\displaystyle J_{\pi_v}(f_v)=\int_{H(F_v)} Trace(\pi_v(h)\pi_v(f_v)) dh_v,\;\;\; f_v\in \mathcal{C}(G(F_v))$$

\noindent By [Beu, \S 8.2], the above integral is absolutely convergent. Choosing models for $G$ and $H$ over $\mathcal{O}_F$, for almost all $v$ if $f_v=\mathbf{1}_{K_v}$, we have

$$\displaystyle J_{\pi_v}(f_v)=\mathcal{L}(\frac{1}{2},\pi_v) vol(H(\mathcal{O}_v))vol(G(\mathcal{O}_v))$$

\noindent (see the introduction for the definition of $\mathcal{L}(s,\pi_v)$). Hence, we define a {\it normalized} spherical character $J_{\pi_v}^\natural$ by

$$\displaystyle J_{\pi_v}^\natural=\frac{1}{\mathcal{L}(\frac{1}{2},\pi_v)}J_{\pi_v}$$

\noindent By [Beu, Theorem 8.2.1] we have

\vspace{4mm}

\hspace{5mm} (1) $\pi_v$ is $H_v$-distinguished if and only if $J_{\pi_v}\neq 0$.

\vspace{4mm}

\noindent Moreover, by [Beu, Corollary 8.6.1], for all parabolic subgroup $P=MU$ of $G_v$ and for all square-integrable representation $\sigma$ of $M$ there is at most one irreducible subrepresentation $\pi\subset i_P^G(\sigma)$ such that $J_\pi\neq 0$. Thus, we are in the situation of \S \ref{analytic families} and the family of distributions $\pi_v\in Temp(G_v)\mapsto J_{\pi_v}$ is analytic.

\vspace{2mm}

\noindent Let $\Pi_v$ be a generic unitary representation of $G'(F_v)$. We may write $\Pi_v=\Pi_{n,v}\boxtimes \Pi_{n+1,v}$ where $\Pi_{n,v}$ and $\Pi_{n+1,v}$ are generic and unitary representations of $GL_n(E_v)$ and $GL_{n+1}(E_v)$ respectively. Let $\mathcal{W}(\Pi_{n,v},\overline{\psi}_E)$ and $\mathcal{W}(\Pi_{n+1,v},\psi_E)$ be the Whittaker models of $\Pi_{n,v}$ and $\Pi_{n+1,v}$ corresponding to the characters $\overline{\psi}_E$ and $\psi_E$ respectively. Set $\mathcal{W}(\Pi_v)=\mathcal{W}(\Pi_{n,v},\overline{\psi}_E)\otimes \mathcal{W}(\Pi_{n+1,v},\psi_E)$. We define linear forms (the local Flicker-Rallis periods)

$$\beta_{n,v}: \mathcal{W}(\Pi_{n,v},\overline{\psi}_E)\to \C, \;\;\; \beta_{n+1,v}: \mathcal{W}(\Pi_{n+1,v},\psi_E)\to \C$$

\noindent and scalar products

$$\theta_{n,v}:\mathcal{W}(\Pi_{n,v},\overline{\psi}_E)\times \mathcal{W}(\Pi_{n,v},\overline{\psi}_E)\to \C,\;\;\; \theta_{n+1,v}:\mathcal{W}(\Pi_{n+1,v},\psi_E)\times \mathcal{W}(\Pi_{n+1,v},\psi_E)\to \C$$

\noindent by

$$\displaystyle \beta_{k,v}(W_k)=\int_{N_{k-1}(F_v)\backslash GL_{k-1}(F_v)} W_k(\epsilon_k(\tau)g_{k-1}) \eta_{E_v/F_v}(\det g_{k-1})^{k-1}dg_{k-1}$$

$$\displaystyle \theta_{k,v}(W_k, W'_k)=\int_{N_{k-1}(E_v)\backslash GL_{k-1}(E_v)} W_k(g_{k-1})\overline{W'_k(g_{k-1})}dg_{k-1}$$

\noindent for all $k=n,n+1$, all $W_n,W'_n\in \mathcal{W}(\Pi_{n,v},\overline{\psi}_E)$ and all $W_{n+1},W'_{n+1}\in\mathcal{W}(\Pi_{n+1,v},\psi_E)$ where $\epsilon_k(\tau)=diag(\tau^{k-1},\tau^{k-2},\ldots,\tau,1)$ (recall that $\tau$ is a fixed nonzero element of $E$ such that $Tr_{E/F}(\tau)=0$). The above integrals are absolutely convergent (see [JS] Propositions 1.3 and 3.16 for the absolute convergence of $\theta_{k,v}$, the proof for $\beta_{k,v}$ is identical). Set $\beta_v=\beta_{n,v}\otimes \beta_{n+1,v}$ and $\theta_v=\theta_{n,v}\otimes \theta_{n+1,v}$. If $E_v/F_v$, $\Pi_v$, $\psi_{E,v}$ are unramified, $\tau$ is a unit in $E_v$ and $W_v\in \mathcal{W}(\Pi_v)$ is the unique $K'_v$-invariant vector such that $W_v(1)=1$, we have (see [JS, Proposition 2.3] and [Zh2,\S 3.2])

$$\displaystyle \beta_v(W_v)=vol(K'_v)L(1,\Pi_{n,v},As^{(-1)^{n-1}})L(1,\Pi_{n+1,v}, As^{(-1)^n})$$

\noindent and

$$\displaystyle \theta_v(W_v)=vol(K'_v)L(1,\Pi_{n,v}\times \Pi_{n,v}^\vee)L(1,\Pi_{n+1,v}\times \Pi_{n+1,v}^\vee)$$

\noindent Hence, we define normalized versions $\beta_v^\natural$ and $\theta_v^\natural$ of $\beta_v$ and $\theta_v$ by

$$\displaystyle \beta_v^\natural=\frac{\beta_v}{L(1,\Pi_{n,v},As^{(-1)^{n-1}})L(1,\Pi_{n+1,v}, As^{(-1)^n})},\;\;\; \theta_v^\natural=\frac{\theta_v}{L(1,\Pi_{n,v}\times \Pi_{n,v}^\vee)L(1,\Pi_{n+1,v}\times \Pi_{n+1,v}^\vee)}$$

\noindent For all $s\in \C$, we also have the local Rankin-Selberg period $\lambda_v(s,.):\mathcal{W}(\Pi_v)\to \C$ defined by

$$\displaystyle \lambda_v(s,W_n\otimes W_{n+1})=\int_{N_n(E_v)\backslash GL_n(E_v)} W_n(g_n)W_{n+1}(g_n) \lv \det g_n\rv_{E_v}^sdg_n$$

\noindent for all $(W_n,W_{n+1})\in \mathcal{W}(\Pi_{n,v},\overline{\psi}_E)\times\mathcal{W}(\Pi_{n+1,v},\psi_E)$, and its normalization $\lambda_v^\natural(s,.)$ given by

$$\displaystyle \lambda_v^\natural(s,.)=\frac{\lambda_v(s,.)}{L(s+\frac{1}{2},\Pi_{n,v}\times \Pi_{n+1,v})}$$

\noindent The integral defining $\lambda_v(s,.)$ is absolutely convergent for $Re(s)\gg 0$ and $\lambda_v^\natural(s,.)$ extends to an entire function on $\C$ (see [JPSS] and [Jac] for the archimedean case). We will set $\lambda_v^\natural=\lambda_v^\natural(0,.)$. Obviously $\lambda_v^\natural$ defines a $H_1'(F_v)$-invariant linear form on $\Pi_v$. Moreover by [JPSS] and [Jac], there exists $W\in \mathcal{W}(\Pi_v)$ such that $\lambda_v^\natural(W)=1$. Hence $\lambda_v^\natural$ defines a nonzero element in $Hom_{H_1'}(\Pi_v,\C)$. If $\Pi_v$ is tempered then $\lambda_v(s,.)$ is absolutely convergent for $Re(s)>-1/2$ and we will set $\lambda_v=\lambda_v(0,.)$.

\vspace{2mm}

\noindent We are now ready to define the (normalized) local spherical character $I^\natural_{\Pi_v}:\mathcal{S}(G'(F_v))\to \C$ attached to $\Pi_v$. Let $f'_v\in \mathcal{S}(G'(F_v))$. If $v$ is nonarchimedean then choose a compact-open subgroup $K_{f'_v}$ of $G'(F_v)$ by which $f'_v$ is biinvariant and let $\mathcal{B}_{\Pi_v}$ be an orthonormal basis of $\Pi_v^{K_{f'_v}}$ for the scalar product $\theta^\natural_v$. If $v$ is archimedean, we let $\mathcal{B}_{\Pi_v}$ be any orthonormal basis of $\Pi_v$ for the scalar product $\theta^\natural_v$ consisting of $C_{K'_v}$-eigenvectors. Then we set

$$\displaystyle I^\natural_{\Pi_v}(f'_v)=\sum_{W\in \mathcal{B}_{\Pi_v}} \lambda^\natural_v(\Pi_v(f'_v)W)\overline{\beta^\natural_v(W)}$$

\noindent The sum is absolutely convergent and does not depend on the choice of $\mathcal{B}_{\Pi_v}$. If moreover $\Pi_v$ is tempered then we define an {\it unnormalized} local spherical character $I_{\Pi_v}:\mathcal{S}(G'(F_v))\to \C$ by using $\theta_v$, $\beta_v$ and $\lambda_v$ instead of $\theta^\natural_v$, $\beta^\natural_v$ and $\lambda^\natural_v$. Finally, the proofs of [JS, Proposition 1.3] and [JPSS, Theorem 2.7] easily show that the family of distributions $\Pi_v\in Temp(G'_v)\mapsto I_{\Pi_v}$ is analytic in the sense of \S \ref{analytic families}.

\subsection{Orbital integrals}\label{orbital integrals}

Consider the action of $H\times H$ on $G$ by left and right translations. Then, an element $\delta\in G$ is said to be {\it regular semisimple} for this action if its orbit is closed and its stabilizer is trivial. Denote by $G_{rs}$ the open subset of regular semisimple element in $G$. Let $v$ be a place of $F$ and $\delta\in G_{rs}(F_v)$ be regular semisimple. We define the (relative) orbital integral associated to $\delta$ as the distribution given by

$$\displaystyle O(\delta,f_v)=\int_{H(F_v)\times H(F_v)} f_v(h\delta h')dhdh',\;\;\; f_v\in \mathcal{S}(G(F_v))$$

\noindent There is another way to see these orbital integrals. For all $f_v\in \mathcal{S}(G(F_v))$, we define a function $\widetilde{f}_v\in \mathcal{S}(U(V)(F_v))$ by

$$\displaystyle \widetilde{f}_v(x)=\int_{H(F_v)} f_v(h(1,x)) dh,\;\;\; x\in U(V)(F_v)$$

\noindent This defines a surjective linear map $\mathcal{S}(G(F_v))\to \mathcal{S}(U(V)(F_v))$. Let us say that an element $x\in U(V)$ is regular semisimple if it is so for the action of $U(W)$ by conjugation i.e. if the $U(W)$-conjugacy class of $x$ is closed and the stabilizer of $x$ in $U(W)$ is trivial. Denote by $U(V)_{rs}$ the open subset of regular semisimple element in $U(V)$. For all $x\in U(V)_{rs}(F_v)$ we define the orbital integral associated to $x$ as the distribution

$$\displaystyle O(x,\varphi)=\int_{U(W)(F_v)} \varphi(hxh^{-1})dh,\;\;\; \varphi_v\in \mathcal{S}(U(V)(F_v))$$

\noindent For all $\delta=(\delta_W,\delta_V)\in G_{rs}$ the element $x=\delta_W^{-1}\delta_V$ is regular semisimple in $U(V)$ and this defines a surjection $G_{rs}\twoheadrightarrow U(V)_{rs}$. Moreover, for all $\delta\in G_{rs}(F_v)$ and all $f\in \mathcal{S}(G(F_v))$, we have the equality

$$\displaystyle O(\delta,f)=O(x,\widetilde{f})$$

\noindent where $x=\delta_W^{-1}\delta_V$.

\vspace{2mm}

\noindent We can also define orbital integrals on the space $\mathcal{S}(\mathfrak{u}(V)(F_v))$. Call an element $X\in \mathfrak{u}(V)$ regular semisimple if it is so for the adjoint action of $U(W)$. Let us denote by $\mathfrak{u}(V)_{rs}$ the open subset of regular semisimple elements. Then, for all $X\in \mathfrak{u}(V)_{rs}(F_v)$ we can define an orbital integral by

$$\displaystyle O(X,\varphi)=\int_{U(W)(F_v)} \varphi(h^{-1}Xh)dh,\;\;\; \varphi\in \mathcal{S}(U(V)(F_v))$$

\noindent The Cayley map $\mathfrak{c}:X\mapsto (1+X)(1-X)^{-1}$ realizes a $U(W)$-equivariant isomorphism between the open subsets $\mathfrak{u}(V)^\circ=\{X\in \mathfrak{u}(V);\; \det(1-X)\neq 0\}$ and $U(V)^0=\{x\in U(V);\; \det(1+x)\neq 0\}$. Assume that $v$ is nonarchimedean and let $\omega\subset \mathfrak{u}(V)^\circ(F_v)$ and $\Omega\subset U(V)^\circ(F_v)$ be open and closed $U(W)(F_v)$-invariant neighborhoods of $0$ and $1$ respectively such that the Cayley map restricts to an analytic isomorphism between $\omega$ and $\Omega$ preserving measures. For all $\varphi\in \mathcal{S}(U(V)(F_v))$, we define a function $\varphi_\natural$ by

$$\displaystyle \varphi_{\natural}(X)=\left\{
    \begin{array}{ll}
        \varphi(\mathfrak{c}(X)) & \mbox{ if } X\in\omega \\
        0 & \mbox{ otherwise}
    \end{array}
\right.
$$

\noindent Then for all $X\in \omega_{rs}=\omega\cap \mathfrak{u}(V)_{rs}(F_v)$ and all $\varphi\in \mathcal{S}(U(V)(F_v))$ we have

$$\displaystyle O(\mathfrak{c}(X),\varphi)=O(X,\varphi_\natural)$$

Consider now the action of $H_1'\times H_2'$ on $G'$ by left and right translations. As before, an element $\gamma\in G'$ is said to be {\it regular semisimple} for this action if its orbit is closed and its stabilizer is trivial. Denote by $G'_{rs}$ the open subset of regular semisimple element in $G'$. Let $v$ be a place of $F$ and $\gamma\in G'_{rs}(F_v)$ be regular semisimple. We define the (relative) orbital integral associated to $\gamma$ as the distribution given by

$$\displaystyle O(\gamma,f'_v)=\int_{H_1'(F_v)\times H_2'(F_v)} f'_v(h_1^{-1}\gamma h_2)\eta(h_2)dh_1dh_2,\;\;\; f'_v\in \mathcal{S}(G'(F_v))$$

\noindent There is another way to see these orbital integrals. Recall that $S_{n+1}(F_v)=\{s\in GL_{n+1}(E_v); s\overline{s}=1\}$ and that we have a surjective map $\nu:GL_{n+1}(E_v)\to S_{n+1}(F_v)$, $\nu(g)=g\overline{g}^{-1}$. For all $f'_v\in \mathcal{S}(G'(F_v))$, we define a function $\widetilde{f}'_v\in \mathcal{S}(S_{n+1}(F_v))$ by

$$\displaystyle \widetilde{f}'_v(s)=\int_{H_1'(F_v)}\int_{GL_{n+1}(F_v)} f'_v(h_1(1,gh_2)) dh_2dh_1,\;\;\; g\in GL_{n+1}(E_v), \; s=\nu(g)$$

\noindent if $n$ is even and

$$\displaystyle \widetilde{f}'_v(s)=\int_{H_1'(F_v)}\int_{GL_{n+1}(F_v)} f'_v(h_1(1,gh_2)) \eta_v'(gh_2)dh_2dh_1,\;\;\; g\in GL_{n+1}(E_v), \; s=\nu(g)$$

\noindent if $n$ is odd. In any case, this defines a surjective linear map $\mathcal{S}(G'(F_v))\to \mathcal{S}(S_{n+1}(F_v))$. The group $GL_n$ acts on $S_{n+1}$ by conjugation and we shall say that an element $s\in S_{n+1}$ is regular semisimple if it is so for this action i.e. if the $GL_n$-conjugacy class of $s$ is closed and the stabilizer of $s$ in $GL_n$ is trivial. We will denote by $S_{n+1,rs}$ the open subset of regular semisimple elements in $S_{n+1}$. For all $s\in S_{n+1,rs}(F_v)$ we define the orbital integral associated to $s$ as the distribution

$$\displaystyle O(s,\varphi')=\int_{GL_n(F_v)} \varphi'(h^{-1}sh) \eta_{E_v/F_v}(h) dh,\;\;\; \varphi'\in \mathcal{S}(S_{n+1}(F_v))$$

\noindent For $\gamma=(\gamma_1,\gamma_2)\in G'_{rs}$ the element $s=\nu(\gamma_1^{-1}\gamma_2)\in S_{n+1}$ is regular semisimple and this defines a surjection $G'_{rs}\twoheadrightarrow S_{n+1,rs}$. Moreover, for all $\gamma\in G'_{rs}(F_v)$ and all $f'\in \mathcal{S}(G'(F_v))$, we have the equality

$$\displaystyle O(\gamma,f')=\left\{
    \begin{array}{ll}
        O(s,\widetilde{f}') & \mbox{ if } n \mbox{ is even,} \\
        \eta'_v(\gamma_1^{-1}\gamma_2)O(s,\widetilde{f}') & \mbox{ if } n \mbox{ is odd.}
    \end{array}
\right.
$$

\noindent where $s=\nu(\gamma_1^{-1}\gamma_2)$.

\vspace{2mm}

\noindent We can also define orbital integrals on the space $\mathcal{S}(\mathfrak{s}_{n+1}(F_v))$. Call an element $X\in \mathfrak{s}_{n+1}$ regular semisimple if it is so for the adjoint action of $GL_n$. Let us denote by $\mathfrak{s}_{n+1,rs}$ the open subset of regular semisimple elements. Then, for all $X\in \mathfrak{s}_{n+1,rs}(F_v)$ we can define an orbital integral by

$$\displaystyle O(X,\varphi')=\int_{GL_n(F_v)} \varphi'(h^{-1}Xh) \eta_{E_v/F_v}(h) dh,\;\;\; \varphi'\in \mathcal{S}(\mathfrak{s}_{n+1}(F_v))$$

\noindent The Cayley map $\mathfrak{c}=\mathfrak{c}_{n+1}:X\mapsto (1+X)(1-X)^{-1}$ realizes a $GL_n$-equivariant isomorphism between the open subsets $\mathfrak{s}^\circ_{n+1}=\{X\in \mathfrak{s}_{n+1};\; \det(1-X)\neq 0\}$ and $S^\circ_{n+1}=\{s\in S_{n+1};\; \det(1+s)\neq 0\}$. Let $\omega'\subset \mathfrak{s}^\circ_{n+1}(F_v)$ and $\Omega'\subset S^\circ_{n+1}(F_v)$ be open and closed $GL_n(F_v)$-invariant neighborhoods of $0$ and $1$ respectively such that the Cayley map restricts to an analytic isomorphism between $\omega'$ and $\Omega'$ preserving measures. For all $\varphi'\in \mathcal{S}(S_{n+1}(F_v))$, we define a function $\varphi'_{\natural}\in \mathcal{S}(\mathfrak{s}_{n+1}(F_v))$ by

$$\displaystyle \varphi'_{\natural}(X)=\left\{
    \begin{array}{ll}
        \varphi'(\mathfrak{c}(X)) & \mbox{ if } X\in\omega' \\
        0 & \mbox{ otherwise}
    \end{array}
\right.
$$

\noindent Then for all $X\in \omega'_{rs}=\omega'\cap \mathfrak{s}_{n+1,rs}(F_v)$ and all $\varphi'\in \mathcal{S}(S_{n+1}(F_v))$ we have

$$\displaystyle O(\mathfrak{c}(X),\varphi')=O(X,\varphi'_\natural)$$

\subsection{Correspondence of orbits and transfer}\label{correspondence}

\noindent We now recall the correspondence between orbits following [Zh1, \S 2.4]. We will denote by $H_1'\backslash G'/H_2'$, $H\backslash G/H$, $S_{n+1}/GL_n$ and $U(V)/U(W)$ the geometric quotients of $G'$, $G$, $S_{n+1}$ and $U(V)$ by $H_1'\times H_2'$, $H\times H$, $GL_n$ and $U(W)$ respectively (the last two actions being given by conjugation). We will also write $\left(H_1'\backslash G'/H_2'\right)_{rs}$, $\left(H\backslash G/H\right)_{rs}$, $(S_{n+1}/GL_n)_{rs}$ and $(U(V)/U(W))_{rs}$ for the regular semisimple loci in these geometric quotients. These are the image of $G'_{rs}$, $G_{rs}$, $S_{n+1,rs}$ and $U(V)_{rs}$ by the natural projections. The maps $(\gamma_1,\gamma_2)\in G'\mapsto \nu(\gamma_1^{-1}\gamma_2)$ and $(\delta_W,\delta_V)\in G\mapsto \delta_W^{-1}\delta_V$ induce isomorphisms

$$\displaystyle H_1'\backslash G'/H_2'\simeq S_{n+1}/GL_n \mbox{ and } H\backslash G/H\simeq U(V)/U(W)$$

\noindent and similarly for the regular semisimple loci. Moreover, there is a natural isomorphism (see [Zh1, \S 3.1])

$$\displaystyle H_1'\backslash G'/H_2'\simeq H\backslash G/H \leqno (1)$$

\noindent which preserves the regular semisimple loci. For all field extension $k$ of $F$ we have $\left(H_1'\backslash G'/H_2'\right)_{rs}(k)=H_1'(k)\backslash G'_{rs}(k)/H_2'(k)$ and $H(k)\backslash G_{rs}(k)/H(k)$ is a subset of $\left(H\backslash G/H\right)_{rs}(k)$. The above isomorphism thus induces injections

$$\displaystyle H(k)\backslash G_{rs}(k)/H(k)\hookrightarrow H_1'(k)\backslash G'_{rs}(k)/H_2'(k) \leqno (2)$$

\noindent and

$$\displaystyle U(V)_{rs}(k)/U(W)(k)\hookrightarrow S_{n+1,rs}(k)/GL_n(k) \leqno (3)$$

\noindent This last map admits the following explicit description. Choosing a basis of $V$ whose last element is $e$ we get an embedding $U(V)(k)\hookrightarrow GL_{n+1}(k\otimes_F E)$. By [Zh3, lemma 2.3] any regular semisimple element $x\in U(V)_{rs}(k)$ is $GL_{n}(k\otimes_F E)$-conjugated to a regular semisimple element of $S_{n+1}(k)$ which is unique up to $GL_n(k)$-conjugation. The $GL_n(k)$-conjugacy class of this element is exactly the image of $x$ by the map (3).

\vspace{2mm}

\noindent We have a similar situation at the level of Lie algebras: we have a canonical isomorphism between geometric quotients

$$\displaystyle \mathfrak{s}_{n+1}/GL_n\simeq \mathfrak{u}(V)/U(W) \leqno (4)$$

\noindent which preserves the regular semisimple loci $(\mathfrak{s}_{n+1}/GL_n)_{rs}=\mathfrak{s}_{n+1,rs}/GL_n$ and $(\mathfrak{u}(V)/U(W))_{rs}=\mathfrak{u}(V)_{rs}/U(W)$. For all field extension $k$ of $F$ this induces an injection

$$\displaystyle \mathfrak{u}(V)_{rs}(k)/U(W)(k)\hookrightarrow \mathfrak{s}_{n+1,rs}(k)/GL_n(k) \leqno (5)$$

\noindent We now define, following [Zh2, \S 4.1], two families of {\it transfer factors} $\Omega_v: G'_{rs}(F_v)\to \C^\times$ and $\omega_v: \mathfrak{s}_{n+1,rs}(F_v)\to \C^\times$, $v$ a place of $F$, satisfying the following conditions:

\vspace{2mm}

\begin{itemize}
\renewcommand{\labelitemi}{$\bullet$}
\item For all $v$ and all $\gamma\in G'_{rs}(F_v)$ (resp. all $X\in \mathfrak{s}_{n+1,rs}(F_v)$) we have $\Omega_v(h_1\gamma h_2)=\eta_v(h_2)\Omega_v(\gamma)$ (resp. $\omega_v(h^{-1}Xh)=\eta_{E_v/F_v}(h)\omega_v(X)$) for all $(h_1,h_2)\in H_1'(F_v)\times H_2'(F_v)$ (resp. for all $h\in GL_n(F_v)$);

\item For all $\gamma\in G'_{rs}(F)$ (resp. all $X\in \mathfrak{s}_{n+1,rs}(F)$), we have the product formula $\prod_v \Omega_v(\gamma)=1$ (resp. $\prod_v \omega_v(X)=1$) where almost all terms in the product are equal to $1$.
\end{itemize}

\vspace{2mm}

\noindent Let $v$ be a place of $F$. For all $s\in S_{n+1,rs}(F_v)$ and all $X\in \mathfrak{s}_{n+1,rs}(F_v)$, we set

$$\displaystyle \Omega_v(s)=\eta'_v\left(\det(s)^{-[\frac{n+1}{2}]} \det(e_{n+1},e_{n+1}s,\ldots,e_{n+1}s^n) \right) \leqno (6)$$

$$\displaystyle \omega_v(X)=\eta'_v\left(\det(e_{n+1},e_{n+1}X,\ldots,e_{n+1}X^n)\right) \leqno (7)$$

\noindent where $e_{n+1}=(0,\ldots,0,1)$ and $\eta'_v$ is the local component at $v$ of the character $\eta':\A_E^\times\to\C^\times$ extending $\eta_{E/F}$ that we fixed at the beginning. Note that (see the proof of [Zh1, lemma 3.5])

$$\displaystyle \Omega_v(\mathfrak{c}(X))=\eta_v(2)^{n(n+1)/2}\omega_v(X)\leqno (8)$$

\noindent for all $X\in \mathfrak{s}_{n+1,rs}(F_v)$ sufficiently close to $0$. Finally for all $\gamma=(\gamma_1,\gamma_2)\in G'_{rs}(F_v)$, we set

$$\displaystyle \Omega_v(\gamma)=\left\{
    \begin{array}{ll}
        \Omega_v(s) & \mbox{ if } n \mbox{ is even,} \\
        \eta_v'(\gamma_1^{-1}\gamma_2)\Omega_v(s) & \mbox{ if } n \mbox{ is odd.}
    \end{array}
\right.
$$

\noindent where $s=\nu(\gamma_1^{-1}\gamma_2)$. For future reference, we record the following formula

$$\displaystyle \Omega_v(\gamma)O(\gamma,f')=\Omega_v(s)O(s,\widetilde{f}') \leqno (9)$$

\noindent for all $f'\in \mathcal{S}(G'(F_v))$, all $\gamma\in G'_{rs}(F_v)$ and where we have set $s=\nu(\gamma_1^{-1}\gamma_2)$.

\vspace{2mm}

\noindent Using the transfer factors we can define the notion of matching functions as follows. Let $v$ be a place of $F$. We say that functions $f'\in \mathcal{S}(G'(F_v))$ and $f\in \mathcal{S}(G(F_v))$ {\it match} each other or that they are {\it smooth transfer} of each other if we have the equality

$$\displaystyle O(\delta,f)=\Omega_v(\gamma)O(\gamma,f')$$

\noindent for every $\delta\in G_{rs}(F_v)$ and $\gamma\in G'_{rs}(F_v)$ whose orbits correspond to each other via the embedding (2). Similarly, we say that functions $\varphi'\in \mathcal{S}(\mathfrak{s}_{n+1}(F_v))$ and $\varphi\in \mathcal{S}(\mathfrak{u}(V)(F_v))$ {\it match} each other or that they are {\it smooth transfer} of each other if we have the equality

$$\displaystyle O(X,\varphi_W)=\omega_v(Y)O(Y,\varphi')$$

\noindent for every $X\in \mathfrak{u}(V)_{rs}(F_v)$ and $Y\in \mathfrak{s}_{n+1,rs}(F_v)$ whose orbits correspond to each other via the embedding (5).

\vspace{2mm}

\noindent If the place $v$ splits in $E$ then the existence of smooth transfer is easy (see [Zh1, Proposition 2.5]). One of the main achievement of [Zh1] was to prove the existence of smooth transfer for nonarchimedean places. In other words, we have the following (see [Zh1, Theorem 2.6]):

\begin{theo}[Zhang]\label{smooth transfer Zhang}
Let $v$ be a nonarchimedean place of $F$.
\begin{enumerate}[(i)]
\item For every function $f'\in \mathcal{S}(G'(F_v))$ there exists a function $f\in \mathcal{S}(G(F_v))$, matching $f'$ and conversely for every function $f\in \mathcal{S}(G(F_v))$ there exists a function $f'\in \mathcal{S}(G'(F_v))$ which matches $f$.

\item For every function $\varphi'\in \mathcal{S}(\mathfrak{s}_{n+1}(F_v))$ there exists a function $\varphi\in \mathcal{S}(\mathfrak{u}(V)(F_v))$ matching $\varphi'$ and conversely for every function $\varphi\in \mathcal{S}(\mathfrak{u}(V)(F_v))$ there exists a function $\varphi'\in \mathcal{S}(\mathfrak{s}_{n+1}(F_v))$ which matches $\varphi$.
\end{enumerate}
\end{theo}

\noindent One of the main ingredient in the proof of Zhang was the following (see [Zh1, Theorem 4.17])

\begin{theo}[Zhang]\label{transfer and FT}
Let $v$ be any place of $F$. If $\varphi\in \mathcal{S}(\mathfrak{u}(V)(F_v))$ and $\varphi'\in \mathcal{S}(\mathfrak{s}_{n+1}(F_v))$ match then so do $\epsilon(\frac{1}{2},\eta_{E_v/F_v},\psi)^{n(n+1)/2}\widehat{\varphi}$ and $\widehat{\varphi'}$.
\end{theo}

\noindent In a recent paper Xue ([Xue]) was able to extend Zhang results to obtain a weak version of smooth transfer at archimedean places (which however is sufficient for many global applications). In order to state Xue's result in the generality that we need, we have to vary the hermitian space $W$. Let us denote momentarily the groups $G$ and $H$ by $G^W$ and $H^W$. To every hermitian space $W'$ of rank $n$ over $E$ we can associate similar groups $G^{W'}$ and $H^{W'}$ and replacing $W$ by $W'$ everywhere in the previous paragraphs we have a notion of matching between test functions in $\mathcal{S}(G'(F_v))$ and test functions in $\mathcal{S}(G^{W'}(F_v))$, $v$ a place of $F$. Then, Xue's result reads as follows: 

\begin{theo}[Xue]\label{smooth transfer Xue}
Let $v$ be an archimedean place of $F$. Then, the space of functions $f'\in \mathcal{S}(G'(F_v))$ admitting a smooth transfer to $\mathcal{S}(G^{W'}(F_v))$ for all hermitian space $W'$ of rank $n$ over $E$ is dense in $\mathcal{S}(G'(F_v))$. Similarly, the space of functions $f\in \mathcal{S}(G^W(F_v))$ such that there exists a function $f' \in \mathcal{S}(G'(F_v))$ matching $f$ and with the property that for all hermitian space $W'$ of rank $n$ over $E$ for which $W'_v\not\simeq W_v$ the function $f'$ match $0\in \mathcal{S}(G^{W'}(F_v))$, is dense in $\mathcal{S}(G(F_v))$.
\end{theo}

\noindent Parallel to the existence of smooth transfer there is also a fundamental lemma for the case at hand. This fundamental lemma has been proved by Yun in (sufficiently large) positive characteristic [Yu] and extended to the characteristic zero case by J.Gordon in the appendix to [Yu]. It can be stated as follows:

\begin{theo}[Yun-Gordon]\label{fundamental lemma}
There exists a constant $c(n)$ depending only on $n$ such that for every place $v$ of $F$ of residual characteristic greater than $c(n)$ the following holds: If $W_v$ admits a self-dual lattice $L_v$ then the function $f'_v=\mathbf{1}_{K'_v}$ match the function $f_v=\mathbf{1}_{G(\mathcal{O}_v)}$ where we have defined a model of $G$ over $\mathcal{O}_v$ using the self-dual lattice $L_v$, otherwise the function $f'_v=\mathbf{1}_{K'_v}$ match the function $f_v=0$.
\end{theo}

\subsection{Transfer of spherical characters, Zhang's conjecture and Ichino-Ikeda conjecture}\label{transfer of spherical characters}

We shall say of a function $f\in \mathcal{S}(G(\A))$ that it is {\it nice} if it satisfies the following conditions:

\vspace{2mm}

\begin{itemize}
\renewcommand{\labelitemi}{$\bullet$}

\item $f$ is factorizable: $f=\otimes_v f_v$;

\item There exists a nonarchimedean place $v_1$ of $F$ and a finite union $\Omega_1$ of cuspidal Bernstein components of $G(F_{v_1})$ such that $f_{v_1}\in \mathcal{S}(G(F_{v_1}))_{\Omega_1}$;
\end{itemize}

\noindent We define the notion of {\it nice} function on $G'(\A)$ similarly. To state the next theorem, we will need to consider more than one pair of hermitian spaces $(W,V)$. Recall that we have an orthogonal decomposition $V=W\oplus^\perp Ee$ where $(e,e)=1$. To any (isomorphism class of) $n$-dimensional hermitian space $W'$ over $E$ we associate the pair $(W',V')$ where $V'=W'\oplus^\perp Ee$. Using such a pair we may construct a new pair $(H^{W'},G^{W'})$ of reductive groups over $F$ where $H^{W'}=U(W')$ and $G^{W'}=U(W')\times U(V')$. Note that if $W'=W$ then $(H^{W'},G^{W'})=(H,G)$. The discussions of the previous paragraphs of course apply verbatim to $(H^{W'},G^{W'})$. In particular we have a notion of matching between functions in $\mathcal{S}(G^{W'}(F_v))$ and $\mathcal{S}(G'(F_v))$, $v$ a place of $F$, and a notion of nice function on $G^{W'}(\A)$. We shall say that a nice function $f'\in \mathcal{S}(G'(\A))$ match a tuple of nice functions $(f^{W'})_{W'}$, $f^{W'}\in \mathcal{S}(G^{W'}(\A))$ and $W'$ running over all isomorphism classes of $n$-dimensional hermitian spaces over $E$, if for all $W'$ and all place $v$ of $F$ the functions $f'_v$ and $f^{W'}_v$ match. Comparing two (simple) global relative trace formulae that have been proposed by Jacquet and Rallis ([JR]), Zhang proves the following (see [Zh1, Proposition 2.10] and [Zh2, Theorem 4.3]):

\begin{theo}[Zhang]\label{comparison RTF}
Let $\pi$ be an abstractly $H$-distinguished cuspidal automorphic representation of $G(\A)$ such that $BC(\pi)$ is cuspidal and for all non-split archimedean place $v$ the representation $\pi_v$ is tempered. Let $f\in \mathcal{S}(G(\A))$ and $f'\in \mathcal{S}(G'(\A))$ be nice functions and assume that there exists a tuple $(f^{W'})_{W'}$, $f^{W'}\in \mathcal{S}(G^{W'}(\A))$, of nice functions matching $f'$ such that $f^W=f$. Then, we have

$$\displaystyle J_\pi(f)=2^{-2}L(1,\eta_{E/F})^{-2}I_{BC(\pi)}(f')$$

\end{theo}

\vspace{2mm}

\begin{rem}
The above theorem differs slightly from [Zh1, Proposition 2.10] and [Zh2, Theorem 4.3] for the following reasons:
\begin{itemize}
\item First we are not assuming that our test functions are supported in the regular semi-simple locus (this is one of the requirement that Zhang imposes on his 'nice' functions);

\item We are using test functions that are not necessarily compactly supported (since they are only rapidly decreasing, i.e. in the Schwartz space, at the archimedean places). This is necessary if we want to apply this theorem in conjunction with Xue's result (Theorem \ref{smooth transfer Xue}) as the dense subspace of 'transferable' test functions that he constructs has no reason to consists of compactly supported functions;

\item We have drop the condition that $\pi$ is supercuspidal at one split places but we add the assumption that $\pi$ must be tempered at all non-split archimedean places.
\end{itemize}

This last point is only minor: dropping the condition of $\pi$ being supercuspidal at one split place can be done by slightly modifying Zhang's original argument when he separates contributions in the spectral side of the simple trace formulas and using the recent extension by Ramakrishnan of his 'mild Tchebotarev theorem for $GL(n)$' [Ra]. Moreover, here we assume that $\pi$ is tempered at all non-split archimedean places in order to have a result independent of the local Gan-Gross-Prasad conjecture for generic $L$-packets at archimedean places (this conjecture is currently only known for {\it tempered} $L$-packets at these places, see \S \ref{local GGP}), this conjecture was granted as a working hypothesis in [Zh2]. The first point is more serious and to get rid of this assumption on the support we have to use the recent works of Zydor [Zy] on regularization of the geometric side of the Jacquet-Rallis trace formulae and Chaudouard-Zydor [CZ] on extending the transfer to singular orbital integrals. Finally, the extension to rapidly decreasing functions (at the archimedean places) is an easy matter using basic estimates on these functions. A convenient way to do this is to introduce some norms on the automorphic quotient $[G]=G(F)\backslash G(\mathbb{A})$. We give definitions and basic properties of these norms in the appendix. Also, for convenience of the reader, we provide in appendix \ref{appendix} a full proof of Theorem \ref{comparison RTF} using the aforementioned results of Zydor and Chaudouard-Zydor
\end{rem}

\vspace{2mm}

\noindent Thanks to the theory of Rankin-Selberg convolution due to Jacquet, Piatetski-Shapiro and Shalika ([JPSS]), for every cuspidal automorphic representation $\pi$ of $G(\A)$ whose base change is cuspidal we know a factorization of the global spherical character $I_{BC(\pi)}$. More precisely, if $f'\in \mathcal{S}(G'(\A))$ is completely factorizable, we have (see [Zh2, Proposition 3.6])

$$\displaystyle I_{BC(\pi)}(f')=L(1,\eta_{E/F})^2 \frac{L(1/2,BC(\pi))}{L(1,\pi,Ad)}\prod_v I^\natural_{BC(\pi)_v}(f'_v) \leqno (1)$$

\noindent An immediate consequence of this factorization of Theorem \ref{comparison RTF} and of the multiplicity one theorems of Aizenbud-Gourevitch-Rallis-Schiffmann [AGRS] of of Jiang-Sun-Zhu [JSZ] is the following:

\begin{cor}\label{cor RTF}
Let $\pi$ be a globally $H$-distinguished (i.e. such that $J_\pi\neq 0$) cuspidal automorphic representation of $G(\A)$ satisfying the two following conditions:

\begin{itemize}
\item For all non-split archimedean place $v$ of $F$, the representation $\pi_v$ is tempered;

\item There exists a split nonarchimedean place $v_0$ such that $\pi_{v_0}$ is supercuspidal.
\end{itemize}

\noindent Then for all place $v_1$ of $F$ different from $v_0$ where $\pi$ is tempered, there exists a constant $C(\pi_{v_1})\in \C$ such that for all pair $(f_{v_1},f'_{v_1})\in \mathcal{S}(G(F_{v_1}))\times \mathcal{S}(G'(F_{v_1}))$ of matching functions with the property that $f'_{v_1}$ has a matching test function $f_{v_1}^{W'}\in \mathcal{S}(G^{W'}(F_{v_1}))$ for all hermitian space $W'$ of rank $n$, we have

$$\displaystyle J_{\pi_{v_1}}(f_{v_1})=C(\pi_{v_1})I_{BC(\pi_{v_1})}(f'_{v_1})$$
\end{cor}

\vspace{2mm}

\begin{rem}
Note that the condition of matching of the function $f'_{v_1}$ is empty if $v_1$ is nonarchimedean (by Theorem \ref{smooth transfer Zhang}) or splits in $E$.
\end{rem}

\vspace{3mm}

\noindent\ul{Proof}: Let $\mathbb{A}^{v_1}$ denote the adeles outside of $v_1$ and let $f^{v_1}=\prod_{v\neq v_1}\in \mathcal{S}(G(\mathbb{A}^{v_1}))$ be a factorizable test function. By the multiplicity one results of [AGRS] and [JSZ] and \ref{local spherical characters} (1), there exists a constant $C\in \C$ such that

$$\displaystyle J_\pi(f^{v_1}\otimes f_{v_1})=C J_{\pi_{v_1}}(f_{v_1})$$

\noindent for all $f_{v_1}\in\mathcal{S}(G(F_{v_1}))$. Since $J_\pi\neq 0$ we may choose the function $f^{v_1}$ so that $C\neq 0$. Moreover, up to replacing $f_{v_0}$ by its projection to $\mathcal{S}(G(F_{v_0}))_{\Omega_0}$, where $\Omega_0$ denotes the Bernstein component of $\pi_{v_0}$, we may assume that $f_{v_0}\in\mathcal{S}(G(F_{v_0}))_{\Omega_0}$. Then, for all $f_{v_1}\in\mathcal{S}(G(F_{v_1}))$ the function $f=f^{v_1}\otimes f_{v_1}$ is nice. By Theorem \ref{smooth transfer Xue} of Xue we can also arrange $f^{v_1}$ such that for all non-split archimedean place $v$ the function $f_v$ admits a transfer $f'_v\in \mathcal{S}(G'(F_v))$ which itself admits a transfer to $\mathcal{S}(G^{W'}(F_v))$ for all rank $n$ hermitian space $W'$ over $E$. By Theorem \ref{smooth transfer Zhang} of Zhang and our choice of $f^{v_1}$ at non-split archimedean places, we can find a factorizable test function ${f'}^{v_1}=\prod_{v\neq v_1} f'_v\in \mathcal{S}(G'(\mathbb{A}^{v_1}))$ matching $f^{v_1}$ and moreover admitting a matching to $\mathcal{S}(G^{W'}(\mathbb{A}^{v_1}))$ for all $W'$. Since $v_0$ splits in $E$, there exists an isomorphism $G'_{v_0}\simeq G_{v_0}\times G_{v_0}$ so that $BC(\pi_{v_0})=\pi_{v_0}\boxtimes \pi_{v_0}^\vee$ and then every function $\varphi'=\varphi_1\otimes \varphi_2\in \mathcal{S}(G'(F_{v_0}))\simeq \mathcal{S}(G(F_{v_0}))\otimes \mathcal{S}(G(F_{v_0}))$ match the function $\varphi_1\ast \varphi_2^\vee\in \mathcal{S}(G(F_{v_0}))$ where $\varphi_2^\vee(g):=\varphi_2(g^{-1})$ and $\ast$ denotes the convolution product. From this we easily infer that we can choose $f'_{v_0}\in \mathcal{S}(G'(F_{v_0}))_{\Omega_0\times\Omega_0^\vee}$, where $\Omega_0^\vee$ is the Bernstein component dual to $\Omega_0$. Then, for all $f'_{v_1}\in \mathcal{S}(G'(F_{v_1}))$ the function $f':={f'}^{v_1}\otimes f'_{v_1}$ is nice. Now the result follows from Theorem \ref{comparison RTF} applied to $f=f^{v_1}\otimes f_{v_1}$ and $f'={f'}^{v_1}\otimes f'_{v_1}$ where $(f_{v_1},f'_{v_1})\in \mathcal{S}(G(F_{v_1}))\times \mathcal{S}(G'(F_{v_1}))$ is any pair of matching functions with $f'_{v_1}$ satisfying the condition of the statement (Note that the assumption that $\pi_{v_0}$ is supercuspidal implies that $BC(\pi)$ is cuspidal). $\blacksquare$

\vspace{3mm}

\noindent In [Zh2, conjecture 4.4], Zhang makes the following conjecture

\begin{conj}\label{Zhang conj}
Let $v$ be a place of $F$ and let $\pi_v=\pi_{n,v}\boxtimes \pi_{n+1,v}$ be an irreducible tempered unitary $H_v$-distinguished representation of $G(F_v)$. Then, for all matching functions $f\in \mathcal{S}(G(F_v))$ and $f'\in \mathcal{S}(G'(F_v))$, we have

$$\displaystyle I_{BC(\pi_v)}(f'_v)=\kappa_v(\pi_v) L(1,\eta_{E_v/F_v})^{-1} J_{\pi_v}(f_v)$$

\noindent the constant $\kappa_v(\pi_v)$ being given by

$$\displaystyle \kappa_v(\pi_v)=\lv \tau\rv_{E_v}^{(d_n+d_{n+1})/2} \left(\frac{\epsilon(1/2,\eta_{E_v/F_v},\psi_v)}{\eta'_v(-2\tau)}\right)^{n(n+1)/2} \omega_{BC(\pi_{n,v})}(\tau)$$

\noindent where $\omega_{BC(\pi_{n,v})}$ denotes the central character of $BC(\pi_{n,v})$ and $d_n=\binom{n}{3}$, $d_{n+1}=\binom{n+1}{3}$. 
\end{conj}

\begin{rem}
The above conjecture actually differs slightly from [Zh2,conjecture 4.4]. Indeed, there is a discrepancy in the definition of the constant $\kappa_v(\pi_v)$. In {\it loc. cit.} there is an extra factor $\eta_v(disc(W))$ and the factor $\eta'_v(-2\tau)$ is replaced by $\eta'_v(\tau)$. This discrepancy seems to come from lemma 9.1 of [Zh2] on the compatibility between the transfer factors on the group and on the Lie algebra (see \S \ref{correspondence}(8)). Of course, this difference has no impact for global applications since in any case if $\pi$ is an automorphic representation then $\prod_v \kappa_v(\pi_v)=1$.
\end{rem}

\noindent Obviously, we may deduce from the conjunction of the above conjecture, of Theorems \ref{comparison RTF}, \ref{smooth transfer Zhang} and \ref{smooth transfer Xue} and of the factorization (1) some instances of the Ichino-Ikeda conjecture (as stated in the introduction). In [Zh2], Zhang was able to verify conjecture \ref{Zhang conj} in certain particular cases (see [Zh2, Theorem 4.6]). More precisely, he proves the conjecture when either

\begin{itemize}
\item The place $v$ splits in $E$;

\item $v$ is nonarchimedean, $\pi_v$ is unramified and the residue characteristic of $v$ is sufficiently large;

\item $v$ is nonarchimedean and $\pi_v$ is supercuspidal.
\end{itemize}

\noindent The main goal of this paper is to prove conjecture \ref{Zhang conj} for all nonarchimedean places $v$. Namely, we prove

\begin{theo}\label{theo principal}
Conjecture \ref{Zhang conj} holds for all nonarchimedean place $v$ of $F$.
\end{theo}

\noindent As in [Zh2], this theorem has consequences for the Ichino-Ikeda conjecture. Namely, we will deduce from it the following

\begin{theo}\label{theo II}
Assume that all the archimedean places of $F$ split in $E$ and let $\pi$ be a cuspidal automorphic representation of $G(\A)$ which is everywhere tempered and such that there exists a nonarchimedean place $v$ of $F$ with $BC(\pi_v)$ is supercuspidal. Then conjecture \ref{II conj} holds for $\pi$.
\end{theo}

\noindent The proofs of Theorems \ref{theo principal} and \ref{theo II} will be given in \S \ref{proof of main theo} and \S \ref{proof of II theo}.

\subsection{A globalization result}\label{globalization section}

\noindent Until the end of this paragraph, we make the following assumption

\begin{center}
the hermitian space $W$ is anisotropic
\end{center}

\noindent This means that $H=U(W)$ is an anisotropic group over $F$. 

\vspace{2mm}

\noindent Let $v_1$ be a nonarchimedean place of $F$ which is inert in $E$ and let $S$ be a finite set of nonarchimedean places of $F$ which split in $E$. Let $\sigma$ be a unitary supercuspidal representation of $G(F_S)$. Recall that $Temp_{H}(G_{v_1})$ denotes the set of (isomorphism classes of) tempered irreducible $H(F_{v_1})$-distinguished representations $\pi_1$ of $G(F_{v_1})$ (see \S \ref{local GGP}). Let $Irr_{v_1,\sigma,H}(G)$ be the set of irreducible representations $\pi_1\in Irr(G_{v_1})$ for which there exists a cuspidal automorphic representation $\pi$ of $G(\A)$ which is globally $H$-distinguished (i.e. such that $J_\pi\neq 0$) such that $\pi_{v_1}\simeq \pi_1$ and $\pi_S\simeq \sigma\otimes \chi$ for some unramified character $\chi\in \Psi_{unit}(G_S)$. The goal of this section is to prove the following result:

\begin{prop}\label{prop glob}
The set $Irr_{v_1,\sigma,H}(G)\cap Temp(G_{v_1})$ is dense in $Temp_{H}(G_{v_1})$.
\end{prop}

\noindent The proof of this proposition follows closely that of Corollary A.8 in [ILM, appendix A]. We will need a lemma which is the analog of lemma A.2 of {\it loc. cit}. Before stating it we need to introduce some more notations.

\vspace{2mm}

\noindent Let $P=MN$ be a parabolic subgroup of $G_{v_1}$ and let $\sigma$ be a square-integrable representation of $M(F_{v_1})$. We will say that the tempered representation $i_P^{G_{v_1}}\sigma$ of $G(F_{v_1})$ is {\it regular} if for all $w\in W(G_{v_1},M)$ we have $w\sigma\not\simeq \sigma$. Recall that this implies that the representation $i_M^{G_{v_1}}\sigma$ is irreducible. We will denote by $Temp_{reg}(G_{v_1})$ the set of all regular tempered representations of $G(F_{v_1})$. It is an open subset of $Temp(G_{v_1})$. Recall that in \S \ref{base change} for all $c>0$ we have defined subsets $Irr_{\leqslant c}(U(W)_{v_1})$ and $Irr_{\leqslant c}(U(V)_{v_1})$ of $Irr(U(W)_{v_1})$ and $Irr(U(V)_{v_1})$ respectively. In what follows, we set $Irr_{\leqslant c}(G_{v_1})=Irr_{\leqslant c}(U(W)_{v_1})\boxtimes Irr_{\leqslant c}(U(V)_{v_1})$ and $Irr_{unit,\leqslant c}(G_{v_1})=Irr_{unit}(G_{v_1})\cap Irr_{\leqslant c}(G_{v_1})$.

\vspace{2mm}

\begin{lem}\label{lemma 1}
Let $0<c<\frac{1}{2}$. Then, $Temp_{reg}(G_{v_1})$ is open in $Irr_{unit,\leqslant c}(G_{v_1})$ (for the Fell topology).
\end{lem}

\noindent\ul{Proof}: The proof is the same as in [ILM, lemma A.2] the key fact being that the real exponents of tempered representations of unitary groups are all half integers. In {\it loc. cit.} the authors use the work of Muic on generic square-integrable representations of classical groups ([Mui]) to deduce this fact for tempered generic representations. However, as already noted in [ILM, Remark A.3], the same result holds for all tempered representation thanks to the work of M\oe{}glin and M\oe{}glin-Tadic on classification of square-integrable representations of classical groups ([Moe1],[MT]). Note that the {\it basic assumption} made by M\oe{}glin and Tadic (see [MT,\S 2] for a precise statement) to prove their classification is now known since it follows from the canonical normalization of intertwining operators for unitary groups due to Mok ([Mok,Proposition 3.3.1]) and Kaletha-Minguez-Shin-White ([KMSW, lemma 2.2.3]) together with the classical reducibility criterion of Silberger and Harish-Chandra ([Sil1,\S 5.4], [Sil2, lemma 1.2; lemma 1.3]). For quasi-split unitary groups, a different proof has been given by M\oe{}glin ([Moe2]) using twisted endoscopy. For a proof of the basic assumption for quasisplit symplectic and orthogonal groups using the normalization of intertwining operators see [Xu,Proposition 3.2]. $\blacksquare$

\vspace{3mm}

\noindent\ul{Proof of Proposition} \ref{prop glob}: By [Beu1, Corollary 8.6.1], the closure of $Temp_H(G_{v_1})$ is an union of connected components of $Temp(G_{v_1})$. Hence $Temp_{H,reg}(G_{v_1}):=Temp_H(G_{v_1})\cap Temp_{reg}(G_{v_1})$ is dense in $Temp_H(G_{v_1})$. Let $\pi_1\in Temp_{H,reg}(G_{v_1})$. It is sufficient to show that $\pi_1$ belongs to the closure of $Irr_{v_1,\sigma,H}(G)\cap Temp(G_{v_1})$. Since $\sigma$ is $H(F_S)$-distinguished (see \S \ref{local spherical characters}), by [SV] Theorem 6.2.1 and Theorem 6.4.1 we know that $\pi_1$ and $\sigma$ belong to the support of the Plancherel measures for $L^2(H(F_{v_1})\backslash G(F_{v_1}))$ and $L^2(H(F_S)\backslash G(F_S))$ respectively. From [SV, Theorem 16.3.2], it follows that there exists a sequence of globally $H$-distinguished automorphic representations $\pi_n$ of $G(\A)$ such that $\pi_{n,v_1}\to \pi_1$ and $\pi_{n,S}\to \sigma$ for the Fell topology. Since $\sigma\otimes \Psi_{unit}(G_S)$ is open in $Irr_{unit}(G_S)$ we have $\pi_{n,S}\in \sigma\otimes \Psi_{unit}(G_S)$ for $n$ sufficiently large. This implies that $\pi_{n,v_1}$ belongs to $Irr_{v_1,\sigma,H}(G)$ and $\pi_n$, $BC(\pi_n)$ are cuspidal for $n$ sufficiently large. 

\vspace{2mm}

\noindent Set $c=\frac{1}{2}-\frac{1}{(n+1)^2+1}$. By Lemma \ref{lemma 2}, $\pi_{n,v_1}$ belongs to $Irr_{unit,\leqslant c}(G_{v_1})$ for $n$ sufficiently large. Hence, by Lemma \ref{lemma 1}, $\pi_{n,v_1}\in Temp_{H,reg}(G_{v_1})$ for $n$ sufficiently large and this ends the proof of the proposition. $\blacksquare$

\section{Proof of Zhang's conjecture}\label{section 3}

In this section we will prove Theorem \ref{theo principal} and Theorem \ref{theo II}. As Theorem \ref{theo principal} has already been proved by Zhang at every split place $v$, we only need to prove it at every nonarchimedean place $v$ of $F$ which is inert in $E$. Fix such a place $v$. We will now drop all the index $v$: $E/F=E_v/F_v$,$G=G_v$, $H=H_v$, $G'=G'_v$, $H_1'=H_{1,v}'$, $H_2'=H_{2,v}'$, $\psi=\psi_v$, $\psi_E=\psi_{E,v}$ and so on. Also, to ease notation we will just write $\mathfrak{s}=\mathfrak{s}_{n+1}$. Finally, we will now use {\it unnormalized} Haar measures (see \S \ref{measures}). In particular, Theorem \ref{theo principal} now takes the following form (see [Zh2, lemma 4.7]):

\begin{theo}\label{theo principal 2}
Let $\pi=\pi_{n}\boxtimes \pi_{n+1}$ be a $H$-distinguished irreducible tempered representation of $G(F)$. Then, for all matching functions $f\in \mathcal{S}(G(F))$ and $f'\in \mathcal{S}(G'(F))$, we have

$$\displaystyle I_{BC(\pi)}(f')=\kappa(\pi) J_{\pi}(f)$$

\noindent where

$$\displaystyle \kappa(\pi)=\lv \tau\rv_{E}^{(d_n+d_{n+1})/2} \left(\frac{\epsilon(1/2,\eta_{E/F},\psi)}{\eta'(-2\tau)}\right)^{n(n+1)/2} \omega_{BC(\pi_{n})}(\tau)$$
\end{theo}

\subsection{A result of Zhang on truncated local expansion of the spherical character $I_\Pi$}\label{truncated local expansion}

\noindent In this section we recall a result of Zhang [Zh2] on the existence of truncated local expansion for the spherical characters $I_\Pi$. This result is the main ingredient in the proof by Zhang of some particular cases of conjecture \ref{Zhang conj}. It will also play a crucial role in the proof of Theorem \ref{theo principal}.

\vspace{2mm}

\noindent Let us set

$$\displaystyle \xi_-=\begin{pmatrix} 0 & \ldots & \ldots & \ldots & 0 \\ \tau & \ddots & & & \vdots \\ 0 & \ddots & \ddots & & \vdots \\ \vdots & \ddots & \ddots & \ddots & \vdots \\ 0 & \ldots & 0 & \tau & 0
\end{pmatrix}\in \mathfrak{s}(F)$$

\noindent It is a regular nilpotent element for the $GL_n(F)$-action by conjugation (see [Zh2, lemma 6.1]). Zhang has defined a {\it regularized} orbital integral $\mu_{\xi_-}$ over the orbit of $\xi_-$ (see [Zh2, definition 6.10]). It is a $GL_n(F)$-invariant linear form $\varphi\in \mathcal{S}(\mathfrak{s}(F))\mapsto \mu_{\xi_-}(\varphi)$ which coincide with the usual orbital integral when the support of $\varphi$ intersect the orbit of $\xi_-$ in a compact set.

\vspace{2mm}

\noindent For all $X=\begin{pmatrix} A & u \\ v & w\end{pmatrix}\in \mathfrak{s}_{n+1}(F)$, we define

$$\displaystyle \Delta_-(X):=\det\left(v,vA,\ldots,vA^{n-1} \right)$$

\noindent Note that (see \S \ref{correspondence} (7))

$$\displaystyle \omega(X)=\eta(-1)^n \eta'(\Delta_-(X)) \leqno (1)$$

\noindent for all $X\in \mathfrak{s}_{rs}(F)$ and

$$\displaystyle \Delta_-(\xi_-)=(-1)^{n(n-1)/2}\tau^{n(n+1)/2} \leqno (2)$$

\vspace{2mm}

\noindent Let $r>m'>m>0$ be positive integers. In [Zh2, definition 8.1], Zhang defines a notion of $(m,m',r)$-admissible test functions on $G'(F)$. They span a finite dimensional subspace of $\mathcal{S}(G'(F))$. In what follows when we say that $(m,m',r)$ is {\it sufficiently large}, we shall mean that $m$ is sufficiently large, that $m'$ is sufficiently large depending on $m$ and that $r$ is sufficiently large depending on $(m,m')$. Recall that in \S \ref{orbital integrals} we have, using a Cayley map, associated to any function $f'\in \mathcal{S}(G'(F))$ a function $f'_\natural$ on $\mathfrak{s}(F)$. Also, in \S \ref{measures} we have defined a certain Fourier transform $\varphi\mapsto \widehat{\varphi}$ on $\mathcal{S}(\mathfrak{s}(F))$. We extract from [Zh2] the two following results (see Lemma 8.8, Theorem 8.5 of [Zh2] and their proofs):

\begin{prop}[Zhang]\label{const orb int}
Let $U$ and $\mathcal{Z}$ be relatively compact neighborhood of $1$ and $0$ in $G'(F)$ and $(\mathfrak{s}/GL_n)(F)$ respectively. Then, if $(m,m',r)$ is sufficiently large, for every $(m,m',r)$-admissible function $f'$ we have $Supp(f')\subseteq U$ and the function $X\in \mathcal{Z}_{rs}\mapsto \eta'(\Delta_-(X))O(X,\widehat{f}'_\natural)$ is constant and equal to $\eta'(\Delta_-(\xi_-)) \mu_{\xi_-}(\widehat{f}'_\natural)$. Moreover, we can find a $(m,m',r)$-admissible function $f'$ such that $\mu_{\xi_-}(\widehat{f}'_\natural)\neq 0$.
\end{prop}

\begin{theo}[Zhang]\label{trunc LE}
Let $\Pi=\Pi_n\boxtimes \Pi_{n+1}$ be an irreducible tempered representation of $G'(F)$.  Then, if $(m,m',r)$ is sufficiently large (depending on $\Pi$) we have the equality

$$\displaystyle I_\Pi(f')=\lv \tau\rv_E^{(d_n+d_{n+1})/2}\omega_{\Pi_n}(\tau)\mu_{\xi_-}(\widehat{f}'_\natural)$$

\noindent for all $(m,m',r)$-admissible function $f'$, where $d_n=\binom{n}{3}$ and $\omega_{\Pi_n}$ denotes the central character of $\Pi_n$.
\end{theo}

\noindent A direct consequence of Proposition \ref{const orb int} and Theorem \ref{trunc LE} is the following:

\begin{cor}\label{existence good functions}
Let $C\subseteq Temp(G')$ be a compact subset and let $U$ and $\mathcal{Z}$ be relatively compact neighborhood of $1$ and $0$ in $G'(F)$ and $(\mathfrak{s}/GL_n)(F)$ respectively.. Then, there exists a test function $f'\in \mathcal{S}(G'(F))$ satisfying the following conditions:

\begin{enumerate}[(i)]
\item $Supp(f')\subseteq U$ and the function $X\in \mathcal{Z}_{rs}\mapsto \eta'(\Delta_-(X)) O(X,\widehat{f}'_\natural)$ is constant and equal to $\eta'(\Delta_-(\xi_-))\mu_{\xi_-}(\widehat{f}'_\natural)$; 

\item $\mu_{\xi_-}(\widehat{f}'_\natural)\neq 0$;

\item For all $\Pi\in C$ we have the equality

$$\displaystyle I_\Pi(f')=\lv \tau\rv_E^{(d_n+d_{n+1})/2}\omega_{\Pi_n}(\tau)\mu_{\xi_-}(\widehat{f}'_\natural)$$
\end{enumerate}
\end{cor}

\noindent\ul{Proof}: For all $r>m'>m>0$, let us denote by $C[m,m',r]$ the set of $\Pi\in C$ such that the equality

$$\displaystyle I_\Pi(f')=\lv \tau\rv_E^{(d_n+d_{n+1})/2}\omega_{\Pi_n}(\tau)\mu_{\xi_-}(\widehat{f}'_\natural)$$

\noindent holds for all $(m,m',r)$-admissible function $f'$. Note that $C[m,m',r]$ is a closed subset of $C$. Obviously, by Proposition \ref{const orb int}, it suffices to show that if $(m,m',r)$ is sufficiently large then $C[m,m',r]=C$ and for that we may assume $C$ to be connected. By Theorem \ref{trunc LE} we have

$$\displaystyle \bigcup_{m_0>0}\bigcap_{m\geqslant m_0}\bigcup_{m_0'>m}\bigcap_{m'\geqslant m_0'}\bigcup_{r_0>m'}\bigcap_{r\geqslant r_0} C[m,m',r]=C$$

\noindent Now, by Baire category theorem, this implies that for $(m,m',r)$ sufficiently large the set $C[m,m',r]$ is not meager i.e. it has nonempty interior (since it is closed). By connectedness of $C$ and analyticity of $\Pi\mapsto I_\Pi$ this implies $C[m,m',r]=C$ and this ends the proof. $\blacksquare$

\subsection{Weak comparison of local spherical characters}\label{weak comparison}

\begin{prop}\label{prop weak comp}
For all $\pi\in Temp_H(G)$ there exists a nonzero constant $C(\pi)\in \C$ such that for all matching functions $f\in \mathcal{S}(G(F))$ and $f'\in \mathcal{S}(G'(F))$ we have

$$\displaystyle J_\pi(f)=C(\pi)I_{BC(\pi)}(f')$$

\noindent Moreover, the function $\pi\in Temp_H(G)\mapsto C(\pi)$ is analytic.
\end{prop}

\vspace{2mm}

\noindent\ul{Proof}: Assume that we have proved the existence of a constant $C(\pi)$ as in the proposition for a dense set of $\pi$ in $Temp_H(G)$. We claim that the proposition can be deduced from this. Indeed, for all $\pi\in Temp_H(G)$ we can define a constant $C(\pi)$ as follows: choose any set $(f,f')\in \mathcal{S}(G(F))\times \mathcal{S}(G'(F))$ of matching functions such that $I_{BC(\pi)}(f')\neq 0$ (the existence of such a pair follows from Theorem \ref{smooth transfer Zhang} and Theorem \ref{trunc LE}) and set $C(\pi)=J_\pi(f)I_{BC(\pi)}(f')^{-1}$. Of course this constant may {\it a priori} depend on the choice of $f$ and $f'$ but it follows from the analyticity of $\pi\mapsto J_\pi$ and $\Pi\mapsto I_\Pi$ and our assumption that in fact it is independent of such a choice. Still by analyticity of the spherical characters the equality of the proposition is true for all $\pi\in Temp_H(G)$ and all pair  of matching functions $(f,f')$ and the function $\pi\in Temp_H(G)\mapsto C(\pi)$ is analytic. Moreover it is nowhere zero since for all $\pi\in Temp_H(G)$ there exists $f\in \mathcal{S}(G(F))$ such that $J_\pi(f)\neq 0$ and there exists a $f'\in \mathcal{S}(G'(F))$ matching $f$ (by Theorem \ref{smooth transfer Zhang}).

\vspace{2mm}

\noindent We now prove the existence of a dense subset of $\pi$ satisfying the proposition. To this end we will use Proposition \ref{prop glob}. We first need to globalize the situation at hand. Let

\begin{itemize}
\item $\E/\F$ be a quadratic extension of number fields such that all archimedean places of $\F$ are nonsplit in $\E$ and $v_1$ be a place of $\F$ such that $\E_{v_1}/\F_{v_1}\simeq E/F$;

\item $\W$ a $n$-dimensional hermitian space over $\E$ such that

\begin{itemize}
\item for all archimedean place $v$ of $\F$ the group $U(\W)_v$ is anisotropic (in particular $\W$ is anisotropic); 
\item $\W_{v_0}\simeq W$.
\end{itemize}
\end{itemize}

\noindent We will set $\V=\W\oplus^\perp \E e$ where $(e,e)=1$ (so that $\V_{v_1}\simeq V$), $\H=U(\W)$ and $\G=U(\V)$. Let $v_0$ be two nonarchimedean places of $\F$ which split in $\E$ and let $\sigma$ be a supercuspidal representations of $\G(\F_{v_0})$. Applying Proposition \ref{prop glob} to $S=\{v_0\}$, we deduce the existence of a dense subset $D\subset Temp_H(G)$ such that for all $\pi\in D$ there exists a globally $\H$-distinguished cuspidal automorphic representation $\Pi$ of $\G(\A)$ such that $\Pi_{v_1}\simeq \pi$ and $\Pi_{v_0}$ is supercuspidal. Applying Corollary \ref{cor RTF} to such representations $\Pi$, which we remark are necessarily tempered at all archimedean places since $U(\W)$ is anisotropic there, we deduce for all $\pi\in D$ there exists a constant $C(\pi)$ as in the proposition. $\blacksquare$

\subsection{A local trace formula}\label{local TF}

Let $f_1,f_2\in \mathcal{S}(G(F))$. Then

\vspace{4mm}

\hspace{5mm} (1) The integral

$$\displaystyle J(f_1,f_2)=\int_{H(F)}\int_{H(F)} \int_{G(F)} f_1(h_1gh_2)f_2(g)dgdh_1dh_2$$

\hspace{12mm} is absolutely convergent.

\vspace{2mm}

\noindent This follows from [Zh1, Lemma A.4].

\vspace{2mm}

\noindent By [Beu1, Proposition 8.2.1(v)], we have

$$\displaystyle J(f_1,f_2)=\int_{Temp_H(G)} J_\pi(f_1) J_{\pi^\vee}(f_2) d\mu_G(\pi) \leqno (2)$$

\noindent where $d\mu_G(\pi)$ denotes the Harish-Chandra-Plancherel measure of $G(F)$. We also have (see \S \ref{orbital integrals} for the definition of $\widetilde{f}_i$)

$$\displaystyle J(f_1,f_2)=\int_{U(W)(F)}\int_{U(V)(F)} \widetilde{f}_1(h^{-1}xh)\widetilde{f}_2(x) dx dh$$

\noindent Let us fix open and closed $U(W)(F_v)$-invariant neighborhoods $\omega\subset \mathfrak{u}(V)(F_v)$ and $\Omega\subset U(V)(F_v)$ of $0$ and $1$ as in \S \ref{orbital integrals}. Assume that $\widetilde{f}_2$ is supported in $\Omega$. Then, we have (see \S \ref{orbital integrals} for the definition of $f_{i,\natural}$)

$$\displaystyle J(f_1,f_2)=\int_{U(W)(F)}\int_{\mathfrak{u}(V)(F)} f_{1,\natural}(h^{-1}Xh)f_{2,\natural}(X) dX dh$$

\noindent By Fourier transform, we also have

$$\displaystyle J(f_1,f_2)=\int_{U(W)(F)}\int_{\mathfrak{u}(V)(F)} \widecheck{f}_{1,\natural}(h^{-1}Xh)\widehat{f}_{2,\natural}(X) dX dh$$

\noindent By [Zh1, Corollary 4.5] this expression is absolutely convergent so that we can switch the two integrals and we finally get

$$\displaystyle J(f_1,f_2)=\int_{\mathfrak{u}(V)(F)} \widecheck{f}_{1,\natural}(X)O(X,\widehat{f}_{2,\natural}) dX \leqno (3)$$

\noindent Summing up, from (2) and (3) we deduce that

$$\displaystyle \int_{Temp_H(G)} J_\pi(f_1) J_{\pi^\vee}(f_2) d\mu_G(\pi)=\int_{\mathfrak{u}(V)(F)} \widecheck{f}_{1,\natural}(X)O(X,\widehat{f}_{2,\natural}) dX \leqno (4)$$

\noindent for all functions $f_1,f_2\in \mathcal{S}(G(F))$ with $Supp(\widetilde{f}_2)\subseteq \Omega$.

\vspace{2mm}

\noindent We will also need the following formula (see [Beu1, Proposition 8.2.1(iv)]):

$$\displaystyle \widetilde{f}(1)=\int_{Temp_H(G)} J_\pi(f) d\mu_G(\pi) \leqno (5)$$

\noindent for all $f\in \mathcal{S}(G(F))$.

\subsection{Proof of Theorem \ref{theo principal 2}}\label{proof of main theo}

We keep the notations of the previous paragraph. Let $f\in \mathcal{S}(G(F))$. Denote by $C\subseteq Temp_H(G)$ the support of the function $\pi\mapsto J_{\pi}(f)$. It is a compact set and so is its dual $C^\vee$. Let us denote by $\mathcal{Y}$ the image of $\Omega$ in $(H\backslash G/H)(F)=(U(V)/U(W))(F)$ and by $\mathcal{Z}$ the image of the support of $\widecheck{f}_\natural$ in $(\mathfrak{u}(V)/U(W))(F)$. We will denote by the same letters the corresponding subsets in $(H_1'\backslash G'/H_2')(F)$ and $(\mathfrak{s}/GL_n)(F)$ respectively (see \S \ref{correspondence} (1) and (4)). By Corollary \ref{existence good functions}, there exists a function $f'\in \mathcal{S}(G'(F))$ such that

\begin{itemize}
\item $f'$ is supported in the inverse image of $\mathcal{Y}$ in $G'(F)$;

\item The function $Y\in \mathcal{Z}_{rs}\mapsto \eta'(\Delta_-(Y))O(Y,\widehat{f}'_\natural)$ is constant and equal $\eta'(\Delta_-(\xi_-))\mu_{\xi_-}(\widehat{f}'_\natural)$;

\item $\mu_{\xi_-}(\widehat{f}'_\natural)\neq 0$;

\item For all $\Pi\in BC(C^\vee)$ we have

$$\displaystyle I_\Pi(f')=\lv \tau\rv_E^{(d_n+d_{n+1})/2}\omega_{\Pi_n}(\tau)\mu_{\xi_-}(\widehat{f}'_\natural)$$
\end{itemize}

\noindent Let $f_2\in \mathcal{S}(G(F))$ be a function matching $f'$ (whose existence is guaranteed by Theorem \ref{smooth transfer Zhang}). Up to multiplying $f_2$ by the characteristic function of $\mathcal{Y}$ we may assume that $\widetilde{f}_2$ is supported in $\Omega$. By \S \ref{correspondence}(8) the functions $\eta(2)^{n(n+1)/2}f'_\natural$ and $f_{2,\natural}$ match. Hence, by Theorem \ref{transfer and FT} so do $\eta(2)^{n(n+1)/2}\widehat{f}'_\natural$ and $\epsilon(\frac{1}{2},\eta_{E/F},\psi)^{n(n+1)/2}\widehat{f}_{2,\natural}$. Thus, by \ref{truncated local expansion}(1) and \ref{truncated local expansion}(2), for all $X\in \Omega_{rs}$, denoting by $Y\in \mathcal{Z}_{rs}$ the corresponding element, we have

\[\begin{aligned}
\displaystyle \epsilon(\frac{1}{2},\eta_{E/F},\psi)^{n(n+1)/2}O(X,\widehat{f}_{2,\natural}) & =\eta(2)^{n(n+1)/2}\eta(-1)^n\eta'(\Delta_-(Y))O(Y,\widehat{f}'_\natural) \\
 & =\eta(2)^{n(n+1)/2}\eta(-1)^n\eta'(\Delta_-(\xi_-))\mu_{\xi_-}(\widehat{f}'_\natural) \\
 & =\eta'(-2\tau)^{n(n+1)/2}\mu_{\xi_-}(\widehat{f}'_\natural)
\end{aligned}\]

\noindent Consequently, we have

\[\begin{aligned}
\displaystyle (1)\;\;\; \left(\frac{\epsilon(\frac{1}{2},\eta_{E/F},\psi)}{\eta'(-2\tau)}\right)^{n(n+1)/2}\int_{\mathfrak{u}(V)(F)} \widecheck{f}_{\natural}(X)O(X,\widehat{f}_{2,\natural}) dX & =\mu_{\xi_-}(\widehat{f}'_\natural)\int_{\mathfrak{u}(V)(F)} \widecheck{f}_{\natural}(X)dX \\
 & =\mu_{\xi_-}(\widehat{f}'_\natural)f_{\natural}(0)=\mu_{\xi_-}(\widehat{f}'_\natural)\widetilde{f}(1) \\
 & =\mu_{\xi_-}(\widehat{f}'_\natural)\int_{Temp_H(G)} J_\pi(f) d\mu_G(\pi)
\end{aligned}\]

\noindent where the last equality follows from \ref{local TF}(5).

\vspace{2mm}

\noindent On the other hand, by Proposition \ref{prop weak comp}, for all $\pi\in C$ we have

$$\displaystyle J_{\pi^\vee}(f_2)=C(\pi^\vee)I_{BC(\pi^\vee)}(f')=C(\pi^\vee)\lv\tau\rv_E^{(d_n+d_{n+1})/2}\omega_{BC(\pi^\vee)_n}(\tau)\mu_{\xi_-}(\widehat{f}'_\natural)$$

\noindent It follows that

\[\begin{aligned}
\displaystyle (2)\;\;\; \int_{Temp_H(G)} J_\pi(f) J_{\pi^\vee}(f_2)  & d\mu_G(\pi)= \\
 & \mu_{\xi_-}(\widehat{f}'_\natural)\lv\tau\rv_E^{(d_n+d_{n+1})/2}\int_{Temp_H(G)} J_\pi(f)C(\pi^\vee)\omega_{BC(\pi^\vee)_n}(\tau)d\mu_G(\pi)
\end{aligned}\]

\noindent Since $\mu_{\xi_-}(\widehat{f}'_\natural)\neq 0$, we deduce from (1), (2) and \ref{local TF}(4) that

$$\displaystyle \int_{Temp_H(G)} J_\pi(f) \left(\kappa(\pi^\vee) C(\pi^\vee)-1\right)d\mu_G(\pi)=0 \leqno (3)$$

\noindent for all $f\in \mathcal{S}(G(F))$ and where

$$\displaystyle \kappa(\pi)=\lv \tau\rv_E^{(d_n+d_{n+1})/2}\left(\frac{\epsilon(\frac{1}{2},\eta_{E/F},\psi)}{\eta'(-2\tau)} \right)^{n(n+1)/2}\omega_{BC(\pi)_n}(\tau)$$

\noindent Let $\mathcal{Z}(G)$ denotes the Bernstein center of $G$ (see [BD]). We may see $\mathcal{Z}(G)$ as a unital subalgebra of the space of continuous functions on $Temp(G)$ which moreover acts on $\mathcal{S}(G(F))$ with the property that $J_\pi(z\star f)=z(\pi)J_\pi(f)$ for all $z\in \mathcal{Z}(G)$, all $f\in \mathcal{S}(G(F))$ and all $\pi\in Temp(G)$. Thus by (3), we get

$$\displaystyle \int_{Temp_H(G)} z(\pi)J_\pi(f) \left(\kappa(\pi^\vee) C(\pi^\vee)-1\right)d\mu_G(\pi)=0 \leqno (4)$$

\noindent for all $f\in \mathcal{S}(G(F))$ and all $z\in \mathcal{Z}(G)$. For all $\pi\in Temp(G)$ let us denote by $\chi_\pi$ the 'infinitesimal character' of $\pi$, that is the algebra homomorphism $\chi_\pi:\mathcal{Z}(G)\to \mathbb{C}$ given by $\chi_\pi(z):=z(\pi)$ for all $z\in \mathcal{Z}(G)$. Set $Y:=Specmax \mathcal{Z}(G)$. Then, the map $Temp(G)\to Y$, $\pi\mapsto \chi_\pi$ is continuous and proper. Let $Y_{temp}\subseteq Y$ be the image of this map and $\mu_{Y}$ be the push-forward of the Plancherel measure $\mu_G$ to $Y_{temp}$. Then, by the disintegration of measures there exists a measurable mapping $\chi\mapsto \mu_\chi$ from $Y_{temp}$ to the space of measures on $Temp(G)$ such that

$$\displaystyle \int_{Temp(G)} \varphi(\pi)d\mu_G(\pi)=\int_{Y_{temp}} \int_{Temp(G)} \varphi(\pi)d\mu_\chi(\pi)d\mu_Y(\chi) \leqno (5)$$

\noindent for all continuous compactly-supported function $\varphi:Temp(G)\to \mathbb{C}$ and such that for all $\chi\in Y_{temp}$, $\mu_\chi$ is supported on $Temp_\chi(G):=\{\pi\in Temp(G)\mid \chi_\pi=\chi \}$. By (4), we get

$$\displaystyle \int_{Y_{temp}} z(\chi)\int_{Temp_{H,\chi}(G)}J_\pi(f) \left(\kappa(\pi^\vee) C(\pi^\vee)-1\right)d\mu_\chi(\pi) d\mu_Y(\chi)=0 \leqno (6)$$

\noindent for all $f\in \mathcal{S}(G(F))$ and all $z\in \mathcal{Z}(G)$ where we have set $Temp_{H,\chi}(G):=Temp_H(G)\cap Temp_\chi(G)$. Since the restriction of $\mathcal{Z}(G)$ to $Y_{temp}$ is self-adjoint (i.e. for all $z\in \mathcal{Z}(G)$ there exists $z^*\in \mathcal{Z}(G)$ such that $z^*(\chi)=\overline{z(\chi)}$ for all $\chi\in Y_{temp}$), separates points and for all $f\in \mathcal{S}(G(F))$ the function $\pi\in Temp(G)\mapsto J_\pi(f)$ is compactly supported, by (6) and the Stone-Weierstrass theorem for $\mu_Y$-almost all $\chi\in Y_{temp}$, we get

$$\displaystyle \int_{Temp_{H,\chi}(G)}J_\pi(f) \left(\kappa(\pi^\vee) C(\pi^\vee)-1\right)d\mu_\chi(\pi)=0$$

\noindent for all $f\in \mathcal{S}(G(F))$. Since $Temp_\chi(G)$ is finite, we have

$$\displaystyle \int_{Temp_{H,\chi}(G)}J_\pi(f) \left(\kappa(\pi^\vee) C(\pi^\vee)-1\right)d\mu_\chi(\pi)=\sum_{\pi\in Temp_{H,\chi}(G)}
J_\pi(f) \left(\kappa(\pi^\vee) C(\pi^\vee)-1\right)\mu_\chi(\pi)$$

\noindent for $\mu_Y$-almost all $\chi\in Y_{temp}$ and all $f\in \mathcal{S}(G(F))$. Finally, as the spherical characters $J_\pi$ for $\pi\in Temp_{H,\chi}(G)$ are linearly independent, we get that

$$\displaystyle \left(\kappa(\pi^\vee) C(\pi^\vee)-1\right)\mu_\chi(\pi)=0$$

\noindent for $\mu_Y$-almost all $\chi\in Y_{temp}$ and all $\pi\in Temp_{H,\chi}(G)$ which by (5) means that $\kappa(\pi)C(\pi)=1$ for $\mu_G$-almost all $\pi\in Temp_H(G)$. Since $\pi\in Temp_H(G)\mapsto \kappa(\pi)C(\pi)$ is analytic and the support of $\mu_G$ is precisely $Temp(G)$, it follows that $\kappa(\pi)C(\pi)=1$ for all $\pi\in Temp_H(G)$ which is what we wanted. $\blacksquare$

\subsection{A first corollary}\label{a first corollary}

In this paragraph we prove the following corollary to Theorem \ref{theo principal 2}. It will be needed for the proof of Theorem \ref{theo II}.

\begin{cor}\label{spectral char of transfer}
Let $f\in \mathcal{S}(G(F))$ and $f'\in\mathcal{S}(G'(F))$. Then $f$ and $f'$ match if and only if we have

$$\displaystyle I_{BC(\pi)}(f')=\kappa(\pi)J_\pi(f)$$

\noindent for all $\pi\in Temp_H(G)$ and where as before we have set

$$\displaystyle \kappa(\pi)=\lv \tau\rv_E^{(d_n+d_{n+1})/2}\left(\frac{\epsilon(\frac{1}{2},\eta_{E/F},\psi)}{\eta'(-2\tau)} \right)^{n(n+1)/2}\omega_{BC(\pi)_n}(\tau)$$
\end{cor}

\vspace{2mm}

\noindent\ul{Proof}: The necessity follows from Theorem \ref{theo principal 2}. Let us prove the sufficiency. Thus, we assume that

$$\displaystyle I_{BC(\pi)}(f')=\kappa(\pi)J_\pi(f)$$

\noindent for all $\pi\in Temp_H(G)$ and we want to prove that $f$ and $f'$ match. Let $f_2\in \mathcal{S}(G(F))$ be a function which matches $f'$ (such a function exists by Theorem \ref{smooth transfer Zhang}). Then, by Theorem \ref{theo principal 2} and the assumption, for all $\pi\in Temp_H(G)$ we have $J_\pi(f)=J_\pi(f_2)$. Thus, by \ref{local TF}(2), for all $f_1\in \mathcal{S}(G(F))$ we have

$$\displaystyle J(f_1,f)=J(f_1,f_2) \leqno (1)$$

\noindent Let $x_0\in U(V)_{rs}(F)$ and choose $f_1$ so that $\widetilde{f}_1$ is supported in a small neighborhood of $x_0$ in $U(V)_{rs}(F)$. Then a formal manipulation, which is justified since everything is absolutely convergent here, yields

$$\displaystyle J(f_1,f)=\int_{U(V)(F)} f_1(x) O(x,f)dx \leqno (2)$$

\noindent and

$$\displaystyle J(f_1,f_2)=\int_{U(V)(F)} f_1(x) O(x,f_2)dx \leqno (3)$$

\noindent Since the functions $x\in U(V)_{rs}(F)\mapsto O(x,f)$ and $x\in U(V)_{rs}(F)\mapsto O(x,f_2)$ are locally constant (see [Zh1, Proposition 3.13]), we may choose $f_1$ such that $\displaystyle \int_{U(V)(F)} f_1(x) O(x,f)dx=O(x_0,f)$ and $\displaystyle \int_{U(V)(F)} f_1(x) O(x,f_2)dx=O(x_0,f_2)$. For such a choice, it follows from (1), (2) and (3) that $O(x_0,f)=O(x_0,f_2)$. As $x_0$ was arbitrary we see that $f$ and $f_2$ have the same regular semisimple orbital integrals and hence $f$ and $f'$ match. $\blacksquare$

\subsection{Proof of Theorem \ref{theo II}}\label{proof of II theo}

We may assume that $\pi$ is abstractly $H(\mathbb{A})$-distinguished (hence for all $v$, $\pi_v$ is $H_v$-distinguished). By the multiplicity one theorems of Aizenbud-Gourevitch-Rallis-Schiffmann [AGRS] of of Jiang-Sun-Zhu [JSZ], there exists a constant $C$ such that

$$\displaystyle J_\pi(f)=C\prod_v J_{\pi_v}^\natural(f_v)$$

\noindent for all factorizable test function $f=\prod_v f_v\in \mathcal{S}(G(\mathbb{A}))$ and we only need to show that $C=4^{-1}\mathcal{L}(\pi,1/2)$. For this, it is sufficient to prove the existence of $f\in \mathcal{S}(G(\mathbb{A}))$ with

$$\displaystyle J_\pi(f)=4^{-1}\mathcal{L}(\pi,\frac{1}{2})\prod_v J_{\pi_v}^\natural(f_v)$$

\noindent and $J_{\pi_v}^\natural(f_v)\neq 0$ for all place $v$. Let $v_1$ be a (nonarchimedean) place of $F$ such that $BC(\pi_{v_1})$ is supercuspidal (such a place exists by assumption). This implies in particular that $BC(\pi)$ is cuspidal and $\pi_{v_1}$ supercuspidal (by Lemma \ref{lemma supercuspidal}). By Theorem \ref{comparison RTF}, Theorem \ref{theo principal} and identity \ref{transfer of spherical characters} (1), it suffices to show that there exist a nice function $f'\in\mathcal{S}(G'(\A))$ matching a tuple of nice functions $(f^{W'})_{W'}$, $f^{W'}\in \mathcal{S}(G^{W'}(\A))$, such that $I_{BC(\pi_v)}(f'_v)\neq 0$ for all $v$. Let $\Omega_1$ be the Bernstein component of $BC(\pi_{v_1})$ in $G'(F_{v_1})$. Then, we can find a function ${f'}^\circ_{v_1}\in\mathcal{S}(G'(F_{v_1}))_{\Omega_1}$ such that $I_{BC(\pi_{v_1})}({f'}^\circ_{v_1})\neq 0$. Let $f'=\prod_v f'_v$ be a factorizable test function in $\mathcal{S}(G'(\A))$ such that $f'_{v_1}={f'}^\circ_{v_1}$ and $I_{BC(\pi_v)}(f'_v)\neq 0$ for all other place $v$. By construction, the function $f'$ is nice. Moreover, by the assumption on archimedean places, Theorem \ref{smooth transfer Zhang} and Theorem \ref{fundamental lemma}, we can find a tuple of functions $(f^{W'})_{W'}$, $f^{W'}\in \mathcal{S}(G^{W'}(\A))$, matching $f'$. Of course, the functions $f^{W'}$ have no reason of being nice. However, by Lemma \ref{lemma supercuspidal} for all $W'$ there exists a finite union $\Omega^{W'}_1$ of cuspidal Bernstein components of $G^{W'}(F_{v_1})$ such that $\Omega^{W'}_1$ contains all irreducible representation of $G^{W'}(F_{v_1})$ whose base change belongs to $\Omega_1$ and by Corollary \ref{spectral char of transfer} up to replacing $f^{W'}_{v_1}$ by its projection $f^{W'}_{v_1,\Omega^{W'}_1}$ onto $\mathcal{S}(G^{W'}(F_{v_1}))_{\Omega^{W'}_1}$ we may assume that $f^{W'}_{v_1}=f^{W'}_{v_1,\Omega^{W'}_1}$. Then, for all $W'$ the function $f^{W'}$ is nice and we are done. $\blacksquare$

\appendix

\section{Compaison of relative trace formulae}\label{appendix}

The goal of this appendix is to provide a proof of Theorem \ref{comparison RTF}. Inspired by [Kott, \S 18], we start by introducing a convenient notion of norms on the adelic points of a variety over $F$. 

\subsection{Norms on adelic varieties}\label{norms}

\noindent We will use the following convenient although not very precise notations. If $f_1$, $f_2$ are positive valued functions on a set $X$ we will write

$$\displaystyle f_1(x)\ll f_2(x),\;\;\mbox{ for all } x\in X$$

\noindent to mean that there exists a constant $C>0$ such that $f_1(x)\leqslant Cf_2(x)$ for all $x\in X$. We will also write

$$\displaystyle f_1(x)\prec f_2(x),\;\;\mbox{ for all } x\in X$$

\noindent or just $f_1\prec f_2$ if there exist constants $C,d>0$ such that $f_1(x)\leqslant Cf_2(x)^d$ for all $x\in X$. Finally, we will write 

$$\displaystyle f_1(x)\sim f_2(x),\;\;\mbox{ for all } x\in X$$

\noindent or simply $f_1\sim f_2$ if $f_1\prec f_2$ and $f_2\prec f_1$.

\vspace{2mm}

\noindent Let $X$ be a set. By an abstract norm on $X$ we will just mean a function $\lVert .\rVert:X\to [1,+\infty[$. Let $\lVert .\rVert_1$ and $\lVert .\rVert_2$ be two abstract norms on $X$. We will say that $\lVert .\rVert_1$ dominates $\lVert .\rVert_2$ if $\lVert x\rVert_2\prec \lVert x\rVert_1$ for all $x\in X$ and we will say that $\lVert .\rVert_1$ and $\lVert .\rVert_2$ are equivalent if $\lVert .\rVert_1$ dominates $\lVert .\rVert_2$ and $\lVert .\rVert_2$ dominates $\lVert .\rVert_1$ i.e. if $\lVert .\lVert_1\sim \lVert.\rVert_2$. Let $f:X\to Y$ be a map between two sets and let $\lVert .\rVert_Y$ be an abstract norm on $Y$. Then, we define an abstract norm $f^*\lVert .\rVert_Y$ on $X$ by

$$\displaystyle f^*\lVert x\rVert_Y:=\lVert f(x)\rVert_Y$$

\noindent for all $x\in X$.

\vspace{2mm}

\noindent Let $F$ be a number field, $\A$ its ring of adeles and for every place $v$ of $F$ we will denote by $F_v$ the corresponding completion. For every finite extension $F'$ of $F$, we will write $\A_{F'}=\A\otimes_F F'$ for the adele ring of $F'$. We fix algebraic closures $\overline{F}$ of $F$ and $\overline{F}_v$ of $F_v$. For every place $v$ of $F$, we will denote by $\lvert .\rvert_v$ the normalized absolute value on $F_v$. This absolute value extends uniquely to an absolute value on $\overline{F}_v$ that we will also denote by $\lvert .\rvert_v$. We define

$$\displaystyle \A_{\overline{F}}=\overline{F}\otimes_F\A=\varinjlim_{F'} \A_{F'}$$

\noindent where the limits is taken over all finite subextension of $\overline{F}/F$. Let $X$ be an algebraic variety over $F$ (i.e. a reduced separated scheme of finite type over $F$). Since $X$ is of finite type we have $\displaystyle X(\A_{\overline{F}})=\varinjlim_{F'} X(\A_{F'})$. We are going to define certain (equivalence classes of) abstract norms on $X(\A_{\overline{F}})$ and $X(\overline{F}_v)$, $v$ a place of $F$. The definition of these abstract norms in mainly inspired by [Kott, \S 18]. First assume that $X$ is affine and choose a set $\{P_1,\ldots,P_k\}$ of generators for the $F$-algebra $F[X]$. For every place $v$ of $F$ we define an abstract norm $\lVert .\rVert_{X_v}$ on $X(\overline{F}_v)$ by

$$\displaystyle \lVert x\rVert_{X_v}:=\max\left(1,\lvert P_1(x)\rvert_v,\ldots,\lvert P_k(x)\rvert_v\right)$$

\noindent for all $x\in X(\overline{F}_v)$. Choosing a different generating set $\{Q_1,\ldots,Q_\ell\}$ would yield another family of abstract norms $(\lVert.\rVert_{X_v}')_v$ with the following properties:

\begin{itemize}
\item For all place $v$, $\lVert.\rVert_{X_v}'\sim \lVert.\rVert_{X_v}$;

\item There exists $d>0$ such that for almost all place $v$, we have

$$\displaystyle \lVert.\rVert_{X_v}^{1/d}\leqslant \lVert.\rVert_{X_v}'\leqslant \lVert.\rVert_{X_v}^{d}$$
\end{itemize}

\noindent In particular, for all $v$ the equivalence class of the abstract norm $\lVert.\rVert_{X_v}$ does not depend on the particular generating set chosen and by a {\it norm} on $X(\overline{F}_v)$ we will mean any abstract norm in this equivalence class. Note that the norms $(\lVert .\rVert_{X_v})_v$ constructed above are Galois invariant in the sense that $\lVert {}^\sigma x\rVert_{X_v}=\lVert x\rVert_{X_v}$ for all $x\in X(\overline{F}_v)$ and all $\sigma\in Gal(\overline{F}_v/F_v)$. This allows us to extend the norm $\lVert .\rVert_{X_v}$ to $X(K)$ for any finite extension $K$ of $F_v$: choosing any embedding $\iota:K\hookrightarrow \overline{F}_v$ we set

$$\displaystyle \lVert x\rVert_{X_v}:=\lVert \iota(x)\rVert_{X_v}$$

\noindent for any $x\in X(K)$.

\vspace{2mm}

\noindent We now define an abstract norm $\lVert .\rVert_{X}$ on $X(\A_{\overline{F}})$ as follows. Let $x\in X(\A_{\overline{F}})$ and choose a finite extension $F'/F$ such that $x\in X(\A_{F'})$. Then, we may write $x$ as a product $\prod_w x_w$, $x_w\in X(F'_w)$, indexed by the set of places of $F'$ and we set

$$\displaystyle \lVert x\rVert_{X}:=\prod_v \left(\prod_{w\mid v}\lVert x_w\rVert_{X_v}^{[F'_w:F_v]}\right)^{1/[F':F]}$$

\noindent where the first product is over the set of places $v$ of $F$ and the second product is over the set of places $w$ of $F'$ above $v$. Note that this definition does not depend on the choice of the finite extension $F'/F$ such that $x\in X(\A_{F'})$. Moreover, choosing a different generating set would give an equivalent abstract norm. By a {\it norm} on $X(\A_{\overline{F}})$ we will mean any abstract norm in this equivalence class. We will assume from now on that for any affine variety $X$ over $F$ we have fixed norms $\lVert .\rVert_X$ on $X(\A_{\overline{F}})$ and norms $\lVert.\rVert_{X_v}$ on $X(\overline{F}_v)$, for all place $v$ of $F$, as above (i.e. by choosing a finite generating set of $F[X]$). In the particular case $X=\A^1$ (the affine line) we will even take

$$\displaystyle \lVert x\rVert_{\A^1_v}= \max(1,\lvert x\rvert_v)$$

\noindent for all place $v$ of $F$ and for all $x\in X(\overline{F}_v)=\overline{F}_v$. Note that by the product formula we then have

$$\lVert x\rVert_{\A^1}=\lVert x^{-1}\rVert_{\A^1} \leqno (1)$$

\noindent for every $x\in \overline{F}^\times$.

\vspace{2mm}

\noindent We continue to assume that $X$ is affine. Let $\mathcal{U}=(U_i)_{i\in I}$ be a finite covering of $X$ by affine open subsets. We can define another abstract norm $\lVert.\rVert_{X_v,\mathcal{U}}$ on $X(\overline{F}_v)$ by

$$\displaystyle \lVert x\rVert_{X_v,\mathcal{U}}:=\min\{\lVert x\rVert_{U_{i,v}};i\in I\mbox{ such that }x\in U_i(\overline{F}_v)\},\;\;\; x\in X(\overline{F}_v)$$

\noindent Then we have (see [Kott, Proposition 18.1(6)])

\begin{itemize}
\item For all place $v$, $\lVert .\rVert_{X_v,\mathcal{U}}\sim \lVert .\rVert_{X_v}$;

\item There exists $d>0$ such that for almost all place $v$, we have

$$\displaystyle \lVert .\rVert_{X_v}^{1/d}\leqslant \lVert .\rVert_{X_v,\mathcal{U}}\leqslant \lVert .\rVert_{X_v}^{d}$$

\item For all place $v$, $\lVert .\rVert_{X_v,\mathcal{U}}$ is Galois invariant.
\end{itemize}

\noindent We can also define an abstract norm $\lVert .\rVert_{X,\mathcal{U}}$ on $X(\A_{\overline{F}})$ by sending $x\in X(\A_{F'})$, $F'/F$ a finite extension, to

$$\displaystyle \lVert x\rVert_{X,\mathcal{U}}:=\prod_v \left(\prod_{w\mid v}\lVert x_w\rVert_{X_v,\mathcal{U}}^{[F'_w:F_v]}\right)^{1/[F':F]}$$

\noindent Then $\lVert .\rVert_{X,\mathcal{U}}$ is a norm on $X(\A_{\overline{F}})$ (i.e. $\lVert .\rVert_{X,\mathcal{U}}\sim\lVert .\rVert_{X}$). This allows us to extend the definition of the abstract norms $\lVert .\rVert_X$ and $\lVert .\rVert_{X_v}$ to any algebraic variety $X$ over $F$ as follows. Let $X$ be such a variety and choose a finite covering $\mathcal{U}=(U_i)_{i\in I}$ of $X$ by affine open subsets. Then the definitions of the abstract norms $\lVert .\rVert_{X,\mathcal{U}}$ and $\lVert .\rVert_{X_v,\mathcal{U}}$ as above still make sense and we will set $\lVert .\rVert_{X}:=\lVert .\rVert_{X,\mathcal{U}}$, $\lVert .\rVert_{X_v}:=\lVert .\rVert_{X_v,\mathcal{U}}$. Choosing a different covering $\mathcal{V}$ of $X$ would give abstract norms $(\lVert .\rVert'_{X_v})_v$ and $\lVert .\rVert_{X}'$ satisfying the following

\begin{itemize}
\item For all $v$, $\lVert .\rVert'_{X_v}\sim \lVert .\rVert_{X_v}$ and there exists $d>0$ such that for almost all $v$ we have

$$\lVert .\rVert_{X_v}^{1/d}\leqslant \lVert .\rVert'_{X_v}\leqslant \lVert .\rVert_{X_v}^{d}$$

\item $\lVert .\rVert_{X}'\sim \lVert .\rVert_{X}$.
\end{itemize}

\noindent In particular the equivalence class of $\lVert .\rVert_X$ (resp. of $\lVert .\rVert_{X_v}$ for $v$ a place of $F$) doesn't depend on the particular choice of $\mathcal{U}$ and by a {\it norm} on $X(\A_{\overline{F}})$ (resp. on $X(\overline{F}_v)$) we will mean any abstract norm in this equivalence class. From now on we assume that every algebraic variety over $F$ has been equipped with a family of norms as above (i.e. by choosing a finite covering $\mathcal{U}$ by affine open subsets). If $X$ is affine we also assume that these norms have been defined using the trivial covering $\mathcal{U}=\{X\}$ so that they coincide with the ones we already fixed. If $G$ is an affine algebraic group over $F$ we also define a norm $\lV .\rVert_{[G]}$ on $[G]=G(F)\backslash G(\A)$ by

$$\displaystyle \lV x\rVert_{[G]}:=\inf_{\gamma\in G(F)} \lV \gamma x\rVert_G$$

\noindent for all $x\in [G]$.

\begin{prop}\label{prop norms}
Let $X$ and $Y$ be algebraic varieties over $F$ and let $G$ be an affine algebraic group over $F$.
\begin{enumerate}[(i)]
\item The function $x\mapsto \lVert x\rVert_{X}$ is locally bounded on $X(\A)$.

\item Let $f:X\to Y$ be a morphism of algebraic varieties. Then $f^*\lVert .\rVert_{Y}\prec \lVert .\rVert_{X}$. In particular we have $\lVert gg'\rVert_G\prec \lVert g\rVert_G\lVert g'\rVert_G$ and $\lVert g^{-1}\rVert_G\sim \lVert g\rVert_G$ for all $g,g'\in G(\A_{\overline{F}})$. If moreover $f$ is a finite morphism (in particular if it is a closed embedding) then $f^*\lVert .\rVert_{Y}\sim \lVert .\rVert_{X}$.

\item Let $f\in F[X]$ and let $X_f=D(f)$ be the principal open subset of $X$ defined by the nonvanishing of $f$. Then, we have

$$\displaystyle \lVert x\rVert_{X_f}\sim \lVert x\rVert_X \lVert f(x)^{-1}\rVert_{\A^1}$$

\noindent for all $x\in X_f(\A_{\overline{F}})$.

\item Let $U\subset X$ be an open subset and assume that $X$ is quasi-affine. Then, we have

$$\displaystyle \lVert x\rVert_U\sim \lVert x\rVert_X$$

\noindent for all $x\in U(\overline{F})$. More generally if $p:X\to Y$ is a regular map and $Y$ is quasi-affine then for all open subset $V\subset Y$ we have

$$\displaystyle \lVert x\rVert_{p^{-1}(V)}\sim \lVert x\rVert_X$$

\noindent for all $x\in p^{-1}(V)(\A_{\overline{F}})$ such that $p(x)\in V(\overline{F})$.

\item If $X$ is quasi-affine then there exists $d>0$ such that

$$\displaystyle \sum_{x\in X(F)} \lVert x\rVert_{X}^{-d}$$

\noindent converges.

\item Let $d_rg$ be a right Haar measure on $G(\A)$. Then there exists $d>0$ such that the two integrals

$$\displaystyle \int_{G(\A)} \lVert g\rVert_G^{-d} d_rg,\;\;\; \int_{[G]} \lVert x\rVert_{[G]}^{-d} dx$$

\noindent converge.

\item Assume that $X$ carries a $G$-action and that we have a regular map $p:X\to Y$ making $X$ into a $G$-torsor over $Y$. Fix a right Haar measure $d_rg$ on $G(\A)$. Then for all $d>0$ there exists $d'>0$ such that

$$\displaystyle \int_{G(\A)} \lVert gx\rVert_X^{-d'} d_rg\ll \lVert p(x)\rVert_Y^{-d}$$

\noindent for all $x\in X(\A)$.

\item Assume that $G$ is connected and reductive and let $\mathfrak{S}\subset G(\A)$ be a Siegel domain (see [MW1, \S I.2.1]). Then, we have

$$\displaystyle \lV g\rVert_G\sim \lV g\rVert_{[G]}$$

\noindent for all $g\in \mathfrak{S}$.

\item Let $H<G$ be a closed subgroup such that $G/H$ is quasi-affine (this is the case if for example $H$ is reductive or if there is no nontrivial morphism $H\to \G_m$). Then, we have

$$\displaystyle \lV x\rVert_{[H]}\sim \lV x\rVert_{[G]}$$

\noindent for all $x\in [H]$. In particular, by (vi) there exists $d>0$ such that the integral

$$\displaystyle \int_{[H]} \lV x\rV_{[G]}^{-d} dx$$

\noindent converges.
\end{enumerate}
\end{prop}

\noindent\ul{Proof}:

\begin{enumerate}[(i)]
\item This follows from the fact that for all $v$ the function $x\in X(F_v)\mapsto \lVert x\rVert_{X_v}$ is locally bounded and the fact that for almost all $v$  we have $\lVert x_v\rVert_{X_v}=1$ for all $x_v\in X(\mathcal{O}_v)$.

\item It suffices to prove the following

\begin{itemize}
\item For all place $v$, we have $f^*\lVert .\rVert_{Y_v}\prec \lVert .\rVert_{X_v}$ and if $f$ is finite $\lVert .\rVert_{X_v}\prec f^*\lVert .\rVert_{Y_v}$;

\item There exists $d>0$ such that for almost all place $v$, we have $f^*\lVert .\rVert_{Y_v}\leqslant \lVert .\rVert_{X_v}^{d}$ and if $f$ is finite $\lVert .\rVert_{X_v}^{1/d}\leqslant f^*\lVert .\rVert_{Y_v}$.
\end{itemize}

\noindent Assume that the norms $(\lVert .\rVert_{X_v})_v$ have been defined using the finite affine open covering $\mathcal{U}=(U_i)_{i\in I}$ of $X$ and that the norms $(\lVert .\rVert_{Y_v})_v$ have been defined using the finite affine open covering $\mathcal{V}=(V_j)_{j\in J}$ of $Y$. Up to refining $\mathcal{U}$, we may assume that for all $j\in J$ there exists a subset $I(j)\subset I$ such that $f^{-1}(V_j)=\bigcup_{i\in I(j)} U_i$. If moreover $f$ is finite then for all $j\in J$, the open subset $f^{-1}(V_j)$ is affine so that we may assume that $\mathcal{U}=(f^{-1}(V_j))_{j\in J}$. This allows us to reduce to the case where both $X$ and $Y$ are affine in which case the statement can be proved much the same way as [Kott, Proposition 18.1(1)].

\item Assume that the family of norms $(\lVert .\rVert_{X_v})_v$ has been defined using the finite affine open covering $\mathcal{U}=(U_i)_{i\in I}$ of $X$. Set $U_{i,f}=U_i\cap X_f$ for all $i\in I$. Obviously, we may assume that the family of norms $(\lVert .\rVert_{X_{f,v}})_v$ has been defined using the affine open covering $\mathcal{U}_f=(U_{i,f})_{i\in I}$ of $X_f$ and that

$$\displaystyle \lVert x\rVert_{U_{i,f,v}}=\max\left( \lVert x\rVert_{U_{i,v}}, \lvert f(x)\rvert_v^{-1}\right)$$

\noindent for all place $v$ of $F$ and all $x\in U_{i,f}(\overline{F}_v)$. Then we have

$$\displaystyle \sqrt{\lVert x\rVert_{U_{i,v}} \max(1,\lvert f(x)\rvert_v^{-1})}\leqslant \lVert x\rVert_{U_{i,f,v}}\leqslant \lVert x\rVert_{U_{i,v}} \max(1,\lvert f(x)\rvert_v^{-1})$$

\noindent for all place $v$ of $F$ and all $x\in U_{i,f}(\overline{F}_v)$. It follows that

$$\displaystyle \sqrt{\lVert x\rVert_{X_v} \max(1,\lvert f(x)\rvert_v^{-1})}\leqslant\lVert x\rVert_{X_{f,v}}\leqslant \lVert x\rVert_{X_v} \max(1,\lvert f(x)\rvert_v^{-1})$$

\noindent for all place $v$ of $F$ and all $x\in X_f(\overline{F}_v)$. Taking the product we get

$$\displaystyle \sqrt{\lVert x\rVert_{X} \lVert f(x)^{-1}\rVert_{\A^1}}\leqslant\lVert x\rVert_{X_f}\leqslant \lVert x\rVert_{X} \lVert f(x)^{-1}\rVert_{\A^1}$$

\noindent for all $x\in X_f(\A_{\overline{F}})$.

\item We prove the second claim which is more general than the first. Let $p:X\to Y$ be a regular map, $V\subset Y$ an open subset and assume that $Y$ is quasi-affine. It already follows from (ii) that we have

$$\displaystyle \lVert x\rVert_X\prec \lVert x\rVert_{p^{-1}(V)}$$

\noindent for all $x\in X(\A_{\overline{F}})$. Hence, it suffices to prove the reverse inequality for all $x\in p^{-1}(V)(\A_{\overline{F}})$ such that $p(x)\in V(\overline{F})$. As $Y$ is quasi-affine, up to replacing $V$ by a finite affine open cover we may assume that $V=Y_f$ for some $f\in F[Y]$. Still denoting by $f$ its image in $F[X]$ we then have $p^{-1}(V)=X_f$. Then by (ii), (iii) and (1) we have

$$\displaystyle \lVert x\rVert_{X_f}\sim \lVert x\rVert_{X}\lVert f(x)^{-1}\rVert_{\A^1}=\lVert x\rVert_{X}\lVert f(x)\rVert_{\A^1}\prec \lVert x\rVert_{X}$$

\noindent for all $x\in X_f(\A_{\overline{F}})$ such that $f(x)\in \overline{F}^{\times}$. This implies the desired inequality.

\item As there exists an open embedding of $X$ into an affine variety, by (iv) we immediately reduce to the case where $X$ itself is affine. Then, we can find a closed embedding $\iota:X\hookrightarrow \A^n$ for some integer $n>0$ and by (ii) we are reduced to prove the statement for $X=\A^n$ and then eventually for $X=\A^1$ in which case the statement is easily checked.

\item  Note that

\[\begin{aligned}
\displaystyle \int_{[G]} \lV x\rVert_{[G]}^{-d} dx \leqslant \int_{[G]} \sum_{\gamma\in G(F)}\lV \gamma x\rVert_{G}^{-d} dx =\int_{G(\A)} \lV g\rVert_G^{-d} d_rg
\end{aligned}\]

\noindent for all $d>0$. Hence it suffices to show that for $d$ sufficiently large the last integral above is convergent.  Assume that $H$ is a closed distinguished subgroup of $G$ isomorphic to $\G_m$ or $\G_a$. We first show that if the statement is true for both $H$ and $G/H$ then it is true for $G$. For this we write

$$\displaystyle \int_{G(\A)} \lV g\rV_G^{-d} d_rg=\int_{(G/H)(\A)}\int_{H(\A)} \lV \dot{g}h\rV_G^{-d}d_rhd_r\dot{g}$$

\noindent for all $d>0$ and where $d_rh$, $d_r\dot{g}$ are suitable right Haar measures on $H(\A)$ and $(G/H)(\A)$ respectively (Note that $(G/H)(\A)=G(\A)/H(\A)$). Let $d_0,d_1>0$. Setting $d=d_0+d_1$, we get

$$\displaystyle \int_{H(\A)} \lV gh\rV_G^{-d} d_rh\leqslant \left(\inf_{h\in H(\A)}\lV gh\rV_G\right)^{-d_0}\int_{H(\A)}\lV gh\rV_G^{-d_1} d_rh$$

\noindent for all $g\in G(\A)$. By (ii), there exists $c>0$ such that $\lV h\rV_H\ll \lV gh\rV_G^c\lV g\rV_G^c$ for all $(h,g)\in H(\A)\times G(\A)$. Hence,

$$\displaystyle \int_{H(\A)}\lV gh\rV_G^{-d} d_rh\ll \left(\inf_{h\in H(\A)}\lV gh\rV_G\right)^{-d_0}\lV g\rV_G^{d_1}\int_{H(\A)}\lV h\rV_H^{-d_1/c} d_rh$$

\noindent for all $g\in G(\A)$. As the left hand side above is, as a function of $g$, invariant by right translation by $H(\A)$ we also get

$$\displaystyle \int_{H(\A)}\lV gh\rV_G^{-d} d_rh\ll \left(\inf_{h\in H(\A)}\lV gh\rV_G\right)^{d_1-d_0}\int_{H(\A)}\lV h\rV_H^{-d_1/c} d_rh$$

\noindent for all $g\in G(\A)$. By assumption for $d_1$ sufficiently large the last integral above is convergent. Thus, it only remains to show that for $d'>0$ sufficiently large the integral

$$\displaystyle \int_{(G/H)(\A)} \left(\inf_{h\in H(\A)}\lV \dot{g}h\rV_G\right)^{-d'}d_r\dot{g}$$

\noindent converges. By (ii), we have $\lV \dot{g}\rV_{G/H}\prec \inf_{h\in H(\A)}\lV \dot{g}h\rV_G$ for all $\dot{g}\in G(\A)/H(\A)$. Consequently, the convergence of the last integral above for $d'$ sufficiently large follows from the assumption on $G/H$.

\vspace{2mm}

\noindent Let $P_0$ be a minimal parabolic subgroup of $G$ over $F$. Then, by the Iwasawa decomposition there exists a compact subgroup $K\subset G(\A)$ such that $G(\A)=P_0(\A)K$. As $K$ is compact, by (i) the norm $\lVert .\rVert_{G}$ is bounded on $K$. Moreover, we have

$$\displaystyle \int_{G(\A)} \lVert g\rVert_G^{-d}d_rg=\int_{P_0(\A)} \int_K \lVert kp_0\rVert_G^{-d} dkd_rp_0$$

\noindent for suitable (right) Haar measures $d_rp_0$ and $dk$ on $P_0(\A)$ and $K$ respectively. By (i) and (ii), it follows that we may assume $G=P_0$. Let $P_0=M_0N_0$ be a Levi decomposition. Then the Haar measure $d_rp_0$ decomposes as $d_rp_0=dn_0dm_0$ according to the decomposition $P_0(\A)=N_0(\A)M_0(\A)$. Moreover we have $\lVert n_0m_0\rVert_{P_0}\sim \lVert m_0\rVert_{M_0}\lVert n_0\rVert_{N_0}$ for all $(m_0,n_0)\in M_0(\A)\times N_0(\A)$. This allows us to reduce to the case where $G=M_0$ or $G=N_0$. If $G=N_0$ then it admits a composition series whose successive quotients are isomorphic to $\G_a$ and we are reduced to the case $G=\G_a$ where the statement can be checked directly. Assume now that $G=M_0$ and denote by $A_0$ the maximal split torus in the center of $G$. Then $A_0$ is isomorphic to a product of $\G_m$ and $M_0/A_0$ is anisotropic. Thus, we only need to treat the cases $G=\G_m$ and $G$ anisotropic. Once again if $G=\G_m$ the statement can be checked directly. Now if $G$ is anisotropic we write

$$\displaystyle \int_{G(\A)} \lVert g\rVert_G^{-d} dg=\int_{G(F)\backslash G(\A)} \sum_{\gamma\in G(F)} \lVert \gamma g\rVert_G^{-d} dg \leqno (2)$$

\noindent By (i), (ii) and (v) if $d$ is sufficiently large the function

$$\displaystyle g\in G(\A)\mapsto \sum_{\gamma\in G(F)}\lVert \gamma g\rVert_G^{-d}$$

\noindent is locally bounded. Moreover by [BHC] the quotient $G(F)\backslash G(\A)$ is compact. The result then follows from (2).

\item Let $d>0$. As $p$ is a $G$-torsor and $Y$ is separated, the action of $G$ on $X$ is free i.e. the regular map

$$\displaystyle G\times X\to X\times X$$
$$\displaystyle (g,x)\mapsto (gx,x)$$

\noindent is a closed embedding. By (ii), it follows that there exists $c>0$ such that $\lV g\rV_G\ll \lV gx\rV_X^c\lV x\rV_X^c$ for all $(g,x)\in G(\A)\times X(\A)$. Let $d_0,d_1>0$. Using the same trick as in the first part of the proof of (vi), we show that for $d'=d_0+d_1$ we have

$$\displaystyle \int_{G(\A)}\lV gx\rV_X^{-d'} d_rg\ll \left(\inf_{g\in G(\A)} \lV gx\rV_X \right)^{d_1-d_0}\int_{G(\A)} \lV g\rV_G^{-d_1/c} d_rg$$

\noindent for all $x\in X(\A)$. By (vi), the last integral above is convergent for $d_1$ sufficiently large. Moreover, by (ii) we have $\lV p(x)\rV_{Y}\prec \inf_{g\in G(\A)} \lV gx\rV_X$ for all $x\in X(\A)$. Thus, the statement follows by choosing $d_0$ sufficiently large (depending on $d_1$).

\item Let $T_0$ be a maximal split torus in $R_{F/\Q} G$. Then, up to conjugating $\mathfrak{S}$ by an element of $G(F)$, there exists a compact subset $\Omega\subset G(\A)$ such that

$$\displaystyle \mathfrak{S}\subseteq T_0(\R)\Omega$$

\noindent Hence, by (i) and (ii) it is sufficient to show that

$$\displaystyle \lV a\rVert_G\sim\lV a\rVert_{[G]}$$

\noindent for all $a\in T_0(\R)$. The inequality $\lV a\rVert_{[G]}\leqslant \lV a\rVert_{G}$ is obvious so that we only need to show that $\lV a\rVert_{G}\prec \lV a\rVert_{[G]}$ for all $a\in T_0(\R)$. Let $\chi_1,\ldots,\chi_n$ be a basis of $X^*(T_0)$. Then we have

$$\displaystyle \lV a\rVert_G\sim \max\left(\lv \chi_1(a)\rvert,\lv \chi_1(a)\rvert^{-1},\ldots,\lv \chi_n(a)\rvert,\lv \chi_n(a)\rvert^{-1}\right)$$

\noindent for all $a\in T_0(\R)$. Thus, it suffices to show that for all character $\chi\in X^*(T_0)$ we have $\lv \chi(a)\rvert \prec \lV a\rVert_{[G]}$ for all $a\in T_0(\R)$. Let $\chi$ be such a character and let $V$ be a rational representation of $R_{F/\Q}G$ containing a nonzero vector $v_0$ such that $a.v_0=\chi(a)v_0$ for all $a\in T_0$. Let $v_1,\ldots,v_r$ be a basis of $V$. Set $V_{\A}=V\otimes_{\Q}\A$ and define a nonnegative function $\lv .\rvert_V$ on $V_{\A}$ by

$$\displaystyle \lv \lambda_1v_1+\ldots+\lambda_rv_r\rvert_V=\prod_v \max(\lv \lambda_{1,v}\rvert_v,\ldots,\lv \lambda_{r,v}\rvert_v)$$

\noindent for all $\lambda_1,\ldots,\lambda_r\in \A$. Note that there exist nonzero vectors $v\in V_{\A}$ such that $\lv v\rvert_V=0$ but that, however, if $v\in V_F=V\otimes_{\Q}F$ is nonzero then $\lv v\rvert_V\geqslant 1$. We have $\lv v\rvert_V\prec \lV v\rV_{V_F}$ for all $v\in V_{\A}$, where $V_F$ is considered as an algebraic variety over $F$. Note that $G$ acts on $V_F$ via the natural embedding $G\hookrightarrow (R_{F/\Q}G)_F$. Hence, by (ii) we have

$$\displaystyle \lv \chi(a)\rvert^d\leqslant \lv \chi(a)\rvert^d \lv \gamma v_0\rvert_V=\lv \gamma av_0\rvert_V\prec \lV \gamma av_0\rVert_{V_F}\prec \lV \gamma a\rVert_G$$

\noindent for all $a\in T_0(\R)$, $\gamma\in G(F)$ and where we have set $d=[F:\Q]$. Taking the infimum over $\gamma$ yields the desired inequality.

\item By (ii), the inequality $\lV x\rVert_{[G]}\prec \lV x\rVert_{[H]}$ is obvious so that we only need to show that $\lV x\rVert_{[H]}\prec \lV x\rVert_{[G]}$ for all $x\in [H]$. We will need the following fact (which is where the assumption $G/H$ quasi-affine is crucial):

\vspace{4mm}

\hspace{5mm} (3) There exists a (set-theoretic) section $s:\left(H\backslash G\right)(\overline{F})\to G(\overline{F})$ such that

\vspace{0.4pt}

\hspace{12mm} $\lVert s(x)\rVert_G\prec \lV x\rVert_{H\backslash G}$ for all $x\in \left(H\backslash G\right)(\overline{F})$.

\vspace{4mm}

\noindent Proof of (3): Let $p:G\to H\backslash G$ be the natural surjection. Since $H\backslash G$ is quasi-affine, by (iv), it suffices to find an open covering $(U_i)_{i\in I}$ of $H\backslash G$ and sections $s_i:U_i(\overline{F})\to p^{-1}(U_i)(\overline{F})$ such that $\lVert s_i(x)\rVert_{p^{-1}(U_i)}\prec \lV x\rVert_{U_i}$ for all $i\in I$ and all $x\in U_i(\overline{F})$. It is even sufficient to construct one non-empty open subset $U\subseteq H\backslash G$ and a section $s_U:U(\overline{F})\to p^{-1}(U)(\overline{F})$ such that $\lVert s_U(x)\rVert_{p^{-1}(U)}\prec \lV x\rVert_{U}$ for all $x\in U(\overline{F})$. Indeed, if such a pair $(U,s_U)$ exists, we can find a finite number of translates $U_i:=U\gamma_i$, $\gamma_i\in G(\overline{F})$, $i\in I$, covering $H\backslash G$ and then the sections $s_i:U_i(\overline{F})\to p^{-1}(U_i)(\overline{F})$ given by $s_i(x):=s(x\gamma_i^{-1})\gamma_i$, for all $i\in I$ and $x\in U_i(\overline{F})$, satisfy the desired condition. As $p:G\to H\backslash G$ is a torsor for the \'etale topology, we can find a non-empty open subset $U\subseteq H\backslash G$ and a finite \'etale map $U'\to U$ such that $U'\times_U G$ is the trivial $G$-torsor over $U'$. In particular there exists a regular section $s_{U'}:U'\to U'\times_U G$. Let $s_0:U(\overline{F})\to U'(\overline{F})$ be any set-theoretic section. Then, by (ii) and since $U'\to U$ is finite, the section $s_U:=pr_2\circ s_{U'}\circ s_0:U(\overline{F})\to p^{-1}(U)(\overline{F})$, where $pr_2$ denotes the projection $U'\times_U G\to G$, satisfies the desired condition.

\vspace{2mm}

\noindent Let $s:(G/H)(\overline{F})\to G(\overline{F})$ be a section as in (3). We have $\lV \gamma \rVert_{G/H}\prec\lV \gamma h\rVert_G$, for all $(\gamma,h)\in G(F)\times H(\A)$ (by (ii)) and thus

$$\displaystyle \inf_{\gamma'\in H(\overline{F})} \lV \gamma' h\rVert_H\leqslant \lV s(\gamma)^{-1}\gamma h\rVert_H\prec \lV s(\gamma)\rVert_G\lVert\gamma h\rVert_G\prec \lV \gamma\rVert_{G/H} \lV \gamma h\rVert_G\prec \lV \gamma h\rVert_G$$

\noindent for all $(\gamma,h)\in G(F)\times H(\A)$. Taking the infimum over $\gamma$ it follows that

$$\displaystyle \inf_{\gamma'\in H(\overline{F})} \lV \gamma' h\rVert_H\prec \inf_{\gamma\in G(F)}\lV \gamma h\rVert_G=\lV h\rVert_{[G]}$$

\noindent for all $h\in H(\A)$. Hence, it suffices to show

$$\displaystyle(4)\;\;\;  \lV h\rVert_{[H]}\prec \inf_{\gamma\in H(\overline{F})} \lV \gamma h\rVert_H$$

\noindent for all $h\in H(\A)$. Denote by $N_H$ the unipotent radical of $H$ and let $L_H$ be a Levi component of $H$ (so that $H=L_H\ltimes N_H$). As $[N_H]$ is compact we are easily infer from (i) and (ii) that

$$\displaystyle \lVert \ell n\rVert_{[H]}\sim \lVert \ell\rVert_{[L_H]} \mbox{ and } \inf_{\gamma\in H(\overline{F})}\lVert \gamma \ell n\rVert_{H}\sim \inf_{\gamma_L\in L_H(\overline{F})}\lVert \gamma_L \ell \rVert_{L_H}$$

\noindent for all $\ell\in L_H(\A)$ and all $n\in N_H(\A)$. We are thus reduced to prove (4) in the case where $H$ is reductive. Denote by $H^0$ the connected component of the identity in $H$. Since $H(F)/H^0(F)$, $H(\overline{F})/H^0(\overline{F})$ are finite and $H(\A)/H^0(\A)$ is compact we may assume that $H=H^0$. Let $T_0$ be a maximal split torus of $R_{F/\Q}H$ and let $\chi\in X^*(T_0)$. By (viii), it is sufficient to show that

$$\displaystyle (5)\;\;\; \lv \chi(a)\rvert\prec \lV \gamma a\rVert_H$$

\noindent for all $a\in T_0(\R)$ and all $\gamma\in H(\overline{F})$. Let $V$ be a rational representation of $R_{F/\Q}H$ containing a nonzero vector $v_0$ such that $a.v_0=\chi(a)v_0$ for all $a\in T_0$. Fix a basis $v_1,\ldots,v_n$ of $V$ and let $\lv .\rvert_V$ be the nonnegative function on $V_{\AF}=V\otimes_{\Q} \AF$ defined by

$$\displaystyle \lv \lambda_1v_1+\ldots+\lambda_rv_r\rvert_V=\prod_v \left(\prod_{w\mid v}\max(\lv \lambda_{1,w}\rvert_v,\ldots,\lv \lambda_{r,w}\rvert_v)^{[F'_w:F_v]}\right)^{1/[F':F]}$$

\noindent for all $\lambda_1,\ldots,\lambda_r\in \A_{F'}$, $F'/F$ a finite extension. Note that $\lv v\rv_V\geqslant 1$ for all nonzero vector $v\in V_{\overline{F}}=V\otimes_{\Q}\oF$ and $\lv v\rv_{V}\prec \lV v\rV_{V_F}$ for all $v\in V_{\AF}$. It follows that

$$\displaystyle \lv \chi(a)\rv^d\leqslant \lv \chi(a)\rv^d\lv \gamma v_0\rv_V=\lv \gamma av_0\rv_V\prec \lV \gamma av_0\rV_{V_F}\prec \lV \gamma a\rV_H$$

\noindent for all $(a,\gamma)\in T_0(\R)\times H(\oF)$. Taking the infimum over $\gamma$ we get (5) and this ends the proof of (ix). $\blacksquare$

\end{enumerate}

\vspace{2mm}

\noindent Let $G$ be a connected reductive group over $F$. Fix a maximal compact subgroup $K_\infty$ of $G(\A_\infty)$ and a Haar measure $dg$ on $G(\A)$. We will denote by $\mathcal{U}(\mathfrak{g}_\infty)$ the universal enveloping algebra of (the complexification of) the Lie algebra of $G(\A_\infty)$. For simplicity we will assume that the split center of $G$ is trivial. Denote by $\mathcal{A}([G])$ the space of automorphic functions on $[G]$ by which we mean functions $\phi:[G]\to\C$ satisfying the following conditions

\begin{itemize}
\item $\phi$ is smooth: there exists a compact-open subgroup $K$ of $G(\A_f)$ such that $\phi$ is right $K$-invariant and for all $g_f\in G(\A_f)$ the function $g_\infty\in G(\A_\infty)\mapsto\phi(g_\infty g_f)$ is $C^\infty$;

\item $\phi$ is uniformly of moderate growth: there exists $d>0$ such that for all $u\in \mathcal{U}(\mathfrak{g}_\infty)$ we have $\lv(R(u)\phi)(g)\rv\ll \lV g\rV_{[G]}^d$ for all $g\in G(\A)$.
\end{itemize}

\noindent Note that we don't impose any condition of $K_\infty$-finiteness or $\mathfrak{z}_\infty$-finiteness (where $\mathfrak{z}_\infty$ denotes the center of $\mathcal{U}(\mathfrak{g}_\infty)$). The space $\mathcal{A}([G])$ is naturally equipped with a LF topology (see [Beu1, appendix A] for basic facts about LF vector spaces). As usual, we define $\mathcal{A}_{cusp}([G])$ to be the subspace of cuspidal functions in the following sense: $\phi\in \mathcal{A}([G])$ is cuspidal if for all proper parabolic subgroup $P=MN$ of $G$ we have

$$\displaystyle \int_{[N]}\phi(ng)dn=0$$

\noindent for all $g\in G(\A)$. The space $\mathcal{A}_{cusp}([G])$ is a closed subspace of $\mathcal{A}([G])$ from which it inherits a LF topology and moreover every cuspidal function $\phi\in \mathcal{A}_{cusp}([G])$ is of rapid decay in the following sense: for all $u\in \mathcal{U}(\mathfrak{g}_\infty)$ and for all $d>0$ we have

$$\displaystyle \lv (R(u)\phi)(g)\rv\ll \lV g\rV_{[G]}^{-d}$$

\noindent for all $g\in [G]$ (see [MW1, Corollary I.2.12]). By the open mapping theorem, for all compact-open subgroup $K$ of $G(\A_f)$ the topology on $\mathcal{A}_{cusp}([G])^K$ is also induced by the family of seminorms

$$\displaystyle \lV \phi\rV_{d,u}=\sup_{g\in [G]}\lv (R(u)\phi)(g)\rv \lV g\rV_{[G]}^{d},\;\;\; d>0,u\in \mathcal{U}(\mathfrak{g}_\infty)$$

\noindent There is another natural family of seminorms inducing the given topology on $\mathcal{A}_{cusp}([G])^K$. Let $C_G\in \mathcal{U}(\mathfrak{g}_\infty)$ and $C_K\in \mathcal{U}(\mathfrak{k}_\infty)$ denote the Casimir elements of $G(\A_\infty)$ and $K_\infty$ respectively and set $\Delta=C_G^2+C_K^2$. Then the family of Sobolev seminorms

$$\displaystyle \lV \phi\rV_k=\lV R(1+\Delta)^k \phi\rV_{L^2([G])},\;\;\; k\geqslant 0,\phi\in \mathcal{A}_{cusp}([G])$$

\noindent where $\lV .\rV_{L^2([G])}$ denotes the $L^2$-norm on $L^2([G])$, induce on $\mathcal{A}_{cusp}([G])^K$ its LF topology (this follows essentially from strong approximation together with the Sobolev lemma). We will denote $L^2_{cusp}([G])$ the completion of $\mathcal{A}_{cusp}([G])$ in $L^2([G])$. It is a unitary representation of $G(\A)$ which decomposes discretely.

\vspace{2mm}

\noindent Let now $f\in \mathcal{S}(G(\A))$ be a Schwartz function on $G(\A)$. We denote as usual by

$$\displaystyle K_f(x,y)=\sum_{\gamma\in G(F)} f(x^{-1}\gamma y),\;\;\; x,y\in [G]$$

\noindent the automorphic kernel of $f$. Note that the sum is absolutely convergent by Proposition \ref{prop norms}. Let $\pi\subset \mathcal{A}_{cusp}([G])$ be a cuspidal automorphic representation and let $\mathcal{B}_\pi$ be an orthonormal basis of (the completion of) $\pi$ for the $L^2$ scalar product. We define

$$\displaystyle K_{f,\pi}(x,y)=\sum_{\phi\in \mathcal{B}_\pi} (R(f)\phi)(x)\overline{\phi(y)},\;\;\; x,y\in [G]$$

\noindent Then $K_{f,\pi}$ is the orthogonal projection of $K_f$, seen as a function in $x$, onto $\pi$ or, what amounts to the same, the orthogonal projection of $K_f$, seen as a function of $y$, onto $\overline{\pi}$. Finally, letting $\mathcal{B}\subset \mathcal{A}_{cusp}([G])$ be an orthonormal basis of $L^2_{cusp}([G])$, we set

$$\displaystyle K_{f,cusp}(x,y)=\sum_{\phi\in \mathcal{B}} (R(f)\phi)(x)\overline{\phi(y)},\;\;\; x,y\in [G]$$

\noindent Note that

$$\displaystyle K_{f,cusp}=\sum_{\pi} K_{f,\pi}$$

\noindent where the sum is over a complete family of orthogonal cuspidal automorphic representations $\pi\subset \mathcal{A}_{cusp}([G])$ (all of them if there is multiplicity one).

\begin{prop}\label{prop period}
Let $H_1,H_2\subset G$ be closed algebraic subgroups such that the quotients $G/H_1$ and $G/H_2$ are quasi-affine. Then the integral

$$\displaystyle \int_{[H_1]}\int_{[H_2]}\sum_\pi \lv K_{f,\pi}(h_1,h_2)\rv dh_1dh_2$$

\noindent the sum running over a complete family of orthogonal cuspidal automorphic representations, converges. We even have the stronger following result: let $K_0$ be a compact-open subgroup of $G(\A_f)$ such that $f$ is right $K_0$-invariant and let $\mathcal{B}\subset \mathcal{A}_{cusp}([G])^{K_0}$ be an orthonormal basis of $L^2_{cusp}([G])^{K_0}$ consisting of functions which are $C_K$ and $C_G$ eigenvectors, then the integral

$$\displaystyle \int_{[H_1]}\int_{[H_2]}\sum_{\phi\in \mathcal{B}} \lv (R(f)\phi)(h_1)\rv \lv \phi(h_2)\rv dh_1dh_2$$

\noindent converges.
\end{prop}

\noindent\ul{Proof}: The second statement is obviously stronger than the first since for every cuspidal automorphic representation $\pi$ we can find an orthonormal basis of $\pi^{K_0}$ consisting of $C_K-$ and $C_G-$eigenvectors. Let $\mathcal{B}$ be an orthonormal basis of $L^2_{cusp}([G])^{K_0}$ as in the proposition. By Proposition \ref{prop norms}(ix) it suffices to prove that for all $d>0$ we have

$$\displaystyle \sum_{\phi\in \mathcal{B}} \lv (R(f)\phi)(x)\rv \lv \phi(y)\rv\ll \lV x\rV_{[G]}^{-d}\lV y\rV_{[G]}^{-d} \leqno (7)$$

\noindent for all $x,y\in [G]$. Let $d>0$. Since the family of norms $(\lV .\rV_k)_k$ generates the topology on $\mathcal{A}_{cusp}([G])^{K_0}$, there exists $k>0$ such that

$$\displaystyle \lv \phi(x)\rv\ll \lV \phi\rV_k \lV x\rV_{[G]}^{-d}$$

\noindent for all $\phi\in \mathcal{A}_{cusp}([G])^{K_0}$ and all $x\in [G]$.

\vspace{2mm}

\noindent For all $\phi\in\mathcal{B}$, let us denote by $\lambda_K(\phi),\lambda_G(\phi)\in \R$ the eigenvalues of $C_K$ and $C_G$ acting on $\phi$. Let $N$ be a positive integer. For all $\phi\in\mathcal{B}$, we have

$$\displaystyle R(f)\phi=(1+\lambda_G(\phi)^2+\lambda_K(\phi)^2)^{-N}R(f^{(N)})\phi$$

\noindent where $f^{(N)}=(1+\Delta)^N f$. Hence, we have

\[\begin{aligned}
\displaystyle \sum_{\phi\in\mathcal{B}} \lv R(f)\phi(x)\rv \lv \phi(y)\rv & =\sum_{\phi\in\mathcal{B}} (1+\lambda_G(\phi)^2+\lambda_K(\phi)^2)^{-N} \lv R(f^{(N)})\phi(x)\rv \lv \phi(y)\rv \\
 & \ll \lV x\rV_{[G]}^{-d} \lV y\rV_{[G]}^{-d}\sum_{\phi\in\mathcal{B}} (1+\lambda_G(\phi)^2+\lambda_K(\phi)^2)^{-N} \lV R(f^{(N)})\phi\rV_k \lV \phi\rV_k \\
 & =\lV x\rV_{[G]}^{-d} \lV y\rV_{[G]}^{-d}\sum_{\phi\in\mathcal{B}} (1+\lambda_G(\phi)^2+\lambda_K(\phi)^2)^{k-N} \lV R(f^{(N+k)})\phi\rV_{L^2} \lV \phi\rV_{L^2} \\
 & \leqslant \lV x\rV_{[G]}^{-d} \lV y\rV_{[G]}^{-d}\lV f^{(N+k)}\rV_{L^1}\sum_{\phi\in\mathcal{B}} (1+\lambda_G(\phi)^2+\lambda_K(\phi)^2)^{k-N}
\end{aligned}\]

\noindent for all $x,y\in [G]$ and where $\lV .\rV_{L^1}$ denotes the $L^1$-norm on $L^1([G])$. By [Mu], for $N\gg 1$ the last sum above converges. This proves (7) and ends the proof of the proposition. $\blacksquare$

\vspace{2mm}

\begin{rem}
\begin{itemize}
\item Fix $f^\infty\in \mathcal{S}(G(\A^\infty))$. Then, the proof of the proposition actually shows that for all $d>0$ there exists a continuous seminorm $\nu_d$ on $\mathcal{S}(G(F_\infty))$ so that

$$\displaystyle \left\lvert K_{f_\infty\otimes f^\infty,cusp}(x,y)\right\rvert\leqslant \nu_d(f_\infty) \lVert x\rVert_{[G]}^{-d} \lVert y\rVert_{[G]}^{-d} \leqno (8)$$

\noindent for all $f_\infty\in \mathcal{S}(G(F_\infty))$ and all $x,y\in [G]$.

\item We can prove the first part of the proposition directly by using the Selberg {\it trick}. Indeed, it suffices to show that the series

$$\displaystyle \sum_\pi K_{f,\pi}$$

\noindent converges absolutely in $\mathcal{A}_{cusp}([G\times G])$ or, what amounts to the same, that it converges absolutely in $\mathcal{A}([G\times G])$. To prove this, we only need to show that the sum

$$\displaystyle \sum_\pi \lv K_{f,\pi}(x,y)\rv$$

\noindent converges absolutely for all $x,y\in [G]$ and is bounded uniformly in $x$, $y$. By a theorem of Dixmier-Malliavin ([DM]), we may write $f$ as a finite sum of convolutions $f_{1,i}\star f_{2,i}$, $f_{1,i},f_{2,i}\in \mathcal{S}(G(\A))$, $i=1,\ldots,k$. By the Cauchy-Schwarz inequality, we have

$$\displaystyle \lv K_{f,\pi}(x,y)\rv\leqslant \sum_{i=1}^k K_{h_{2,i},\pi}(x,x)^{1/2}K_{h_{1,i},\pi}(y,y)^{1/2}$$

\noindent for all $\pi$, all $x,y\in [G]$ and where we have set $h_{j,i}=f_{j,i}^*\star f_{j,i}$ where by definition $f_{j,i}^*(g)=\overline{f_{j,i}(g^{-1})}$. Thus, by another application of Cauchy-Schwarz, we get

$$\displaystyle \sum_\pi \lv K_{f,\pi}(x,y)\rv\leqslant \sum_{i=1}^k K_{h_{2,i},cusp}(x,x)^{1/2}K_{h_{1,i},cusp}(y,y)^{1/2}$$

\noindent for all $x,y\in [G]$ and the right hand side is uniformly bounded (even of rapid decay).
\end{itemize}
\end{rem}

\subsection{Relative trace formulae}\label{RTF}

\noindent We now return to the situation considered in section \ref{section 2}. We will use the same notations and normalization of measures as there (in particular our global Haar measures are Tamagawa measures). We define the following 'bases' (geometric quotients): $\mathcal{B}:=H\backslash G/H$ and $\mathcal{B}':=H'_1\backslash G'/H'_2$. These are affine varieties (as $H$, $H'_1$ and $H'_2$ are reductive). We set $p:G\to \mathcal{B}$ and $p':G'\to \mathcal{B}'$ for the natural projections and $p_F$, $p'_F$ for the corresponding maps at the level of $F$-points. For all test function $f\in \mathcal{S}(G(\A))$, $f'\in \mathcal{S}(G'(\A))$, $\delta\in \mathcal{B}(F)$, $\gamma\in \mathcal{B}'(F)$, we set

$$\displaystyle K_{f,\delta}(x,y):=\sum_{\underline{\delta}\in p_F^{-1}(\delta)} f(x^{-1}\underline{\delta}y),\;\;\; x,y\in G(\mathbb{A})$$

$$\displaystyle K_{f',\gamma}(x,y):=\sum_{\underline{\gamma}\in (p'_F)^{-1}(\gamma)} f'(x^{-1}\underline{\gamma}y),\;\;\; x,y\in G'(\mathbb{A})$$

\noindent Note that these sums are absolutely convergent by Proposition \ref{prop norms} (v) and that

$$\displaystyle K_f(x,y)=\sum_{\delta\in \mathcal{B}(F)} K_{f,\delta}(x,y)$$

$$\displaystyle K_{f'}(x,y)=\sum_{\gamma\in \mathcal{B}'(F)} K_{f',\gamma}(x,y)$$

\noindent Whenever convergent, we define the following 'global orbital integrals'

$$\displaystyle O(\delta,f):=\int_{[H]}\int_{[H]} K_{f,\delta}(h_1, h_2)dh_1dh_2$$

$$\displaystyle O(\gamma,f'):=\int_{[H_1']}\int_{[H_2']} K_{f',\gamma}(h_1,h_2)\eta(h_2)dh_2dh_1$$

\noindent and the following two expressions

$$\displaystyle J(f)=\int_{[H]}\int_{[H]} K_f(h_1,h_2)dh_1dh_2$$

$$\displaystyle I(f')=\int_{[H_1']}\int_{[H_2']} K_{f'}(h_1,h_2)\eta(h_2)dh_2dh_1$$

\begin{prop}\label{prop RTF}
\begin{enumerate}[(i)]
\item Assume that $f\in \mathcal{S}(G(\A))$ is a nice function (see \S \ref{transfer of spherical characters}). Then the expressions defining $J(f)$ and $O(\delta,f)$, $\delta\in \mathcal{B}(F)$, are absolutely convergent and we have the equalities

$$\displaystyle \sum_{\delta\in \mathcal{B}(F)} O(\delta,f)=J(f)=\sum_{\pi} J_\pi(f)$$

\noindent where the left sum is absolutely convergent and the right sum is over the set of cuspidal automorphic representations $\pi$ of $G(\A)$.

\item Assume that $f'\in \mathcal{S}(G'(\A))$ is a nice function (see \S \ref{transfer of spherical characters}). Then the expressions defining $I(f')$ and $O(\gamma,f')$, $\gamma\in \mathcal{B}'(F)$, are absolutely convergent and we have the equalities

$$\displaystyle \sum_{\gamma\in \mathcal{B}'(F)} O(\gamma,f')=I(f')=\sum_{\Pi} 2^{-2}L(1,\eta_{E/F})^{-2}I_\Pi(f')$$

\noindent where the left sum is absolutely convergent and the right sum is over the set of cuspidal automorphic representations $\Pi$ of $G'(\A)$ whose central character is trivial on $Z_{H_2'}(\A)$.
\end{enumerate}
\end{prop}

\vspace{2mm}

\noindent\ul{Proof}: We only prove (ii) the proof of (i) being similar.

\vspace{2mm}

\noindent Set $\widetilde{G}'=G'/Z_{H'_2}$ and define $\widetilde{f}'\in \mathcal{S}(\widetilde{G}'(\A))$ by

$$\displaystyle \widetilde{f}'(\widetilde{g})=\int_{Z_{H_2}'(\A)}f'(z\widetilde{g}) dz$$

\noindent Then we have, at least formally,

$$\displaystyle I(f')=\int_{[H_1']}\int_{[H_2'/Z_{H_2'}]}K_{\widetilde{f}'}(h_1,h_2)\eta(h_2) dh_2dh_1$$

\noindent As $f'$ is a nice function we have $K_{\widetilde{f}'}=K_{\widetilde{f}',cusp}$. Thus, by Proposition \ref{prop period}, it follows that the expression defining $I(f')$ is absolutely convergent and that

$$\displaystyle I(f')=\sum_\Pi \int_{[H_1']}\int_{[H_2'/Z_{H_2'}]}K_{\widetilde{f}',\Pi}(h_1,h_2)\eta(h_2) dh_2dh_1$$

\noindent where the sum is over the set of all cuspidal automorphic representations $\Pi$ of $G'(\A)$ with a central character trivial on $Z_{H_2'}(\A)$. We would like to identify the term indexed by $\Pi$ above with the global spherical character $I_\Pi(f')$. However, we don't have equality on the nose because the scalar products used to define $I_\Pi(f')$ and $K_{\widetilde{f}',\Pi}$ are not the same. More precisely, $I_\Pi(f')$ is defined using the Petersson scalar product $(.,.)_{Pet}$ of section \ref{section 2} whereas in the definition of $K_{\widetilde{f}',\Pi}$ we have used the scalar product

\[\begin{aligned}
\displaystyle (\phi,\phi')_{L^2([\widetilde{G}'])} & =\int_{[\widetilde{G}']} \phi(\widetilde{g})\overline{\phi'(\widetilde{g})} d\widetilde{g} \\
 & =vol\left(Z_{G'}(F)Z_{H_2'}(\A)\backslash Z_{G'}(\A)\right) (\phi,\phi')_{Pet}
\end{aligned}\]

\noindent Thus, we get

$$\displaystyle \int_{[H_1']}\int_{[H_2'/Z_{H_2'}]}K_{\widetilde{f}',\Pi}(h_1,h_2)\eta(h_2) dh_2dh_1=vol\left(Z_{G'}(F)Z_{H_2'}(\A)\backslash Z_{G'}(\A)\right)^{-1} I_\Pi(f')$$

\noindent By \ref{measures}(1), we have

$$\displaystyle vol\left(Z_{G'}(F)Z_{H_2'}(\A)\backslash Z_{G'}(\A)\right)=vol\left(E^\times \A^\times\backslash \A_E^\times \right)^2=2^2L(1,\eta_{E/F})^2$$

\noindent and the second equality of (ii) follows.

\vspace{2mm}

\noindent The first equality follows from standard formal manipulations. To justify these manipulations, we need to establish that the following expression is absolutely convergent (as a triple integral)

$$\displaystyle \int_{[H_1']}\int_{[H_2']} \sum_{\gamma\in \mathcal{B}'(F)}K_{f',\gamma}(h_1,h_2)\eta(h_2)dh_2dh_1=\int_{[H_1']}\int_{[H_2'/Z_{H_2'}]} \sum_{\gamma\in \mathcal{B}'(F)}K_{\widetilde{f}',\gamma}(h_1,h_2)\eta(h_2)dh_2dh_1$$

\noindent where for all $\gamma\in \mathcal{B}'(F)$, $K_{\widetilde{f}',\gamma}$ is defined the same way as $K_{f',\gamma}$. For this, by Proposition \ref{prop norms} (v) and (ix), it suffices to show that for all $d>0$ we have an inequality

$$\displaystyle \left\lvert K_{\widetilde{f}',\gamma}(h_1,h_2)\right\rvert\ll \lVert \gamma\rVert_{\mathcal{B}'}^{-d} \lVert h_1\rVert_{[\widetilde{G}']}^{-d} \lVert h_2\rVert_{[\widetilde{G}']}^{-d} \leqno (1)$$

\noindent for all $\gamma\in \mathcal{B}'(F)$ and all $h_1\in [H'_1]$, $h_2\in [H'_2/Z_{H_2'}]$. Since $f'$ is factorizable, we may write $f'=f'_\infty\otimes {f'}^\infty$ with $f'_\infty\in \mathcal{S}(G'(F_\infty))$ and ${f'}^\infty\in \mathcal{S}(G'(\mathbb{A}_F^\infty))$. Let $\mathcal{B}'\hookrightarrow V=\mathbb{A}^r$ be a closed embedding of $\mathcal{B}'$ into some affine space. For all $\varphi\in C_c^\infty(V(F_\infty))$, we define $f'_\varphi:=(\varphi f'_\infty)\otimes {f'}^\infty$ where we identify $\varphi$ with a function on $G'(F_\infty)$ by composition with the projection $G'(F_\infty)\to \mathcal{B}'(F_\infty)$ and the embedding $\mathcal{B}'(F_\infty)\hookrightarrow V(F_\infty)$. Then, $f'_\varphi$ is again a nice function so that $K_{\widetilde{f}'_\varphi}=K_{\widetilde{f}'_\varphi,cusp}$ and, by \ref{norms}(8), for all $d>0$ there exists a continuous semi-norm $\nu_d$ on $\mathcal{S}(G'(F_\infty))$ such that

$$\displaystyle \left\lvert K_{\widetilde{f}'_\varphi}(x,y)\right\rvert\leqslant \nu_d(\varphi f'_\infty)\lVert x\rVert_{[\widetilde{G}']}^{-d} \lVert y\rVert_{[\widetilde{G}']}^{-d} \leqno (2)$$

\noindent for all $\varphi \in C_c^\infty(V(F_\infty))$ and all $x,y\in [G]$. Let $\Gamma$ denote the intersection of $\mathcal{B}'(F)$ with the projection of the support of ${f'}^\infty$ (a compact subset of $\mathcal{B}'(\mathbb{A}_F^\infty)$). Note that for all $\varphi \in C_c^\infty(V(F_\infty))$, all $\gamma\in \mathcal{B}'(F)$ and all $(h_1,h_2)\in H_1'(\mathbb{A}_F)\times H_2'(\mathbb{A}_F)$, we have $K_{\widetilde{f}'_\varphi,\gamma}(h_1,h_2)=\varphi(\gamma)K_{\widetilde{f}',\gamma}(h_1,h_2)$ and $K_{\widetilde{f}',\gamma}(h_1,h_2)=0$ if $\gamma\notin \Gamma$. Hence, by (2), to show (1) it suffices to construct a family $(\varphi_\gamma)_{\gamma\in \Gamma}$ of functions in $C_c^\infty(V(F_\infty))$ satisfying the following two conditions

\begin{itemize}
\item For all $\gamma\in \Gamma$, we have $Supp(\varphi_\gamma)\cap \Gamma=\{\gamma \}$ and $\varphi_\gamma(\gamma)=1$;

\item For all $d>0$, the function $\gamma\in \Gamma\mapsto \lVert \gamma\rVert_{\mathcal{B}'}^d\nu_d(\varphi_\gamma f'_\infty)$ is bounded.
\end{itemize}

\noindent There exists a lattice $L\subset V(F_\infty)$ containing $\Gamma$. Fix a function $\varphi_0\in C_c^\infty(V(F_\infty))$ with the property that $Supp(\varphi_0)\cap L=\{0 \}$ and $\varphi_0(0)=1$. For all $v\in L$, define $\varphi_v\in C_c^\infty(V(F_\infty))$ by $\varphi_v(x)=\varphi_0(x-v)$. We claim that the family $(\varphi_\gamma)_{\gamma\in \Gamma}$ satisfies the two conditions above. Indeed, the first condition is clear and for all $d>0$, there exist $k>0$ and two finite families $(u_i)_{i\in I}$ and $(v_i)_{i\in I}$ of elements of $\mathcal{U}(\mathfrak{g}'_\infty)$ such that

$$\displaystyle \nu_d(\varphi f'_\infty)\leqslant \sup_{i\in I}\sup_{g\in G'(F_\infty)} \left\lvert (R(u_i)\varphi)(g) (R(v_i)f'_\infty)(g)\right\rvert \lVert g\rVert_{G'}^k$$

\noindent for all $\varphi\in C_c^\infty(V(F_\infty))$. Then, for all $\ell>0$ we have

$$\displaystyle \nu_d(\varphi f'_\infty)\leqslant \left(\sup_{i\in I}\sup_{g\in G'(F_\infty)} \left\lvert (R(u_i)\varphi)(g)\right\rvert \lVert g\rVert_{G'}^{-\ell} \right) \left(\sup_{i\in I}\sup_{g\in G'(F_\infty)} \left\lvert (R(v_i)f'_\infty)(g)\right\rvert \lVert g\rVert_{G'}^{k+\ell} \right)$$

\noindent for all $\varphi\in C_c^\infty(V(F_\infty))$. Thus, it suffices to show that for all $u\in \mathcal{U}(\mathfrak{g}'_\infty)$ there exists $\ell>0$ so that the function

$$\displaystyle \gamma\in \Gamma\mapsto \lVert \gamma\rVert_{\mathcal{B}'}^d\sup_{g\in G'(F_\infty)} \left\lvert (R(u)\varphi_\gamma)(g)\right\rvert \lVert g\rVert_{G'}^{-\ell}$$

\noindent is bounded. Fix $u\in \mathcal{U}(\mathfrak{g}'_\infty)$. Then, we can find a finite family $(r_j)_{j\in J}$ of regular functions in $\mathbb{C}\left[R_{F/\mathbb{Q}}G'  \right]$ and a finite family $(X_j)_{j\in J}$ of elements of the symmetric algebra of $V(F_\infty)$ such that, denoting by $\partial(X_j)$ the corresponding constant coefficients differential operators on $V(F_\infty)$, we have

$$\displaystyle R(u)\varphi=\sum_{j\in J} r_j \partial(X_j)\varphi$$

\noindent for all $\varphi\in C_c^\infty(V(F_\infty))$. Since for each $j\in J$ the absolute value of $r_j$ is bounded by a constant times a power of the norm $\lVert .\rVert_{G'}$, we are reduced to prove that for all $X$ in the symmetric algebra of $V(F_\infty)$ there exists $\ell>0$ such that

$$\displaystyle \lVert \gamma\rVert_{\mathcal{B}'}^d\sup_{g\in G'(F_\infty)} \left\lvert \partial(X)\varphi_\gamma(g)\right\rvert \lVert g\rVert_{G'}^{-\ell}\ll 1$$

\noindent for all $\gamma\in \Gamma$. By our choice of the functions $(\varphi_\gamma)_{\gamma\in \Gamma}$ and Proposition \ref{prop norms}(ii) we have

$$\displaystyle \lVert \gamma\rVert_{\mathcal{B}'}\prec\lVert g\rVert_{G'}$$

\noindent for all $\gamma\in \Gamma$ and all $g\in G'(F_\infty)$ with $\varphi_\gamma(g)\neq 0$. Hence, we just need to show that for all $X$ in the symmetric algebra of $V(F_\infty)$ the function

$$\displaystyle \gamma\in \Gamma\mapsto \sup_{g\in G'(F_\infty)} \left\lvert \partial(X)\varphi_\gamma(g)\right\rvert=\sup_{v\in \mathcal{B}'(F_\infty)} \left\lvert \partial(X)\varphi_\gamma(v)\right\rvert$$

\noindent is bounded. But this is obvious by the way we have defined the functions $(\varphi_\gamma)_{\gamma\in \Gamma}$. $\blacksquare$
 
\subsection{Proof of Theorem \ref{comparison RTF}}\label{proof of comparison RTF}

\noindent Let $f\in \mathcal{S}(G(\A))$ and $f'\in \mathcal{S}(G'(\A))$ be nice functions and assume that there exists a tuple $(f^{W'})_{W'}$, $f^{W'}\in \mathcal{S}(G^{W'}(\A))$, matching $f'$ and such that $f^W=f$. By Theorem \ref{fundamental lemma}, we may assume that $f^{W'}=0$ for almost all $W'$. There are natural isomorphisms $H^{W'}\backslash G^{W'}/H^{W'}\simeq \mathcal{B}\simeq \mathcal{B}'$ for all $W'$ (see \S \ref{correspondence}). In order to compare the trace formulas of Proposition \ref{prop RTF}, we need to know that for all $\delta\in \mathcal{B}(F)$ and all $\gamma\in \mathcal{B}'(F)$ corresponding to each other via the previous bijection we have

$$\displaystyle \sum_{W'}O(\delta,f^{W'})= O(\gamma,f') \leqno{(1)}$$

\noindent Note that for $\gamma$ regular semi-simple (i.e. such that the fiber over $\gamma$ in $G'(F)$ consists of regular semi-simple elements), this is a direct consequence of the fact that $(f^{W'})_{W'}$ matches $f'$ (in this case there is at most one nonzero contribution in the left sum). To treat the general case, we need to use recent results of Zydor and Chaudouard-Zydor. More precisely, in [Zy] Zydor has defined, for $\delta\in \mathcal{B}(F)$ and $\gamma\in \mathcal{B}'(F)$, certain distributions

$$\displaystyle g^{W'}\in \mathcal{S}(G^{W'}(\A))\mapsto O^Z(\delta,g^{W'})$$

$$\displaystyle g'\in \mathcal{S}(G'(\A))\mapsto O^Z(\gamma,g')$$

\noindent The definition of $O^Z(\gamma,.)$ is (roughly) as follows. let $A$ be the standard maximal split torus in $GL_n$ and set $\mathfrak{a}:=X_*(A)\otimes \mathbb{R}$. Then for $T\in \mathfrak{a}$, $g'\in \mathcal{S}(G'(\mathbb{A}))$ and $\gamma\in \mathcal{B}'(F)$, Zydor defines a certain 'truncated' kernel $K_{g',\gamma}^T$ on $[H_1']\times [H_2']$ (see [Zy] \S 5.5, note that the $\mathfrak{o}$ of {\it loc. cit.} corresponds to our $\gamma$ and that the function $f$ of {\it loc. cit.} corresponds not to $g'$ but rather to its descent $\widetilde{g}'$ to $S_{n+1}(\mathbb{A})$ as in \S \ref{orbital integrals}) and he shows that for $T$ in a certain cone the integral

$$\displaystyle O^{Z,T}(\gamma,g'):=\int_{[H_1']\times [H_2']} K_{g',\gamma}^T(h_1,h_2)\eta(h_2) dh_2dh_1$$

\noindent converges absolutely (see [Zy,Theorem 5.9]). The definition of $K_{g',\gamma}^T$ is as a sum the main term being $K_{g',\gamma}$ and the remaining terms depending only on $g'_P$ for certain proper parabolic subgroup $P=MU$ of $G'$ where $\displaystyle g'_P(x):=\int_{U(\A)}g'(xu)du$. Since $f'$ is a nice function, we have $f'_P=0$ for all proper parabolic subgroup and thus $K_{f',\gamma}^T=K_{f',\gamma}$ and $O^{Z,T}(\gamma,f')=O(\gamma,f')$ for all $T$ (remark that this also reproves the absolute convergence of $O(\gamma,f')$ of Proposition \ref{prop RTF}). Finally, still for $T$ in a certain cone, Zydor shows that the function $T\mapsto O^{Z,T}(\gamma,g')$ is an exponential-polynomial whose purely polynomial term is constant ([Zy, Theorem 5.9]) and he defines $O^Z(\gamma,g')$ to be this constant. Since $O^{Z,T}(\gamma,f')=O(\gamma,f')$ for all $T$, we also have $O^{Z}(\gamma,f')=O(\gamma,f')$. The definition of $O^Z(\delta,.)$ is similar and since the functions $f^{W'}$ are nice the same argument shows that $O^{Z}(\delta,f^{W'})=O(\delta,f^{W'})$ for all $W'$ and all $\delta\in \mathcal{B}(F)$. Finally, the main result of [CZ] is that if $g'\in \mathcal{S}(G'(\A))$ match a tuple of functions $(g^{W'})_{W'}$, $g^{W'}\in \mathcal{S}(G^{W'}(\A))$, then we have

$$\displaystyle \sum_{W'}O^Z(\delta,g^{W'})= O^Z(\gamma,g')$$

\noindent for all $\delta\in \mathcal{B}(F)$ and all $\gamma\in \mathcal{B}'(F)$ corresponding to each other (strictly speaking [CZ] only considers compactly supported functions but the proofs applies verbatim to Schwartz functions). Together with the previous equalities this shows (1). 

\vspace{2mm}

\noindent Now (1) together with Proposition \ref{prop RTF} leads to the identity

$$\displaystyle \sum_{W'}\sum_{\pi_{W'}} J_{\pi_{W'}}(f^{W'})=\sum_{\Pi} 2^{-2}L(1,\eta_{E/F})^{-2} I_\Pi(f')\leqno (2)$$

\noindent where $\pi_{W'}$ runs over the set of all cuspidal automorphic representations of $G^{W'}(\A)$ and $\Pi$ runs over the set of all cuspidal automorphic representations of $G'(\A)$ whose central character is trivial on $Z_{H_2'}(\A)$. Fix a maximal compact subgroup $K^{W'}=\prod_v K^{W'}_v$ of $G^{W'}(\A)$ for all $W'$ and let $\Sigma$ be the infinite set of places $v$ of $F$ which split in $E$ and where $\pi$, $f$ and $f'$ are unramified. By Theorem \ref{fundamental lemma}, we may assume that for all $W'$ and all $v\in \Sigma$ the function $f^{W'}_v$ is unramified (i.e. it equals $vol(K^{W'}_v)^{-1}\mathbf{1}_{K^{W'}_v}$). Then, in the equality (2) only the $\pi_{W'}$ and the $\Pi$ which are unramified at all places in $\Sigma$ contribute. Define the Hecke algebra $\mathcal{H}_{G,\Sigma}=C_c(G(\A_\Sigma)//K_\Sigma)$ of compactly supported and $K_\Sigma$-biinvariant functions on $G(\A_\Sigma)$. This is the restricted tensor product over $v\in \Sigma$ of the local Hecke algebras $\mathcal{H}_{G,v}=C_c(G(F_v)//K_v)$. We define similarly the Hecke algebra $\mathcal{H}_{G',\Sigma}=C_c(G'(\A_\Sigma)//K'_\Sigma)$ and the local Hecke algebras $\mathcal{H}_{G',v}=C_c(G'(F_v)//K'_v)$. Note that for all $n$-dimensional hermitian space $W'$ over $E$ we have an isomorphism $G(\A_\Sigma)\simeq G^{W'}(\A_\Sigma)$ canonical up to conjugation which induces a canonical isomorphism $\mathcal{H}_{G,\Sigma}\simeq C_c(G^{W'}(\A_\Sigma)//K^{W'}_\Sigma)$. There is a base change homomorphism $\mathcal{H}_{G',\Sigma}\to \mathcal{H}_{G,\Sigma}$, $h\mapsto h^{bc}$ and for all $v\in \Sigma$, all $W'$ and all $h_v\in \mathcal{H}_{G',v}$, $h_v\star f'_v$ and $h^{bc}_v\star f^{W'}_v$ match each other (see [Zh1, Proposition 2.5]). For all irreducible unitary representation $\Pi$ of $G'(\A)$ which is unramified at all places in $\Sigma$ let us denote by $h\mapsto \widehat{h}(\Pi)$ the corresponding character of the Hecke algebra $\mathcal{H}_{G',\Sigma}$. Then, for all $W'$, all cuspidal automorphic representation $\pi_{W'}$ which is unramified at all places in $\Sigma$ and all $h\in \mathcal{H}_{G',\Sigma}$ the element $h^{bc}\in \mathcal{H}_{G,\Sigma}$ acts on $\pi_{W'}^{K^{W'}_\Sigma}$ by $\widehat{h}(BC(\pi_{W'}))$. Let $h\in \mathcal{H}_{G',\Sigma}$. Since the functions $h\star f'$ and $(h^{bc}\star f^{W'})_{W'}$ are nice and match each other we can apply equality (2) to these functions to get 

$$\displaystyle \sum_{W'}\sum_{\pi_{W'}} \widehat{h}(BC(\pi_{W'}))J_{\pi_{W'}}(f^{W'})=\sum_{\Pi} 2^{-2}L(1,\eta_{E/F})^{-2} \widehat{h}(\Pi) I_\Pi(f') \leqno (3)$$

\noindent Let $Irr_{unit,\Sigma}(G'(\A))$ be the set of all irreducible unitary representations of $G'(\A)$ which are unramified at all places in $\Sigma$. Then the functions $\Pi\in Irr_{unit,\Sigma}(G'(\A))\mapsto \widehat{h}(\Pi)$, $h\in \mathcal{H}_{G',\Sigma}$, are bounded and we have $\widehat{h^*}=\overline{\widehat{h}}$ where $h^*(g)=\overline{h(g^{-1})}$. Hence by the Stone-Weierstrass theorem, from (3) we deduce

$$\displaystyle \sum_{W'}\sum_{\pi_{W'}} J_{\pi_{W'}}(f^{W'})=\sum_{\Pi} 2^{-2}L(1,\eta_{E/F})^{-2} I_\Pi(f') \leqno (4)$$

\noindent where this time $\pi_{W'}$ and $\Pi$ run over the sets of cuspidal automorphic representations of $G^{W'}(\A)$ and $G'(\A)$ such that $BC(\pi_{W',v})=\Pi_v=BC(\pi_v)$ for all $v\in \Sigma$. Recall the following {\it automorphic-Cebotarev-density} theorem due to Ramakrishnan ([Ra]):

\begin{theo}[Ramakrishnan]\label{theo ram}
Let $\Pi_1$, $\Pi_2$ be two isobaric automorphic representations of $GL_d(\A_E)$ such that $\Pi_{1,v}\simeq \Pi_{2,v}$ for almost all places $v$ of $F$ that are split in $E$. Then, $\Pi_1=\Pi_2$.
\end{theo}

\noindent As $BC(\pi_{W'})$ is always isobaric it follows from this theorem that the right hand side of (4) reduces to $2^{-2}L(1,\eta_{E/F})^{-2} I_{BC(\pi)}(f')$ and that if $\pi_{W'}$ contributes to the left hand side then $BC(\pi_{W'})=BC(\pi)$. In particular, $\pi_{W'}$ and $\pi$ belong to the same (global) Vogan $L$-packet. By the local Gan-Gross-Prasad conjecture (see \S \ref{local GGP}), and since by assumption $\pi$ is tempered at all archimedean places, we know that there is at most one abstractly $H^{W'}$-distinguished representation in this $L$-packet. By assumption, $\pi$ is such a representation. Hence, the left hand side of (4) reduces to $J_\pi(f)$ and this ends the proof of Theorem \ref{comparison RTF}. $\blacksquare$

\bigskip

\textbf{Bibliography}

\bigskip

[AMS] A.-M. Aubert, A. Moussaoui, M. Solleveld: {\it Generalizations of the Springer correspondence and cuspidal Langlands parameters}, prepublication 2015

[AGRS] A. Aizenbud, D. Gourevitch, S. Rallis, G. Schiffmann: \textit{Multiplicity one theorems}, Ann. of Math. (2) 172 (2010), no. 2, 1407-1434

[BC] N. Bergeron, L. Clozel: {\it Spectre automorphe des vari\'et\'es hyperboliques et applications topologiques}, Ast\'erisque No. 303 (2005), xx+218 pp

[BD] J. Bernstein, P. Deligne: {\it Le "centre" de Bernstein}, In "Representations des groupes r\'eductifs sur un corps local, Travaux en cours" (P.Deligne ed.), Hermann, Paris, 1-32 (1984)

[BK] J. N. Bernstein and B. Kr\"otz: {\it Smooth Fr\'echet globalizations of Harish-Chandra modules}, Israel J. Math. 199 (2014), no. 1, 45-111

[Beu1] R. Beuzart-Plessis: {\it A local trace formula for the Gan-Gross-Prasad conjecture for unitary groups: the archimedean case}, prepublication 2015

[Beu2]----------------------: {\it Endoscopie et conjecture locale raffin\'ee de Gan-Gross-Prasad pour les groupes unitaires}, Compos. Math. 151 (2015), no. 7, 1309-1371

[BHC] A. Borel, Harish-Chandra: {\it Arithmetic subgroups of algebraic groups}, Ann. of Math. 75 (1962), 485-535

[Ca] W. Casselman: {\it Canonical extensions of Harish-Chandra modules to representations of $G$}, Canad. J. Math. 41 (1989), no. 3, 385-438

[CZ] P.-H. Chaudouard, M. Zydor: {\it Le transfert singulier pour la formule des traces de Jacquet-Rallis}, preprint 2016, arXiv:1611.09656

[Cl1] L. Clozel: {\it Orbital integrals on p-adic groups: a proof of the Howe conjecture}, Ann. of Math. (2) 129 (1989), no. 2, 237-251

[Cl2]--------------: {\it Characters of nonconnected, reductive p-adic groups}, Canad. J. Math. 39 (1987), no. 1, 149-167

[De] P. Deligne: {\it Le support du caract\`ere d'une repr\'esentation supercuspidale}, C. R. Acad. Sci. Paris S\'er. A-B 283 (1976), no. 4, Aii, A155-A157

[DM] J. Dixmier, P. Malliavin: {\it Factorisations de fonctions et de vecteurs ind\'efiniment diff\'erentiables}, Bull. Sci. Math. (2) 102 (1978), no. 4, 307-330

[GGP] W. T. Gan, B. H. Gross, D. Prasad: \textit{Symplectic local root numbers, central critical $L$-values and restriction problems in the representation theory of classical groups}, in "Sur les conjectures de Gross et Prasad. I" Ast\'erisque No. 346 (2012), 1-109

[GI] W. T. Gan, A. Ichino: {\it The Gross-Prasad conjecture and local theta correspondence}, prepublication 2014

[Ha] R. Neal Harris: {\it A refined Gross-Prasad conjecture for unitary groups}, arXiv:1201.0518

[HT] M. Harris, R. Taylor: {\it The geometry and cohomology of some simple Shimura varieties}, volume 151 of Annals of Mathematics Studies. Princeton University Press, Princeton, NJ, 2001. With an appendix by Vladimir G. Berkovich

[He] G. Henniart: {\it Une preuve simple des conjectures de Langlands pour $GL(n)$ sur un corps $p$-adique}, Invent. Math., 139(2), 439-455, 2000

[II] A. Ichino, T. Ikeda: {\it On the periods of automorphic forms on special orthogonal groups and the Gross-Prasad conjecture}, Geom. Funct. Anal. 19 (2010), no. 5, 1378-1425

[ILM] A. Ichino, E. Lapid, Z. Mao: {\it On the formal degrees of square-integrable representations of odd special orthogonal and metaplectic groups}, prepublication 2014

[Jac] H. Jacquet: {\it Archimedean Rankin-Selberg integrals}, Automorphic forms and L-functions II. Local aspects, 57-172, Contemp. Math., 489, Amer. Math. Soc., Providence, RI, 2009

[JPSS] H. Jacquet, I. I. Piatetski-Shapiro, J. A. Shalika: {\it Rankin- Selberg convolutions}, Amer. J. Math. 105 (1983), no. 2, 367-46

[JR] H. Jacquet, S. Rallis: {\it On the Gross-Prasad conjecture for unitary groups}, in "On certain L-functions", 205-264, Clay Math. Proc., 13, Amer. Math. Soc., Providence, RI, 2011

[JS] H. Jacquet, J. A. Shalika: {\it On Euler products and the classification of automorphic representations I}, Amer. J. Math. 103 (1981), no. 3, 499-558

[JSZ] D. Jiang, B. Sun, C.-B. Zhu: \textit{Uniqueness of Bessel models: The archimedean case}, Geom. Funct. Anal. 20 (2010), no. 3, 690-709

[KMSW] T. Kaletha, A. Minguez, S. W. Shin, P.-J. White: {\it Endoscopic Classification of Representations: Inner Forms of Unitary Groups}, prepublication 2014

[Kott] R. E. Kottwitz : {\it Harmonic analysis on reductive $p$-adic groups and Lie algebras}, in "Harmonic analysis, the trace formula, and Shimura varieties", 393522, Clay Math. Proc., 4, Amer. Math. Soc., Providence, RI, 2005

[Lan] R. P. Langlands: {\it On the classification of irreducible representations of real algebraic groups}, in "Representation theory and harmonic analysis on semisimple Lie groups", 101-170, Math. Surveys Monogr., 31, Amer. Math. Soc., Providence, RI, 1989

[LRS] W. Luo, Z. Rudnick, P. Sarnak: {\it On the generalized Ramanujan conjecture for $GL(n)$},  in "Automorphic forms, automorphic representations, and arithmetic" (Fort Worth, TX, 1996), Proc. Sympos. Pure Math., vol. 66, Amer. Math. Soc., Providence, RI, 1999, pp. 301-310

[Moe1] C. M\oe{}glin: {\it Sur la classification des s\'eries discr\`etes des groupes classiques p-adiques: param\`etres de Langlands et exhaustivit\'e}, J. Eur. Math. Soc. (JEMS) 4 (2002), no. 2, 143-200

[Moe2]:--------------------: {\it Classification et changement de base pour les s\'eries discr\`etes des groupes unitaires $p$-adiques},  Pacific J. Math. 233 (2007), no. 1, 159-204

[MT] C. M\oe{}glin, M. Tad\`ic: {\it Construction of discrete series for classical $p$-adic groups}, J. Amer. Math. Soc. 15 (2002), no. 3, 715-786

[MW1] C.M\oe{}glin, J.-L. Waldspurger: {\it Spectral decomposition and Eisenstein series. Une paraphrase de l'\'Ecriture}, Cambridge Tracts in Mathematics, 113. Cambridge University Press, Cambridge, 1995. xxviii+338 pp

[MW2]---------------------------------: {\it Le spectre r\'esiduel de $GL(n)$}, Ann. Sci. \'Ecole Norm. Sup. (4) 22 (1989), no. 4, 605-674

[Mok] C. P. Mok: {\it Endoscopic Classification of representations of Quasi-Split Unitary Groups}, Mem. Amer. Math. Soc. 235 (2015), no. 1108

[Mou] A. Moussaoui: {\it Centre de Bernstein enrichi pour les groupes classiques}, prepublication 2015

[Mui] G. Mu\`ic: {\it Some results on square integrable representations; irreducibility of standard representations}, Internat. Math. Res. Notices (1998), no. 14, 705-726

[Mu] W. M\"uller: {\it The trace class conjecture without the $K$-finiteness assumption}, C. R. Acad. Sci. Paris S\'er. I Math. 324 (1997), no. 12, 1333-1338

[MS] W. M\"uller, B. Speh: {\it Absolute convergence of the spectral side of the Arthur trace formula for $GL_n$}, Geom. Funct. Anal. 14 (2004), no. 1, 58-93, With an appendix by E. M. Lapid

[Ra] D. Ramakrishnan: {\it A theorem on GL(n) a la Tchebotarev}, Preprint 2012

[SV] Y. Sakellaridis, A. Venkatesh: {\it Periods and harmonic analysis on spherical varieties}, prepublication 2012

[Sil1] A. Silberger: {\it Introduction to harmonic analysis on reductive $p$-adic groups}, Mathematical Notes, vol. 23, Princeton University Press, Princeton, N.J., 1979, Based on lectures by Harish-Chandra at the Institute for Advanced Study, 1971-1973

[Sil2]-----------------: {\it Special representations of reductive p-adic groups are not integrable}, Ann. of Math., 111, 1980, pp. 571-587

[Tad] M. Tad\'ic: {\it Geometry of dual spaces of reductive groups (non-Archimedean case)}, J. Analyse Math. 51 (1988), 139-181

[Wald1] J.-L. Waldspurger: \textit{La formule de Plancherel pour les groupes $p$-adiques (d'apr\`es Harish-Chandra)}, J. Inst. Math. Jussieu 2 (2003), no. 2, 235-333

[Wald2]------------------: {\it La formule des traces locales tordue}, prepublication 2013

[Wald3]------------------: {\it Sur les valeurs de certaines fonctions $L$ automorphes en leur centre de sym\'etrie}, Compositio Math. 54 (1985), no. 2, 173-242

[Wall] N. Wallach: {\it Real reductive groups II}, Academic Press 1992

[Xu] B. Xu: {\it On the cuspidal support of discrete series for $p$-adic quasisplit $Sp(N)$ and $SO(N)$}, prepublication 2015

[Xue] H. Xue: {\it On the global Gan-Gross-Prasad conjecture for unitary groups: approximating smooth transfer of Jacquet-Rallis}, prepublication 2015

[Yu] Z. Yun: {\it The fundamental lemma of Jacquet and Rallis} With an appendix by Julia Gordon. Duke Math. J. 156 (2011), no. 2, 167-227

[Zh1] W. Zhang: {\it Fourier transform and the global Gan-Gross-Prasad conjecture for unitary groups}, Ann. of Math. (2) 180 (2014), no. 3, 971-1049

[Zh2]------------: {\it Automorphic period and the central value of Rankin-Selberg L-function}, J. Amer. Math. Soc. 27 (2014), no. 2, 541-612

[Zh3]------------: {\it  On arithmetic fundamental lemmas}, Invent. Math., Volume 188, Number 1 (2012), 197-252

[Zy] M. Zydor: {\it Les formules des traces relatives de Jacquet-Rallis grossi\`eres}, preprint 2015, arXiv:1510.04301

\end{document}